\documentclass[amsppt,11pt]{amsart}
\usepackage{curves}
\usepackage{amsmath}
\def\Line#1{\hbox to\textwidth{#1}}%

\newtheorem{theorem}{Theorem}[section]
\newtheorem{lemma}{Lemma}[section]
\newtheorem{definition}{Definition}[section]
\newtheorem{cor}{Corollary}
\newtheorem{prop}{Proposition}
\newcommand{\stoss}{Stosszahlansatz }
\newcommand{\x}{\times}

\renewcommand{\>}{\rangle}

\renewcommand{\b}{\beta}
\newcommand{\bs}{\bigskip}

\newcommand{\ds}{\displaystyle}
\renewcommand{\d}{\delta}
\newcommand{\D}{\Delta}
\newcommand{\e}{\varepsilon}
\newcommand{\g}{\gamma}
\newcommand{\G}{\Gamma}
\renewcommand{\i}{\infty}

\renewcommand{\L}{\Lambda}

\renewcommand{\o}{\omega}

\newcommand{{\z}}{\mathbb Z}
\newcommand{\R}{\mathbb R}
\newcommand{\N}{\mathbb N}

\newcommand{\ue}{u^{\epsilon}}

\renewcommand{\b}{\beta}
\renewcommand{\d}{\delta}

\newcommand{\s}{\sigma}

\renewcommand{\o}{\omega}

\renewcommand{\i}{\infty}

\newcommand{\macp}{\beta}
\newcommand{\micp}{\alpha}
\newcommand{\ald}{\alpha'}

\newcommand{\cH}{{\mathcal H}}
\newcommand{\cF}{{\mathcal F}}
\renewcommand{\P} {{\mathcal P}}
\newcommand{\M} {{\mathcal M}}
\newcommand{\cN} {{\mathcal N}}
\newcommand{\supp}{\mathop{\mathrm{supp}}}

\begin{document}

\begin{center}
{\Large\bf The kinetic limit of a system of coagulating Brownian 
particles}\bs\bs\\
{\sc Alan Hammond$^1$}  \\   Department of
Statistics\\ University of California\\
Berkeley, California 94720 \bs \\
{\sc Fraydoun Rezakhanlou\footnote{Research
supported in part by NSF grant DMS0307021}}
  \\ Department of
Mathematics\\ University of California\\
Berkeley, California 94720--3840

\end{center}

\begin{section}{Introduction}

Understanding the evolution in time of macroscopic quantities such as
pressure or temperature is a central task in non-equilibrium
statistical mechanics. We study this
problem rigorously for a model of mass-bearing
Brownian particles that are prone to coagulate when they are close,
where the macroscopic quantity in this case is the density of
particles of a given mass.
Brownian motion arises in the particles of a colloid, due to the
random agitation of the much smaller molecules that form the ambient
environment. As such, our model could be considered as one of a
colloid, where the dominant interaction between particles is that of
coagulation. A theoretical discussion of coagulation in colloids was
undertaken by Smoluchowski in \cite{smol}.

In the model that we consider, a large number $N$ of particles, each
carrying some integer-valued mass, are,
at some initial time, scattered in $\mathbb{R}^d$, whose dimension $d$
satisfies $d \geq 3$ for the purposes of this paper.
These particles then perform Brownian motions. There is an
$N$-dependent parameter $\epsilon$ that specifies the range of
interaction of any particle in the model: a pair of particles is
liable to coagulate (to form a new particle that combines the mass of
the old two) when the distance between the two is of order
$\epsilon$. The choice of $\epsilon$ as a function of $N$ is dictated
by insisting that the so-called mean free path is bounded away from
zero and infinity in the limit $N \to \infty$
of high particle number. (The mean
free path is the mean time until the first collision of a particle
drawn uniformly at random at the initial time. A scaling that
produces a bounded mean free path is called a 
kinetic limit.)
Our model
incorporates a significant degree of physical realism, absent from
earlier work on this type of problem, in the sense
that we permit the diffusion rate of the particle to depend on its
mass, including the case where this rate is taken to be decreasing in
the mass (it is physically reasonable to suppose that the
diffusion rate of a Brownian particle is inversely proportional to the
mass). As we will later describe in precise terms, we also introduce
a parameter into the mechanism of reaction which allows us to study
such reactions over a natural range of their strengths.

We study the macroscopic evolution of this particle system
by
measuring the density of particles of a given mass $m$ in the vicinity
of a macroscopic location $x$ and at some time $t$. We will prove
that, when the initial number of particles is chosen to be high, this
density typically evolves as the solution of the Smoluchowski system
of PDE,
\setcounter{equation}{0}
\begin{equation}\label{syspde}
\frac{\partial{f_n}}{\partial t}(x,t)  = d(n)\Delta f_n (x,t)  +
Q^n_1(f) (x,t) -
Q^n_2(f) (x,t) \, \qquad \qquad n = 1,2, \ldots
\end{equation}
with initial data $f_n(\cdot,0) = h_n(\cdot)$, to be specified in more
detail shortly. The first term on the
right-hand-side of (\ref{syspde}) corresponds to the diffusion among
particles of mass $n$, with $d(n)$ being one-half of the diffusion
rate of such particles.
The terms in (\ref{syspde}) corresponding to the interaction of pairs 
of particles are
given by the gain term
\begin{equation}\label{gainterm}
Q^n_1(f) (x,t)  = \frac 12  \sum_{m=1}^n \macp (m,n-m) f_m(x,t)
f_{n-m}(x,t),
\end{equation}
and the loss term
\begin{equation}\label{lossterm}
    Q^n_2(f) =   f_n (x,t) \sum_{m=1}^{\infty} \macp (m,n) f_m (x,t).
\end{equation}
Here, the collection of constants $\macp:\mathbb{N}^2 \to (0,\infty)$
quantify the macroscopic propensity of mass at a pair of values to
combine. As well as deriving the system (\ref{syspde}) as the typical
macroscopic profile of our random model, we prove in this paper the
precise relation between the macroscopic constants $\macp$ and the
microscopic mechanism of reaction.

We will be concerned with weak solutions of the system (\ref{syspde}), defined by the equality of the left- and right-hand sides of (\ref{syspde}) after multiplication by $J_n:\R^d \times [0,\infty) \to [0,\infty)$, integration in space-time, and integration by parts. The equality is demanded over all choices of sequences of compactly supported smooth functions $\big\{ J_n: n \in \N \big\}$
such that only finitely many terms in the sequence are not identically zero.

We now give a precise definition of the microscopic process. We in
fact define a sequence of such models,  indexed according to the
initial number $N$ of particles in them. We define a range of
interaction, $\epsilon$, according to $N \epsilon^{d-2} = Z$, where
the exact value of the positive constant $Z$ will shortly be
given. (We will explain heuristically at the beginning of Section \ref{sketpr}
why this relation between $N$ and $\epsilon$ determines the regime of
bounded mean free path). In defining the model, the main elements to describe
are the initial random choice of particle locations and masses, the
diffusive dynamics, and the mechanism for coagulation.

To describe each of these, we require notation for labelling the
particles in this time-dependent model. Let
a countable set $I$ of symbols be given. A configuration $q$ of
particles
is an
$\mathbb{R}^d \times \mathbb{N}$-valued function on a finite subset
$I_q$ of $I$. For any $i \in I_q$, $q(i)$ may be written as
$(x_i,m_i)$. The particle labelled by $i$ has mass $m_i$ and location
$x_i$. In practice, the index set $I_q$ will be a function of time,
with a change occurring only at collision events, of which there are
finitely many in any given sample of one of the random models.

As for the dynamics of the process, the action on $F$ of the
infinitesmal generator $\mathbb{L}$  is given by
\begin{equation}\label{copr}
    (\mathbb{L} F)(q)  = \mathbb{A}_0 F (q) +  \mathbb{A}_C F (q),
\end{equation}
where $F :
    \{ \mathbb{R}^d \times\mathbb{N}  \}^I \to [0,\infty)$ denotes a
    smooth function, its domain being given the product topology.
In (\ref{copr}), the diffusion operator $\mathbb{A}_0$ is given by
\begin{equation}\label{diff}
\mathbb{A}_0 F (q) =
\sum_{i \in I_q}{d(m_i) \Delta_{x_i} F } ,
\end{equation}
while the collision operator $\mathbb{A}_C$ is specified by
\begin{eqnarray}\label{dyna}
\mathbb{A}_C F (q) & = &\frac 12
    \sum_{i,j \in I_q}{\epsilon^{-2} V \Big( \frac{ x_i - x_j
    }{\epsilon} \Big) \micp(m_i,m_j)}  \\
    & & \qquad \qquad \times \
\bigg[ \frac{m_i}{m_i +
    m_j} F \big( S^1_{i,j} q \big) + \frac{m_j}{m_i +
    m_j} F \big( S^2_{i,j} q \big) - F(q) \bigg]. \nonumber
\end{eqnarray}
Here,
\begin{itemize}
\item the collection of constants $\micp: \mathbb{N}^2 \to
[0,\infty)$ are the parameters of strength of interaction between
pairs of particles of given integer mass, to which we earlier
alluded.
\item the function $V: \mathbb{R}^d \to [0,\infty)$ is assumed to be
continuous,  of
compact support, and  with $\int_{\mathbb{R}^d}{V(x) dx} = 1$. Its
role is to include among the models we consider a rule for
coagulation time that may be rather arbitrary, beyond the insistence
that it be Markovian and cause reaction of
a pair of particles at
some time when this pair are to be found within an order of the range
of interaction $\epsilon$.
\item we denote by $S^1_{i,j}q$ that configuration formed from $q$ by
removing the indices $i$ and $j$ from $I_q$, and adding a new index
from $I$ to which  $S^1_{i,j}q$ assigns the value $(x_i,m_i + m_j)$. The
configuration  $S^2_{i,j}q$ is defined in the same way, except that it
assigns the value  $(x_j,m_i + m_j)$ to the new index. The specifics
of the collision event then are that the new particle
appears in one of the locations of the two particles being removed,
with the choice being made randomly with weights proportional to the
mass of the two colliding particles.
\end{itemize}
While this choice of new particle location causes a
pleasing cancellation of terms at numerous times in our proof, it is
inessential: the same of method of proof yields our main theorem
for any choice of location for a new particle that is microscopically
close - that is, within an order of $\epsilon$ - from the place where
the two old particles reacted.

We explain further our reasons for choosing the form of the collision
term in (\ref{dyna}), and the
interpretation of the various terms. Suppose that two particles
$(x_i,m_i)$ and $(x_j,m_j)$ are such that $x_i - x_j$ enters the
support of $V_{\epsilon}( \cdot ) : = \epsilon^{-2} V \big(
\frac{\cdot}{\epsilon}\big)$ at time $t_0$. They are likely to remain
at a displacement of order $\epsilon$ for a time period of order
$\epsilon^2$. Given that the dimension $d$ is at least three, it is
likely that after a time that is large compared to
$\epsilon^2$, the particles have moved apart to a greater distance, and
that they will not return to an $\epsilon$-vicinity of each other at
any later time. This means that the integral
\begin{displaymath}
    I_T = \int_{t =t_0}^{T}{\micp(m_i,m_j) V_{\epsilon} \big( x_i - x_j
\big)} dt,
\end{displaymath}
whose derivative in $T$ is the infinitesimal rate at time $T$ for coagulation
between the pair, and which increases as $T$ rises, is likely to have
gained its eventual value as a function of $T$ when $T$ is of the
order of $t_0 + C \epsilon^2$. Thus, the normalisation by
$\epsilon^{-2}$ that appears in (\ref{dyna}) ensures that the integral
    $I_T$ reaches an eventual value that is of unit order. Neglecting
the effects of other particles, the probability of collision between
the pair before time $T$ is equal to $1 - \exp\big\{ - I_T \big\}$.
We are thus prescribing a dynamics in which the fraction of instances
of pairs of particles of displacement of
order $\epsilon$ that result in a coagulation of that pair is
bounded away from $0$ and $1$, uniformly in small values of
$\epsilon$. The role of the constants $\micp$ is to determine
whether this fraction is high or not.
The reason then for introducing the constants $\micp$ is that they
range over the right scale to determine precisely how liable
mass-pairs are to react microscopically. The choice of the value of
these constants transmits to the macroscopic reaction propensities
appearing in the interaction terms (\ref{gainterm}) and 
(\ref{lossterm}) of the Smoluchowski PDE. The recipe for determining 
$\macp$ from $\micp$ is as
follows: there exists a solution $u = u_{n,m}:
\mathbb{R}^d \to (0,\infty)$ of the equation
\begin{equation}\label{pdew}
    \Delta u_{n,m}(x)  \ = \  \frac{\micp(n,m)}{d(n) + d(m)} V(x) \Big[
1 + u_{n,m}
    (x) \Big]
\end{equation}
that is unique subject to the decay condition
$u_{n,m}(x) =O(|x|^{2-d})$ as $\vert x \vert \to
\infty$.
The quantities $\macp: \mathbb{N} \times \mathbb{N} \to
(0,\infty)$ in (\ref{syspde}) are then specified by the formula
\begin{equation}\label{recp}
  \macp (n,m)= \micp(n,m) \int_{\mathbb{R}^d}
    V(x)(1 + u_{n,m} (x)) dx.
\end{equation}

Our main result is conveniently expressed in terms of the empirical
measures on the locations of particles of a given mass. For each $n
\in \mathbb{N}$ and $t\in [0,\infty)$, we write
$g_n(dx,t)$ for the random measure on $\mathbb{R}^d$ given by
$$
g_n(dx,t)
=\e^{d-2}\sum_{i\in I_{q(t)}}\d_{x_i(t)}(dx)1\!\!1\big(m_i(t)=n\big).$$
We write $g$ for the random measure on space-mass-time $\mathbb{R}^d \times \N \times [0,\infty)$ such that, for each $t \geq 0$, its time-$t$ marginal $g(\cdot,t)$ is given by
$$
g(\cdot,t)
=\e^{d-2}\sum_{i\in I_{q(t)}}\d_{\big(x_i(t),m_i(t)\big)}.
$$

We also require a mild hypothesis on the diffusion coefficients
$d:\mathbb{N} \to (0,\infty)$. Namely, we
suppose that there exists a function $\g : \mathbb{N}^2 \to 
(0,\infty)$ such that $\micp \le \g$,
with $\g$ satisfying
\begin{equation}\label{hypo}
n_2 \g \big( n_1, n_2 + n_3 \big) \max{ \Big\{ 1, \Big[ \frac{d (
n_2 + n_3 )}{d(n_2)}
\Big]^{\frac{3d - 2}{2}},\Big[ \frac{d (
n_2 + n_3 )}{d(n_2)}
\Big]^{2d-1}\Big\}} \leq \big( n_2 + n_3 \big) \g(n_1,n_2).
\end{equation}

The initial random configuration of $N$ particles is formed by
scattering particles of numerous masses independently in
$\mathbb{R}^d$ according to densities that are prescribed for each
mass. These densities will be chosen as continuous functions
$\{ h_n: \mathbb{R}^d \to [0,\infty), n \in \mathbb{N} \}$, and should
satisfy some fairly weak bounds. To be specific, we insist that
\begin{itemize}
\item $\bar k*(-\Delta )^{-1}k\in L^{\i}_{loc}({\mathbb{R}^d})$ where
$k:=\sum_{n=1}^{\i} nh_n$ and $\bar k(x)=k(-x)$.
\item For every $m$, $\sum_{n=1}^{\infty}d(n)^{d/2-1/3}\gamma(m,n){\hat h_n} \in
L^\infty_{loc}({\mathbb{R}^d})$ where
$\hat h_n(x)=\int  h_n(y)|x-y|^{-d+2/3}dy$.
\item For every m, $\sum_{n=1}^{\infty}d(n)^{d/2-1/4}\gamma(m,n){\tilde h_n} \in
L^\infty_{loc}({\mathbb{R}^d})$ where
$\tilde h_n(x)=\int  h_n(y)|x-y|^{-d+1/2}dy$.
\end{itemize}
We then set
$Z=\sum_{n=1}^{\infty}\int_{\mathbb{R}^d}{h_n} \in (0,\infty)$ and choose $N$
points
in $\mathbb{R}^d  \times \mathbb{N}$ independently according to a law
whose density at $(x,n)$ is equal to $h_n(x)/Z$.
Selecting arbitrarily a set of $N$ symbols $\{i_j: j \in \{1,\ldots,
N\} \}$ from $I$, we define the initial configuration $q_0$ by insisting
that $q_0(i_j)$ is equal to the $j$-th of the randomly chosen members
of  $\mathbb{R}^d  \times \mathbb{N}$.

\noindent
{\bf Remark.}
It is not hard to show that our assumptions on the initial data $\{h_n\}$
are satisfied if $k$ is bounded, $k$ has a bounded support,
$d(\cdot)$ is bounded, and $\gamma(m,n)\le C(m)n$ for a function
$C(\cdot)$. Indeed if $k$ is bounded and has a bounded support, then
$\bar k*(-\Delta)^{-1}k\in L^\i_{loc}$ and $\hat h_n,\tilde h_n\in 
L^\i_{loc}$ for every
$n$.  It is worth mentioning that if $k$
belongs to the negative Sobolev Space $H^{-1}=W^{-1,2}$, then $\bar
k*(-\Delta)^{-1}k\in L^\i$.

The next statement is the main result of this paper.
\setcounter{theorem}{0}
\begin{theorem}\label{thmo}
Let $d \geq 3$ and suppose that~(\ref{hypo}) holds. Let
$\mathcal P_N$ denote the law on measures on $\R^d \times \N \times [0,\infty)$ given by the law of $g$ under $\mathbb P_N$;  recall that $\epsilon$ is related to $N$ by means of the
formula $N
\epsilon^{d-2} = Z$, with the constant $Z \in (0,\infty)$ being
given by the expression $Z = \sum_{n \in
\mathbb{N}}{\int_{\mathbb{R}^d}h_n}$.

The sequence $\{\mathcal P_N\}$ is 
tight. Moreover, any limit point $\mathcal P$ of the sequence $\{\mathcal P_N\}$ is concentrated on the space of measures taking the form $\sum_{n=0}^{\infty} f_n(x,t) dx \times \delta_n \times dt$ where 
$\{f_n: n \in \mathbb N\}$ ranges over weak solutions of (\ref{syspde}) that satisfy
the initial condition $f_n(\cdot,0) = h_n(\cdot)$; recall that  the collection of
constants $\macp: \mathbb{N}^2 \to [0,\infty)$ is given by (\ref{recp}).
\end{theorem}

Theorem~\ref{thmo} describes the evolution of the density profiles of particles of various masses in the limit of large particle number by means of the Smoluchowski PDE. 
Note that, if the weak solution of this system of PDE is not known to be unique, we merely demonstrate convergence in a subsequential sense to the space of solutions. For example, admitting the possibility that the system~(\ref{syspde}) has two distinct weak solutions $\big\{ f'_n: n \in \N \big\}$ and $\big\{ \hat{f}_n: n \in \N \big\}$ with initial condition $f_n(\cdot,0) = h_n(\cdot)$, each of the following behaviours is consistent with Theorem~\ref{thmo}:
\begin{itemize}
 \item the empirical densities under the microscopic models $\mathbb{P}_N$ may converge weakly to the solution $\big\{ f'_n: n \in \N \big\}$ as $N \to \infty$ along the subsequence of even integers, and to $\big\{ \hat{f}_n: n \in \N \big\}$ as $N \to \infty$ along the subsequence of odd integers;  
 \item it may be that evolution of these densities is accurately approximated by flipping a fair coin, with the densities converging weakly to  $\big\{ f'_n: n \in \N \big\}$ as $N \to \infty$ should the outcome be heads, and to   $\big\{ \hat{f}_n: n \in \N \big\}$ as $N \to \infty$ should the outcome be tails.
\end{itemize}

These peculiar scenarios are excluded if uniqueness of solutions to~(\ref{syspde}) is known. Some conditions for uniqueness are furnished by Proposition 2.6 of \cite{wrzd}; subsequently to the first appearance of the present work,  \cite{momentbounds} provided uniqueness under quite weak hypotheses. Indeed, as Remark 1.2 of~\cite{momentbounds} discusses, the next proposition is a consequence of Theorems 1.1, 1.2, 1.3 and 1.4 of~\cite{momentbounds}. 
\begin{prop}\label{propuniqueness}
Let the dimension satisfy $d \geq 1$. 
For $a,b > 0$ such that $a + b < 1$, and for positive constants $c_1$ and $c_2$, assume that 
$\macp(n,m)\le c_1(n^a+m^a)$ and $d(n)\ge c_2 n^{-b}$ for
all $n,m\in \N$. Also assume that $d:\N \to (0,\infty)$ is non-increasing. There exists $e > 0$ such that $\sum_nn^e \|h_n\|_{L^\infty(\R^d)}<\i$
 and $\|\sum_nn^e h_n\|_{L^1(\R^d)}<\i$ imply that~(\ref{syspde}) has a unique weak solution. (In fact,  this solution conserves mass, in the sense that  $I:[0,\infty) \to [0,\infty)$ given by $I(t)  = \sum_{m \in \mathbb{N}} m \int_{\mathbb{R}^d} f_m(x,t)dx$ satisfies $I(t) = I(0)$ for all $t \in [0,\infty)$.)
\end{prop}

It is a simple corollary of Theorem~\ref{thmo} and Proposition~\ref{propuniqueness} that convergence to~(\ref{syspde}) in fact holds in the following stronger sense.

\begin{cor}\label{coro}
Let $d \geq 3$ and suppose that~(\ref{hypo}) and the assumptions of Proposition~\ref{propuniqueness} are in force.
Let $J: \mathbb{R}^d 
\times [0,\infty) \to \mathbb{R}$ be a 
 bounded and continuous test function. Then, for each
$n \in \mathbb{N}$ and $t \in (0,\infty)$,
\begin{equation}\label{resu}
\limsup_{N \to \infty } \mathbb{E}_N  \bigg\vert  \int_{\mathbb{R}^d}
J(x,t) \big( g_n(dx,t) - f_n (x,t) dx \big) \bigg\vert = \ 0,
\end{equation}
where again $N
\epsilon^{d-2} = Z$, with $Z = \sum_{n \in
\mathbb{N}}{\int_{\mathbb{R}^d}h_n}$. In (\ref{resu}),
    $\{ f_n: \mathbb{R}^d \times [0,\infty) \to [0,\infty), n \in \mathbb{N} \}$
denotes the unique weak solution to the system of
partial differential equations (\ref{syspde}), with $\macp: \mathbb{N}^2 \to [0,\infty)$ again given by (\ref{recp}).
\end{cor}
\noindent{\bf Remarks.} 
\begin{itemize} 
\item
Included in the space of parameter values that satisfy
(\ref{hypo}) is the case where the diffusion rate $d$ is a decreasing
function of the mass, and the coagulation propensities $\micp$
satisfy $\micp \big( n , m \big) \leq C nm$. In fact for a 
nonincreasing $d(\cdot)$, the condition (\ref{hypo})
is equivalent to saying that $\micp(n,m)\le C(n)m$ for a function $C(n)$.
Also, if the microscopic coagulation
rate $\micp$ is identically constant, then the condition (\ref{hypo}) 
is equivalent to saying that the function
$d(n)n^{1/(2-3d)}$ is nonincreasing.
\item
We believe that the derivation that we have undertaken is of
particular interest, because the technique of proof that we adopt
is robust under numerous changes in the details of the mechanism of
interaction. An alternative model which we may treat is that where the
mechanism of interaction is {\it hard-core}. In this variant, the
reaction mechanism is  deterministic, a pair of particles coagulating at the moment when
their locations differ by $2\epsilon$. The hard-core model with constant
diffusion rates  $d(\cdot) \equiv 1$ was treated in dimension $d=3$ in
\cite{Lang}, with the macroscopic reaction propensity $\macp$ being
identically given by the capacity $ 2 {\rm cap} (B_0(1))$ of the unit-ball. Although we have not written the proof in this paper
using a notation that includes the case of hard-core interaction, the
method of proof yields the results of \cite{Lang}. Indeed, we can
prove the analogue of Theorem \ref{thmo} for hard-core interaction in
any dimension $d \geq 3$, provided that we may invoke the hypothesis 
\begin{equation}\label{hardcorehyp}
n_2  \max{ \Big\{ 1, \Big[ \frac{d (
n_2 + n_3 )}{d(n_2)}
\Big]^{\frac{3d - 2}{2}},\Big[ \frac{d (
n_2 + n_3 )}{d(n_2)}
\Big]^{2d-1}\Big\}} \leq  n_2 + n_3  .
\end{equation}
The relation $\macp(n,m) =  \big( d(n) + d(m) \big)
{\rm cap} \big( B_0(1) \big)$ holds.  The hard-core model is anyway approximated 
in terms of its macroscopic evolution by
the stochastic one, with high choices for the constants $\micp$, as we
will see in Theorem  \ref{th6.2}:  if the  support of $V$ is the unit ball $B_0(1)$ and $n,m\in
\mathbb{N}$ are given, then the limit $\lim_{\micp(n,m)
\to \infty}{\macp(n,m)}$ exists and equals $ \big( d(n) + d(m) \big)
{\rm cap} \big(B_0(1)\big)$,
the limit being taken
with $d(n)$ and $d(m)$ fixed. 
\item Our technique of proof also yields a kinetic limit derivation
for the model in which particles are assumed to have a range of
interaction that is mass-dependent. To give an example of such
variants, suppose that each
particle of mass $m$ has a radius $r(m)$, where $r(m) = m^{\chi}$,
with $\chi$ a fixed value satisfying 
$\chi \in \big( 0, 1/(d-2) \big)$.
 Assume, for the sake of brevity, that the diffusion rate $d$ is uniformly bounded above
and below. We stipulate that particles of mass $m$ and
$n$ are liable to react when their displacement reaches the order of
$\big( r(m) + r(n) \big) \epsilon$. More precisely, we modify the
definition (\ref{dyna}) of the collision operator $\mathbb{A}_{C}$ by
replacing the appearance of $V $ by $ \big( r(n) + r(m) \big)^{-2}  V \big(
\cdot/(r(n) + r(m))\big)$, (the factor that multiplies $V$ being
introduced so that, roughly speaking, the altered collision mechanism
respects 
 the spatial-temporal scaling of Brownian motion). The statement of Theorem
\ref{thmo} is altered in that the macroscopic reaction propensities are
 given by (1.8) with $u_{n,m}$ solving (1.7) where $V$ is replaced
 by $ \big( r(n) + r(m) \big)^{-2}  V \big(
\cdot/(r(n) + r(m))\big)$. (The condition that $\chi$ be bounded above
as stated 
is required in the
proof of Lemma \ref{lemthreep}. Varying the details of this proof is
the principal step required in extending our investigation to this
altered model.) 
We mention also that investigating this
model in the case where the radius $r$ is chosen to be a more rapidly
increasing function of mass may be highly interesting, because it may
be used to describe interacting systems of polymers having a fractal
structure, in which a gelation effect occurs macroscopically, in the
sense that the system of equations (\ref{syspde}) lacks a well-defined
solution after some finite time before which particles of high mass acculumate.
\item In Section \ref{sketpr}, we will discuss how the product of densities
that appears in the interaction terms (\ref{gainterm})
and (\ref{lossterm}) of the macroscopic evolution is a consequence of 
the asymptotic independence for each
positive time of particle locations in the large random models.
As we will discuss, this independence fails to hold in a large
microscopic vicinity of the location of a given  particle: it is less
probable that at any given positive time,
a particle is present at a distance of the order of $\epsilon$
  from a given particle than at a randomly chosen point in space,
  because the mechanism of coagulation has the effect of clearing space
  about any given particle. The endeavour of deriving the relation 
(\ref{recp}) between the
  macroscopic and microscopic reaction propensities has been a
  formidable one, since in so doing,
  we have determined precisely the nature of this negative microscopic
  correlation which forms a correction to the overall asymptotic
  independence of particle locations at positive times. The authors of
  \cite{ESY} accomplished an analogous task for 
the quantum mechanical problem of a Bose-Einstein condensate with a
quite singular choice of pair potential,
in which the negative microscopic correlation corrects the
Born approximation. An assumption that removed interaction between
particles in small regions where several were present was employed.
\end{itemize}

We conclude this introduction by mentioning some previous work
related to the problem.

As we just mentioned, the project of deriving the kinetic limit of the particle
densities that we undertake was treated for the hard-core model, in the case
where the diffusion rates are identically equal to a constant, in
\cite{Lang}. The approach of this paper was that by which Lanford
\cite{Lanford} had derived for short times Boltzmann's equation as a
similar limit of
a system of balls that move at constant velocities and experienced
elastic collisions. That is, the authors verified that the
correlation functions of the particle locations satisfied the BBGKY
hierarchy, a system of equations that may be derived formally for
such systems
by
assuming sufficient independence in the location of
particles. Similarly to Lanford's approach, the derivation in
\cite{Lang} was made on
a short time interval. However, an iteration allowed it to be extended
for all time.

A related problem has been studied by Sznitman
\cite{sznitman}. Brownian spheres diffuse and annihilate as soon as
they touch. The partial differential equation by which the density of
particles evolves was derived for the kinetic limit, in each dimension
$d \geq 2$.

The stochastic coalescent is a random model in which
particles - elements of an abstract set - are assigned a non-negative
mass lying in $[0,\infty)$, and
pairs of particles of mass $x$ and $y$ coagulate at an exponential
rate $K(x,y)$. Norris \cite{norrisone} showed that for many choices of
interaction rate $K$, in the limit of
large initial particle number, the density of particles evolves in
time according to a spatially homogeneous analogue of (\ref{syspde}),
where the densities depend on time alone and the Laplacian term is
absent. The coagulation rate $\macp$ in this case is the same as in
the random model: $\macp(x,y)=K(x,y)$. The derivation was valid
provided that this system of partial differential equations had an
unique solution, sufficient conditions for which were proved.

In \cite{GKO}, a derivation of a system similar to (\ref{syspde}) is
made from a microscopic model with a finite number of possible mass
types, and a mechanism for interaction that is over a long-range compared to
the microscopic scale.

In the next section, we give an overview of the proof of Theorem
\ref{thmo}, describing the organisation of the paper as we do so.

\medskip

\noindent{\bf Acknowledgment.} The first author would like to thank James Norris for introducing him to the topic of diffusive coagulating systems and for valuable discussions.

\end{section}

\begin{section}{A sketch of the proof}\label{sketpr}

An essential element in the proof of Theorem \ref{thmo} is to verify
that the distribution of particles of two distinct masses
is independent in the high $N$ limit. This is the so-called
\stoss of Boltzmann. Before sketching how we perform this
task, we illustrate how assuming a form of the \stoss
tells us that  the choice of scaling $N =
\epsilon^{2-d} Z$ in the model identifies the regime of the bounded mean
free path: that
is, why the number of collisions experienced per unit time by a
typical particle is of unit order, when $N$ is taken to be very
large.  In this regard, recall that the volume swept out by a ball of 
radius $\epsilon$
about the location of a particle in the first unit of time is that of
the `Wiener sausage'  of parameter $\epsilon$, and is of the order of
$\epsilon^{d-2}$.
The number of collisions is of the order of the number of instances
that, in this unit of time, other particles lie in this swept-out
volume. The other particles are distributed independently at later
times - at least, if we admit some form of the \stoss -
so that we expect the particle
to encounter of the order of $N \epsilon^{d-2}$ other particles. The
scaling of $N$
and $\epsilon$ ensures that this remains bounded above and away from zero as
$N$ increases.

In the form in which we prove it,  the \stoss asserts that
the total `propensity to coagulate' between particles of mass $M_1$
and $M_2$, weighted spatially by some time-dependent test functions
$J,\overline{J} :
\mathbb{R}^d \to [0,\infty) \times \mathbb{N} \times [0,\infty)$ in
the first $T$ units of time,
\setcounter{equation}{0}
\begin{equation}\label{quez}
\frac{\epsilon^{d-2}}2 \int_{0}^{T} \sum_{i,j \in I_q}{
    V_{\epsilon} \big( x_i - x_j \big)  \micp(m_i,m_j)  J(x_i,M_1,t)
    \overline{J}(x_j,M_2,t) 1 \! \! 1 \Big\{ m_i = M_1,m_j=M_2 \Big\} dt},
\end{equation}
is close, when $N$ is high enough, to an expression involving a
space-time integral of products of the candidate densities:
\begin{equation}\label{tapd}
 \frac{\macp \big( M_1, M_2 \big)}{2} \int_0^T {\int_{\mathbb{R}^d}{
f^{\delta}(M_1,x,t)
f^{\delta}(M_2,x,t) J(x,M_1,t) \overline{J}(x,M_2,t) dx dt}}.
\end{equation}

The precise statement of this claim appears in Section \ref{stfive} as
Proposition \ref{szprop}.
In seeking to verify the claim, we introduce the random
variables $Q(z)$, or, more fully, $Q(z,M_1,M_2,t)$, given by
\begin{eqnarray}\label{sfd}
    & &\frac{ \epsilon^{d-2}}2 \sum_{i,j \in I_q}{
    V_{\epsilon} \big( x_i(t) - x_j(t) + z \big)  \micp \big(m_i(t),m_j(t)
    \big) } \\
     & & \qquad \qquad \qquad J
    \big( x_i(t),M_1,t \big)
    \overline{J} \big( x_j(t),M_2,t \big) 1 \! \! 1 \Big\{ m_i (t) = M_1,
    m_j (t) = M_2 \Big\}, \nonumber
\end{eqnarray}
where $z \in \mathbb{R}^d$.
In the case of $z=0$, $Q(0)$ measures the
average propensity per particle at time $t$ for particles of mass
$M_1$ and $M_2$ to
interact (with a spatial weighting due to the test functions). It is
an average due to the appearance of the factor of
$\epsilon^{d-2}$ in (\ref{sfd}), whose reciprocal is of the order of
the number of particles in the system initially. The total number of
particles is expected to remain of this order at any given later time,
because a typical particle collides with only finitely many others per
unit time. This means that, ignoring a constant factor, we may think
of $Q(0)$ as the average propensity to interact for any given value
$t$ of the time parameter. We also consider expressions
$\overline{Q}$, defined by the
formula in (\ref{sfd}), with the interaction kernel $V_{\epsilon}$ being
replaced by another, $\hat{V}_{\epsilon}$, having similar properties.
If $z$ is fixed at a small and
constant value, while the parameter $\epsilon$ is taken to be much
smaller, then any such $\overline{Q}(z)$ measures the extent of
appearances at time $t$ of
particles of mass $M_2$ lying in a microscopic ball of radius
$\epsilon$ centred about a small macroscopic displacement $z$ from
some mass $M_1$ particle.
Instances of such pairs of particles are
tentative candidates for collision after the passing of a short period
of macroscopic time. As such, the assertion that,
for fixed $t \in (0,\infty)$, the expression $\overline{Q}(z_1) -
\overline{Q}(z_2)$ is small
if $z_1$ and $z_2$ are small, captures the near independence of the
locations of
particles of distinct masses that is in essence the
\stoss . Indeed, we will demonstrate that, for a particular
variant $\overline{Q}$ of $Q$,
\begin{equation}\label{pkj}
    \lim_{z \downarrow 0}\mathbb{E}_N \Big\vert \int_0^T \Big( \overline{Q}(z) -
     Q(0) \Big) dt \Big\vert = 0,
\end{equation}
    the limit being taken in such a way that $\epsilon/\vert z \vert  \to
    0$ as $z \to 0$. Then, choosing $\eta^\delta$ to be a sequence that
    approximates the Dirac $\delta$ function, we find that
    $\int_0^T Q(0) dt$ is likely to be close to
\begin{equation}\label{laste}
    \int_0^T \int_{\mathbb{R}^d}  \int_{\mathbb{R}^d} \overline{Q}(z_1 -
z_2) \eta^\delta \big( z_1 \big)
    \eta^\delta \big( z_2 \big) dz_1 dz_2 dt.
\end{equation}
By substituting the expression in (\ref{sfd}) (evaluated with the
choice $z = z_1 - z_2$, and with $\hat{V}_{\epsilon}$ replacing
$V_{\epsilon}$) for $\overline{Q}(z_1 - z_2)$,
we find that the expression in (\ref{laste}) takes the form appearing in
(\ref{tapd}), that is, a time-averaged product of candidate
densities. (The details, which appear in the proof of Proposition
\ref{szprop}, are given in (\ref{bhap}) ).
We now turn to sketch how a statement of the form (\ref{pkj}) is derived.

In a first effort to understand the difference in the behaviour of the
time integral of the coagulation rate $Q(0)$ and its $z$-displaced
counterpart $Q(z)$, we write down
a $z$-dependent functional $S_z$ on configuration space for
which the action of the diffusion operator $\mathbb{A}_0$ on $X_z$ produces,
among others, the term $Q(z)$. For each pair $(n,m) \in
\mathbb{N}^2$, we define $\phi_{n,m}^{\epsilon}:
\mathbb{R}^d \to (0,\infty)$ by the
stipulation that \mbox{$\Delta \phi^\epsilon_{n,m} (x) =  \ald \big( n,m
\big) \epsilon^{-d}
V(z/\epsilon)$}, where we write  
\begin{equation}\label{aldo}
\ald(n,m) = \micp(n,m)/\big(d(n) +
d(m)\big) . 
\end{equation}
We then set
\begin{equation}\label{exz}
S_z (q) = \epsilon^{2(d-2)} \sum_{i,j \in I_q}{
\phi_{M_1,M_1}^{\epsilon} (x_i- x_j + z) J (x_i,M_1,t)
\overline{J}(x_j,M_2,t)} 1 \! \! 1 \big\{ m_i = M_1,m_j = M_2\big\}.
\end{equation}
When the operator
$\mathbb{A}_0$ acts, we find
that the term $Q(z) - Q(0)$ is obtained when the operator $\mathbb{A}_0$ acts
on $S_z - S_0$ with each of the two derivatives of the Laplacian falling
on $\phi^\epsilon$ rather than the test functions.
Consider the identity
\begin{eqnarray}\label{hac}
    \big( S_z - S_0 \big) \big( T \big) & = &
    \big( S_z - S_0 \big) \big( 0 \big) +
\int_{0}^{T}{ \left(\frac{\partial}{\partial t}+\mathbb{A}_0 
\right)(S_z - S_0) (t)dt} \\
    & & \qquad + \,
\int_{0}^{T}{ \mathbb{A}_C (S_z - S_0) (t)dt} \, + \, M_T, \nonumber
\end{eqnarray}
and note that the process $\big\{ M_T : T \geq 0 \big\}$ is a martingale.
We know that two terms lie inside this identity that remain of unit
order as we take a limit in low $z$
(with $\epsilon$
tending to $0$ much more quickly). These are $\int_0^T Q(z) dt$ and $-
\int_0^T Q(0) dt$.
Were the only terms to remain of
unit order to be those two, we would deduce
a conclusion of the form we seek, namely (\ref{pkj}) where $\overline{Q}$
is simply $Q$. However, (\ref{pkj}) cannot be valid for this choice of
$\overline{Q}$. Indeed, imagine that $\micp(M_1,M_2)$ is a very large
constant. Particles of order $\epsilon$ apart would typically
coagulate in a very small fraction of the order $\epsilon^2$ of time
in which they are likely to remain at a distance of order
$\epsilon$. At the moment that the pair of particles vanishes, it
ceases to contribute to the expression $Q(0)$, meaning that, at any
particular time $t \in [0,\infty)$, $Q(0)$ need not be large, even if
the constant $\micp(M_1,M_2)$ is. However, for fixed $z$, pairs of
particles, one of which lies in an $\epsilon$-ball about the
$z$-displacement of the other, may remain in this state for a
significant constant multiple of $\epsilon^2$ of upcoming time. This
means that $Q(z)$ is likely to be huge, if $\micp(M_1,M_2)$ is.

In fact, there is another term that remains of unit order in
(\ref{hac}) in the low $z$ limit. This term arises from the action of
the collision operator $\mathbb{A}_C$ on $S_0$. Two or three particles  may
index the summands of
such terms. Two such particles are those that may coagulate,
appearing in the sum associated to the collision
operator. There may or may not be another occurring as a third particle
in the term $S_0$. It is the case of the term corresponding to only two
particles that remains of unit order. It takes the form
\begin{equation}\label{geezee}
     \epsilon^{2(d-2)}  \micp(M_1,M_2) \sum_{k,l \in
I_q}{  V_{\epsilon}\big( x_k - x_l \big)}
\phi_{M_1,M_2}^{\epsilon}(x_k -
x_l ) J(x_k,M_1)  \overline{J}(x_l,M_2).
\end{equation}
This term witnesses the abrupt curtailment of the propensity to
coagulate of a pair of particles at that moment when they do
coagulate. Given the argument just presented against the possibility
that $\overline{Q}(0)$ is equal to $Q(0)$ in (\ref{pkj}), it is not
surprising that this extra term should arise in this way.

There are numerous other terms in (\ref{hac}). It turns out that each
of them vanishes in the limit of low
$z$. Accepting for the time being that this is so,
we have however not quite arrived at a conclusion of the form we seek. What we
have found is that $\int_{0}^{T} Q(z) dt$ is close, when $z$ is small,
to the expression
\begin{eqnarray}\label{geeque}
     &  &  \epsilon^{2(d-2)}  \micp(M_1,M_2) \int_0^T dt \sum_{i,j \in
I_q}{ V_{\epsilon}(x_k - x_l)}
\\
    & & \qquad \qquad \qquad \times \
     \Big[
1 - \phi_{M_1,M_2}^{\epsilon} \big( x_i - x_j \big) \Big] J(x_i,M_1)
\overline{J}(x_j,M_2) 1 \! \! 1 \big\{ m_i = M_1,m_j = M_2 \big\} . \nonumber
\end{eqnarray}
The `$1$' that appears in the square bracket corresponds to $Q(0)$,
and the other term to the unit order term (\ref{geezee}). It was our aim
to show that $\int_{0}^{T} Q(0) dt$ was close to an expression of the
form  $\int_{0}^{T} \overline{Q}(z) dt$, but our first try gave us
such a statement with the roles of $Q$ and $\overline{Q}$ reversed. We are led
therefore to consider an altered functional $X_z$ of the configuration
space:
\begin{equation}\label{exzo}
X_z (q) = \epsilon^{2(d-2)} \sum_{i,j \in I_q}{
\ue_{M_1,M_2} (x_i- x_j + z) J (x_i,M_1,t) \overline{J}(x_j,M_2,t)}  1
\! \! 1 \big\{ m_i = M_1,m_j = M_2 \big\} .
\end{equation}
We have replaced the appearance of $\phi^{\epsilon}_{m_i,m_j} \big(x_i- x_j
+ z\big)$ by that of
$\ue_{m_i,m_j} (x_i- x_j + z)$ to move from (\ref{exz}) to (\ref{exzo}).
Note that the analogue of the expression (\ref{geeque}) in this
case is given by
\begin{eqnarray}
     &  &  \epsilon^{2(d-2)} \sum_{i,j \in
I_q}{   \Big[ \big( d(M_1) + d(M_2) \big) \Delta
\ue_{M_1,M_2} \big( x_i - x_j \big) - \micp \big( M_1,M_2 \big)
    V_{\epsilon} \big( x_i - x_j \big) u^{\epsilon} \big( x_i - x_j) \Big]
}
\\
    & & \qquad \qquad \qquad \qquad \qquad   J(x_i,M_1) 
\overline{J}(x_j,M_2) 1 \!
\! 1 \Big\{ m_i
    =M_1, m_j =M_2\Big\}. \nonumber
\end{eqnarray}
We aim to define the function $\ue_{M_1,M_2}$ in such a way that this
term is equal to $Q(0)$. This means that
$\ue_{M_1,M_2}$ should satisfy
\begin{equation}\label{rgddg}
    \big( d(M_1) + d(M_2) \big) \Delta \ue_{M_1,M_2} (x) = \micp
\big(M_1,M_2 \big) V_\epsilon (x) \ue (x) +
    \micp \big( M_1,M_2 \big) V^{\epsilon} (x)
\end{equation}
for each $x \in \mathbb{R}^d$. Writing $\ue_{M_1,M_2}(x)=
\epsilon^{2-d} u_{M_1,M_2} \big(
x/{\epsilon} \big)$ and recalling the definition (\ref{aldo}) of the constants
$\ald(M_1,M_2)$, the equation (\ref{rgddg}) for any, or all, $\epsilon >
0$, is equivalent to the requirement
that $u_{M_1,M_2}: \mathbb{R}^d \to \mathbb{R}$ satisfies
\begin{equation}\label{hgj}
    \Delta u_{M_1,M_2} (x) = \ald \big(M_1,M_2 \big) V (x) \Big[ 1 +
u_{M_1,M_2} (x) \Big].
\end{equation}
In Section \ref{sec6}, we prove that the equation (\ref{hgj}) has an
unique solution
in $C^2(\mathbb{R}^d)$. This piece of theory allows us to define the
functional $X_z$ so that it has the desired effect: the term
$\int_{0}^{T} Q(0) dt$ is close to
$\int_{0}^{T} \overline{Q}(z) dt$ for small $z$. The form of
$\overline{Q}(z)$ is given by
\begin{equation}\label{ppo}
       \epsilon^{2(d-2)} \sum_{i,j \in
I_q}{ \big( d(M_1) + d(M_2) \big) \Delta \ue_{M_1,M_2}
\big( x_i - x_j + z \big)
     J(x_i,M_1)  \overline{J}(x_j,M_2) 1 \!\! 1 \big\{
     m_i(t)=M_1,m_j(t)=M_2 \big\} }.
\end{equation}
This expression explains the form of the modification as we pass from
the microscopic propensity to coagulate, given by $\micp(M_1,M_2)$, to
its macroscopic analogue $\macp(M_1,M_2)$. The recipe (\ref{recp}) for
computing $\macp$ is nothing other than
$$
    \macp \big(M_1,M_2 \big)  = \big( d(M_1) + d(M_2) \big)
    \int_{\mathbb{R}^d}{\Delta u_{M_1,M_2} (x) dx}
$$
an equation that arises because for small $z$, the significant behaviour in
(\ref{ppo}) is given by that of $u(x-y)$ on its diagonal $x=y$.

Our task then becomes that of examining the analogue of the identity
(\ref{hac}) in the case where $S_z$ is replaced by $X_z$.
Estimates are required to prove that the remaining terms, including
that of the martingale, vanish
in the low $z$-limit. It would suffice in this regard to demonstrate
the following form of independence in the behaviour of the
particles: that there exists a large constant $C$ such that, for any
$n \in \mathbb{N}$ and $t \in (0,\infty)$, if $n$
particles are picked independently at random from the initial
configuration, then the
probability that all of them lie in a given small ball of radius $r$
at time $t$ is bounded above by $C r^n$.
Although this statement may seem no easier to prove than the
\stoss, there is a special case for which it may be obtained
by a coupling argument. This case is the one treated by \cite{Lang},
when all of the particles are supposed to have the same diffusion rate
(although the argument presented here does not require the mechanism of
coagulation to be of hard-core type). With this assumption, consider
an altered model in which there are particles of two types, those living
and those perished. Each particle is living at the beginning. Collisions may
occur only between living particles. When a collision occurs, one particle
remains alive and the other perishes, with each particle having an
equal chance of surviving, independently of other randomness.
Particles of either type diffuse
according to independent Brownian motions at one given rate.
The ensemble of all particles is a system of independently evolving
Brownian motions, and is at any given time a superset of the living
particles, which evolve according to our collision dynamics. The
estimate mentioned follows straightaway from this fact.

In general, however, no such approach is available. After labelling
the various terms arising from (\ref{hac}) at the beginning of
Section \ref{stwo},
we instead establish various weakened statements in Subsection
\ref{stone}, of which that in Lemma \ref{lemthreep}
is an example whose proof is now discussed. It asserts that
\begin{equation}\label{qax}
\mathbb{E}_N{\int_0^T{dt \sum_{k,l,i \in I_q} K_{n_1,n_2,n_3}(x_k,x_l,x_i)
1 \! \! 1 \big\{ m_k = n_1, m_i =
n_3 \big\} }} \leq C \epsilon^{6 - 2d},
\end{equation}
where $K$ is an instance of a smooth function mapping
$\mathbb{R}^{3d}$ to $[0,\infty)$, and where $C$ is a constant that
depends on $n_1$, $n_3$ and $K$. The approach to proving such
statements as (\ref{qax}) has much in common with the global one by
which (\ref{pkj}) is being obtained. We define a function
    $A_{n_1,n_2,n_3}^{\epsilon}: \mathbb{R}^{3d} \to [0,\infty)$ on which
    the action of the diffusion operator $\mathbb{A}_0$ yields the generic
    summand in (\ref{qax}). That is, we insist that
\begin{displaymath}
    \big( d(n_1 )\Delta_1 + d(n_2) \Delta_2 + d(n_3) \Delta_3 \big)
    A_{n_1,n_2,n_3}^{\epsilon} (x_1,x_2,x_3)
    = -  K_{n_1,n_2,n_3}^{\epsilon} (x_1,x_2,x_3)   .
\end{displaymath}
We consider the action of the whole Markov generator $\mathbb{A}_0 +
\mathbb{A}_C$ on
the sum $A^\epsilon_{n_1,n_2,n_3}(x_i,x_j,x_k)$ over those triples of
indices $(i,j,k)$ for which $m_i = n_1$ and $m_k = n_3$. We work
inductively in the lexicographical ordering in $\mathbb{N}^2$ for the
pair $(n_1,n_3)$. The action of $\mathbb{A}_0$ on this sum naturally provides
the expression that we seek to bound in (\ref{qax}). The action of the
collision operator $\mathbb{A}_C$ produces three types of terms, that might be
called gains, losses or internal terms. The gains and losses arise
from collisions that produce or destroy particles of mass $n_1$ or
$n_3$ that contribute to the expression in (\ref{qax}) in its first or
third arguments. The gain term appears in an upper bound for this
expression, and so must be bounded from above. Such a bound is
provided by the sum of the loss terms from particles of lower mass,
terms which are themselves bounded by the inductive hypothesis.
The internal term arises from coagulations that remove two particles
that contribute terms in the left-hand-side of (\ref{qax}) by means of
the unrestricted second argument of $K$, and produce a new particle
that itself contributes such a term. We argue that the internal term
is non-positive, because the term corresponding to the new particle is
less than the sum of the terms corresponding to the two that are
disappearing. It is at this point that the hypothesis (\ref{hypo}) on
the model's parameters is used.

Some additional lemmas presenting straightforward bounds on the
function $\ue$ are given in Subsection \ref{sttwo}.
For the analogue of (\ref{hac}) for the
process $X_z - X_0$, estimates
showing that numerous terms are negligible in the low
$\epsilon$ limit are given in
Subsection \ref{stthree}. The martingale term is similarly shown to be small in
Subsection \ref{stfour}.
In this way, the statement (\ref{pkj}) is derived. The derivation of
the \stoss, as described after (\ref{pkj}), is given in Subsection
\ref{stfive}.

We aim to show that the macroscopic mass densities evolve according to
the system of partial differential equations (\ref{syspde}) as
follows.
In Section \ref{derpde}, we define the
functional
\begin{displaymath}
    Y_{q} = \epsilon^{d-2} \sum_{i \in I_q}{ J \big( x_i,m_i,t \big) },
\end{displaymath}
where $J:   \mathbb{R}^d \times \mathbb{N}\times[0,\infty) \to
\mathbb{R}$ denotes a test function. Consider the identity
\begin{equation}\label{betlab}
    Y_{q}\big( T \big) = Y_{q} \big( 0 \big) + \int_0^T \mathbb{A}_C
    \big( Y_{q} \big) (t) dt   + \int_0^T 
\left(\frac{\partial}{\partial t}+\mathbb{A}_0
\right)
    \big( Y_{q} \big) (t) dt  \, + \, M_T,
\end{equation}
with the process $\big\{ M_T: T \geq 0 \big\}$ being a
martingale. After demonstrating that the martingale term is negligible
by arguments similar to those employed for the martingale arising from
the previous functional $X_z$, we consider the identity obtained from
(\ref{betlab}) in the limit of low $\epsilon$. We expect to derive the
system (\ref{syspde}) in weak form, with the function $J$ playing the
role of the test function. As discussed after (\ref{pkj}), the
\stoss enables us to replace the term involving the collision
operator in (\ref{syspde}) by a time-averaged product of microscopic
candidates for the density, that take the form of the empirical
measure for particles of a given mass mollified on some scale
$\delta$. To pass to the limit in low $\delta$ after $\epsilon$ has
vanished and in this way obtain the weak form of (\ref{syspde}), we
make use of the fact that these empirical measures concentrate in the
limit of low $\epsilon$ on measures having densities that are
uniformly bounded with respect to Lebesgue measure.
This fact is proved in Section \ref{unifb}.

We conclude this overview by
elaborating the analogy between this work and that of \cite{ESY} that was
mentioned in a remark after Theorem \ref{thmo}. The authors of
\cite{ESY} consider a Bose-Einstein condensate with a potential
between each pair of particles 
of the form $N^2 V \big( N(x_i - x_j) \big)$, these particles having
locations $x_i$ and $x_j$ in $3$-space. It is proved that the densities
matrices in the problem satisfy an infinite BBGKY hierarchy of
equations known as the Gross-Pitaevskii hierarchy. 
In our context, the assumption that the
\stoss holds in the naive form that $Q(z) - Q(0)$ is small
for small $z$ would lead to a derivation of the Smoluchowski PDE
(\ref{syspde}) with macroscopic reaction rates $\macp(n,m)$ being
incorrectly computed as being equal to the microscopic ones
$\micp(n,m)$. In \cite{ESY}, an analogously crude assumption that the
two-particle density matrices may be expressed as a product of one-particle
matrices yields the GP-hierarchy with an incorrect coupling constant
equal to that which arises by making use of the Born approximation. 
In each of these contexts,
it is necessary to understand the nature of the microscopic adjustment
to particle independence in the vicinity of each particle to determine
correctly the revelant macroscopic constants. The form of the
adjustment is very similar between the two problems. It would be of
much interest to determine why this is so.

\end{section}

\begin{section}{Establishing the \stoss}\label{stwo}

Recall that, for any given pair $(n,m) \in \mathbb{N}^2$,  $\ue =
\ue_{n,m}: \mathbb{R}^d \to [0,\infty)$ is the unique
function whose existence is ensured by Theorem \ref{th6.1} that lies in $C^2
\big( \mathbb{R}^d \big)$ satisfying
$\lim_{x \to \infty}{\ue_{n,m} \big( x \big)} = 0$
and
\begin{displaymath}
    \Delta \ue_{n,m} \big( x \big)  = \ald \big( n,m \big) \Big[
V_{\epsilon} \big(
    x \big) \ue_{n,m} \big( x \big) +  V^{\epsilon} \big( x \big) \Big],
\end{displaymath}
where the constants $\ald$ were defined in (\ref{aldo}).
We are using the notations
\begin{displaymath}
    V^{\epsilon}(x)  = \epsilon^{-d} V \Big( \frac{x}{\epsilon} \Big),\ \ \ \ 
    V_{\epsilon}(x)  = \epsilon^{-2} V \Big( \frac{x}{\epsilon} \Big).
\end{displaymath}
Note that $\ue_{n,m} : \mathbb{R}^d
\to [0,\infty)$ is given by
\setcounter{equation}{0}
\begin{equation}\label{formphi}
    \ue_{n,m} \big( x \big) = - c_0(d) \int_{\mathbb{R}^d}{\frac{
    V^{\epsilon}(y)}{\vert x - y \vert^{d-2}}  \ald \big( n,m \big)
\Big[ u_{n,m}
\Big( \frac{y}{\epsilon} \Big) + 1 \Big] dy},
\end{equation}
where $c_0=c_0(d)=(d-2)^{-1}\o_d^{-1}$ with $\o_d$ denoting the surface area of
the unit sphere $S^{d-1}$. We are denoting by $u_{n,m}$ the function obtained
by the choice $\epsilon
= 1$.
We present the conditions on the two test functions $J,\overline{J}:
\mathbb{R}^d \times \mathbb{N} \times [0,\i)
\to \R$ that appear in (\ref{exz}). It suffices to work with
functions that take
non-zero values for only one value in the second argument, such
functions measuring the presence of particles of a given mass. By a temporary
abuse of notation, we write
\begin{displaymath}
J(x,M_1,t)  = J(x,t) 1 \! \! 1  \big\{ m = M_1 \big\} \ \ \
\textrm{and} \ \ \ 
\overline{J}(x,M_1,t)  =  \overline{J}(x,t) 1 \! \! 1  \big\{ m =
M_2 \big\}, 
\end{displaymath}
where on the right-hand-side, $J$ and $\overline{J}$ denote smooth
maps from $\mathbb{R}^d \times [0,\infty)$ to $\R$ of compact support. We will
suppress the appearance of the $t$-variable when writing the arguments
of $J$ and $\overline{J}$.

Numerous terms arise when the operators $\mathbb{A}_0$ and
$\mathbb{A}_C$ act on the expression $X_z - X_0$ (recall that the
functions of configurations $X_z$, indexed by $z \in \mathbb{R}^d$,
were defined in (2.8) ). We now label these terms.
Unless stated otherwise, we will adopt a notation whereby all the
index labels appearing in sums should be taken to be distinct. This
includes the case of multiple sums. For example, $\sum_{k,l \in I_q}
\sum_{i \in I_q} f(x_k,x_l,x_i)$ denotes the sum of the evaluation of the
function $f$ over all arguments that are triples $(x_k,x_l,x_i)$ where
$k$, $l$ and $i$ are distinct indices in $I$. Note also that, unless
otherwise stated, whenever
the symbol $\ue$ appears in a summand, we mean $\ue_{M_1,M_2}$.

Firstly, we label those terms arising
from the action of the diffusion operator. To do so, note that, for a
time-dependent functional $F$ of the configuration space, this action
is given by
\begin{displaymath}
    \Big( \frac{\partial}{\partial t} +  \mathbb{A}_0 \Big) F =
\frac{\partial}{\partial t} F + \sum_{i \in I_q}{d \big(m_i)
    \Delta_{x_i} } F.
\end{displaymath}
Thus, we label as follows:
$$
\left(\frac{\partial}{\partial t}+\mathbb{A}_0\right) (X_z - X_0)( 
q(t) ) = H_{11} + H_{12} + H_{13} + H_{14}
+ H_2 + H_3 + H_4,
$$
with
\begin{eqnarray}
H_{11} & = & \epsilon^{2(d-2)} \sum_{i,j \in I_q}{\micp(m_i,m_j)
    \Big[ V^{\epsilon} \big( x_i - x_j + z  \big) -
V^{\epsilon} \big(   x_i - x_j  \big) \Big]} J(x_i,m_i,t)
    \overline{J}(x_j,m_j,t) \\
H_{12} & = & - \epsilon^{2(d-2)} \sum_{i,j \in I_q}{\micp(m_i,m_j)
    V_{\epsilon} \big( x_i - x_j  \big) \ue \big( x_i - x_j
    \big)  J(x_i,m_i,t) \overline{J}(x_j,m_j,t) } \nonumber \\
H_{13} & = &  \epsilon^{2(d-2)} \sum_{i,j \in I_q}{\micp(m_i,m_j)
     V_{\epsilon} \big( x_i - x_j + z \big) \ue \big( x_i - x_j+z
    \big)  J(x_i,m_i,t) \overline{J}(x_j,m_j,t) }, \nonumber \\
H_{14} & = &  \epsilon^{2(d-2)} \sum_{i,j \in I_q}{
     \Big[ \ue \big( x_i - x_j + z \big) -  \ue \big(
    x_i - x_j \big) \Big]} \nonumber \\
    & & \qquad \quad \Big[    J_t \big(x_i,m_i,t\big)
    \overline{J}(x_j,m_j,t) +   J \big(x_i,m_i,t\big)
      \overline{J}_t(x_j,m_j,t) \Big] , \nonumber
\end{eqnarray}
along with
\begin{eqnarray}
H_2 & =  & 2 \epsilon^{2(d-2)} \sum_{i,j \in I_q}{ d(m_i)
      \overline{J}(x_j,m_j,t)}
    \Big[  \ue_x (x_i - x_j + z ) - \ue_x
(x_i - x_j ) \Big]
    \cdot J_x(x_i,m_i,t), \\
H_3  & = &   2 \epsilon^{2(d-2)} \sum_{i,j \in I_q}{ d(m_j)
      J(x_i,m_i,t)}  \Big[  \ue_x (x_i - x_j + z ) -
     \ue_x (x_i - x_j ) \Big]
    \cdot \overline{J}_x (x_j,m_j,t) , \nonumber
\end{eqnarray}
and
\begin{eqnarray}
H_4 & = & \epsilon^{2(d-2)} \sum_{i,j \in I_q}{} \Big[
    \ue (x_i - x_j + z ) - \ue (x_i - x_j ) \Big]
\\
    & &  \qquad \qquad \quad  \Big[ d(m_i) \Delta_x J(x_i,m_i,t)
    \overline{J}(x_j,m_j,t) +  d(m_j) J(x_i,m_i,t)
    \Delta_x \overline{J}(x_j,m_j,t) \Big] , \nonumber
\end{eqnarray}
where $f_x$ denotes the gradient of $f$, and $\cdot$ the scalar
product. As for those terms arising from
the action of the collision operator,
\begin{displaymath}
\mathbb{A}_C (X_z - X_0)(q) = G_z (1) + G_z (2) - G_0 (1) - G_0 (2),
\end{displaymath}
where $G_z (1)$ is set equal to
\begin{eqnarray}
& & \frac 12\sum_{k,l \in
I_q}{\micp(m_k,m_l) V_{\epsilon}( x_k - x_l)} \epsilon^{2(d-2)} \sum_{i
\in I_q}{} \\
    & & \quad \bigg\{
    \frac{m_k}{m_k + m_l} \Big[  \ue ( x_k - x_i
    + z) J( x_k , m_k + m_l,t)
    \overline{J}(x_i,m_i,t) \nonumber \\
    & & \qquad \qquad \qquad  + \, \ue ( x_i - x_k
    + z) J( x_i , m_i ,t) \overline{J}(x_k,m_k + m_l,t)  \Big]
     \nonumber \\
    & & \quad \, \, + \, \frac{m_l}{m_k + m_l}  \Big[ \ue ( x_l - x_i
    + z) J( x_l , m_k + m_l,t)
    \overline{J}(x_i,m_i,t) \nonumber \\
    & &  \qquad \qquad \qquad  + \, \ue ( x_i - x_l
    + z) J( x_i , m_i ,t) \overline{J}(x_l,m_k + m_l,t) \Big]
     \nonumber \\
    & & \qquad  - \ \
    \Big[  \ue (x_k - x_i + z) J( x_k , m_k,t)
    \overline{J}(x_i,m_i,t) \nonumber \\
    & &  \qquad \qquad \qquad \qquad + \,  \ue (x_i - x_k + z)
    J( x_i , m_i ,t) \overline{J}(x_k,m_k,t) \Big]
     \nonumber \\
    & & \qquad - \ \
     \Big[  \ue (x_l - x_i + z) J( x_l , m_l,t)
    \overline{J}(x_i,m_i,t) \nonumber \\
    & & \qquad \qquad \qquad \qquad + \,  \ue (x_i - x_l + z) J( x_i ,
    m_i ,t) \overline{J}(x_l,m_l,t) \Big] \bigg\}
     \nonumber ,
\end{eqnarray}
and where
\begin{equation}\label{eqnzeta}
G_z (2)  =  - \epsilon^{2(d-2)} \sum_{k,l \in
I_q}{ \micp(m_k,m_l)  V_{\epsilon}(x_k - x_l)}
\ue (x_k -
x_l + z ) J(x_k,m_k,t)  \overline{J}(x_l,m_l,t).
\end{equation}
The expression in (\ref{geezee}) is an analogue of $G_0(2)$. The
terms in $G_z(1)$
arise from the changes in the functional $X_z$ when a collision
occurs due to the influence of the appearance and disppearance of
particles on other particles that are not directly involved. Those in
$G_z(2)$ are due to the absence after collision of the summand in
$X_z$ indexed by the colliding particles.

Note that
\begin{equation}\label{igt}
    H_{12} + G_0 \big( 2 \big) = 0.
\end{equation}
The process $\big\{ \big( X_z - X_0 \big)(t): t \geq 0 \big\}$ satisfies
\begin{eqnarray}\label{hact}
    \big( X_z - X_0 \big) \big( T \big) & = &
    \big( X_z - X_0 \big) \big( 0 \big) +
\int_{0}^{T}{ \Big( \frac{\partial}{\partial t} + \mathbb{A}_0 \Big)
\big(X_z - X_0\big) (t)dt} \\
    & & \qquad + \,
\int_{0}^{T}{ \mathbb{A}_C (X_z - X_0) (t)dt} \, + \, M(T), \nonumber
\end{eqnarray}
with $\big\{ M(t) : t \geq 0 \big\}$ being a martingale.
By using the labels for the various terms
that we just introduced, we find from (\ref{hact}) by use of (\ref{igt}) that
\begin{eqnarray}\label{abc}
& & \Big\vert
\int_{0}^{T}{ H_{11} \big(t\big) dt} +
\int_{0}^{T}{ H_{13} \big(t\big) dt} \Big\vert \\
& \leq &  \big\vert X_z - X_0 \big\vert \big( q(T) \big) + \big\vert
X_z - X_0 \big\vert \big( q(0) \big) \nonumber \\
& & + \int_{0}^{T}{\big\vert H_{14}  \big\vert (t) 
dt}+\int_{0}^{T}{\big\vert H_2  \big\vert (t) dt} +
\int_{0}^{T}{\big\vert  H_3  \big\vert (t) dt}  +
\int_{0}^{T}{\big\vert  H_4 \big\vert (t) dt} \nonumber \\
& & +
\int_{0}^{T}{  \big\vert G_z(1) - G_0(1) \big\vert (t) dt} +
\int_{0}^{T}{  \big\vert G_z(2) \big\vert (t) dt} \, + \, \big\vert
M(T) \big\vert. \nonumber
\end{eqnarray}
Since $J$ is of compact support, we have that $X_z(q(T))=0$ for
$T$ sufficiently large.

We denote by $\mathbb{P}_N$ the measure on functions from
$t \in [0,\infty)$ to the configurations determined by the process at
time $t$. Its expectation will be denoted $\mathbb{E}_N$.
We aim to prove the following estimates: for each $T > 0$,
\begin{eqnarray}\label{jayeye}
& & \int_{0}^{T}{ \mathbb{E}_N \vert H_{14} \vert (t) dt} \leq C
\vert z \vert^{\frac{2}{d+1}}, \\
& & \int_{0}^{T}{ \mathbb{E}_N \vert H_2 \vert (t) dt} \leq C \big\vert z
\big\vert^{\frac{1}{d + 1}},
\nonumber \\
& & \int_{0}^{T}{ \mathbb{E}_N \vert H_3 \vert (t) dt} \leq C \big\vert z
\big\vert^{\frac{1}{d + 1}}, \nonumber \\
& & \int_{0}^{T}{ \mathbb{E}_N \vert H_4 \vert (t) dt} \leq  C
\vert z \vert^{\frac{2}{d+1}}, \nonumber \\
& & \int_{0}^{T}{ \mathbb{E}_N \vert G_z(1) - G_0 (1) \vert (t) dt} \leq
C  \big\vert z \big\vert^{\frac{2}{d + 1}} , \nonumber \\
& & \int_{0}^{T}{ \mathbb{E}_N \vert G_z (2) \vert (t) dt} \leq C \Big(
\frac{\epsilon}{\vert z \vert} \Big)^{d-2},
\nonumber\\
& &\mathbb{E}_N|X_z-X_0(0)|\le C|z|. \nonumber
\end{eqnarray}
Later we apply the limit $\vert z \vert \to 0$  in such a way that
$\frac{\epsilon}{\vert z \vert} \to 0$. We will also show that, for
each $T \in (0,\infty)$,
\begin{equation}\label{martest}
\mathbb{E}_N \Big[ M \big( T \big)^2 \Big] \leq
C\e^{d-2}
  .
\end{equation}
\subsection{Lemmas bounding collision propensity}\label{stone}
We prove  three lemmas:
\begin{lemma}\label{lembc}
For any $T \in [0,\infty)$,
$$
\epsilon^{d-2} \mathbb{E}_N \int_0^T dt \sum_{i,j \in I_q}{\micp(m_i,m_j)}
    V_{\epsilon} \big( x_i - x_j \big) \leq 2Z.
$$
\end{lemma}
\begin{lemma}\label{lemtwop}
Let $J: \mathbb{R}^d \to [0,\infty)$ be continuous, and let $H:
\mathbb{R}^d \to [0,\infty)$ be a twice differentiable function such that
$ - \Delta H = J$.
Then we have the following inequalities,
\begin{displaymath}
\mathbb{E}_N \int_0^T \epsilon^{2(d-2)} \sum_{i,j \in I_q}{J(x_i - x_j)
m_i m_j (d(m_i) + d(m_j))} dt
\leq  \epsilon^{2(d-2)} \mathbb{E}_N \sum_{i,j \in I_{q(0)}}{H(x_i - x_j)
m_i m_j},
\end{displaymath}
\begin{displaymath}
\mathbb{E}_N \int_0^T\epsilon^{2(d-2)} \sum_{i,j \in I_q}{ V_{\epsilon} ( x_i
- x_j ) \micp(m_i,m_j) m_i m_j H(x_i - x_j)}
dt
\leq  \epsilon^{2(d-2)} \mathbb{E}_N  \sum_{i,j \in I_{q(0)}}{H(x_i - x_j)
m_i m_j}.
\end{displaymath}
\end{lemma}

\noindent
{\bf Remark} The fact that the dimension $d$ is at least three implies
that $H \geq 0$.
\begin{lemma}\label{lemthreep}
Assume that the function $\g:\mathbb{N}^2\to (0,\i)$ satisfies
\begin{equation}\label{hypo2}
n_2 \g \big( n_1, n_2 + n_3 \big) \max{ \Big\{ 1, \Big[ \frac{d (
n_2 + n_3 )}{d(n_2)}
\Big]^{\frac{3d - 2}{2}}\Big\}} \leq \big( n_2 + n_3 \big) \g(n_1,n_2).
\end{equation}
There exists a collection of constants $C: \mathbb{N}^2 \to (0,\infty)$,
such that, for any smooth function $J : \mathbb{R}^{2d} \to [0,\infty)$,
and any given $n_1,n_3 \in \mathbb{N}$,
\begin{eqnarray}
& &
\mathbb{E}_N{\int_0^T{dt \sum_{k,j,i \in I_{q(t)}}{\g (m_i,m_j)
V_{\epsilon} \big(
    x_i - x_j \big) J ( x_i , x_k ) 1 \! \! 1 \{ m_i = n_1, m_k =
n_3 \} }}} \\
    & \leq & \epsilon^{3(2-d)} C_{n_1,n_3}
    \sum_{n_2 \in \mathbb{N}}{\int{A^{\epsilon}_{n_1,n_2,n_3} \big( x_1,x_2,x_3
    \big)h_{n_1}(x_1)h_{n_2}(x_2)h_{n_3}(x_3) dx_1 dx_2 dx_3}}, \nonumber
\end{eqnarray}
where, also given $\epsilon > 0$ and $n_2 \in \mathbb{N}$, the function
$ A_{n_1,n_2,n_3}^{\epsilon}: \mathbb{R}^{3d} \to
[0,\infty) $ is defined by
\begin{eqnarray}
    & & \big( d(n_1 )\Delta_{x_1} + d(n_2) \Delta_{x_2} + d(n_3)
\Delta_{x_3} \big)
    A_{n_1,n_2,n_3}^{\epsilon} (x_1,x_2,x_3) \\
    & & \qquad = -   \g (n_1,n_2)
    V_{\epsilon} ( x_1 - x_2 ) J ( x_1 ,
    x_3 ) .  \nonumber
\end{eqnarray}
\end{lemma}
{\bf Proof of Lemma \ref{lembc}} Let $X(t)$ denote the total number of
particles present in the model at time $t$. Then, for any $T > 0$,
$$
\mathbb{E}_NX (T) = \mathbb{E}_NX (0) + \int_0^T{}
\mathbb{E}_N \mathbb{A}_0 X (t)dt +  \int_0^T{}
\mathbb{E}_N \mathbb{A}_C X (t)dt.
$$
   From $A_0 X = 0$, it follows that
$$
    - \int_0^T{} \mathbb{E}_N A_C X (t)dt =  \mathbb{E}_N \int_0^T dt
\sum_{i,j \in I_q}{\micp(m_i,m_j)} V_{\epsilon} (
x_i - x_j )  \leq \mathbb{E}_NX(0) = N.
$$
The result follows from the fact that $N = Z \epsilon^{2-d}$. $\Box$ \\

{\bf Proof of Lemma \ref{lemtwop}}
Set
$$
X_q = \epsilon^{2(d-2)} \sum_{i,j \in I_q}{H(x_i - x_j)m_i m_j}.
$$
Recall the mechanism of the dynamics at collision:  the location of
the newly created particle is one of the
two locations of the colliding particles, with weights proportional to
the masses of the incident particles. We see that when $\mathbb{A}_C$ acts on
$X_q$, all those terms indexed by pairs of particles one of which is
not involved in the collision cancel. Thus,
\begin{equation}\label{locc}
\mathbb{A}_C X = - \epsilon^{2(d-2)} \sum_{i,j \in I_q}{ V_{\epsilon} ( x_i
- x_j ) \micp(m_i,m_j) m_i m_j H(x_i - x_j)}.
\end{equation}
Note also that
\begin{equation}\label{lokc}
\mathbb{A}_0 X =  \epsilon^{2(d-2)} \sum_{i,j \in I_q}{ \Delta  H (x_i
- x_j) m_i m_j (d(m_i) + d(m_j))}.
\end{equation}
   From the non-positivity of $\mathbb{A}_C X$, apparent from (\ref{locc}),
and the non-negativity of $X$, follows
$$
- \mathbb{E}_N \int_0^T \mathbb{A}_C X(t) dt
- \mathbb{E}_N \int_0^T \mathbb{A}_0 X(t) dt \leq \mathbb{E}_N X(0).
$$
Substituting the form for $\mathbb{A}_0 X$ and $\mathbb{A}_C X$
in (\ref{lokc}) yields the
result. $\Box$\\

{\bf Proof of Lemma \ref{lemthreep}} Note that
\begin{eqnarray}\label{estar}
    & & A_{n_1,n_2,n_3}^{\epsilon}(x_1,x_2,x_3) \\
    & = &
    c_0 \big( 3d \big) \int_{\mathbb{R}^d}  \int_{\mathbb{R}^d}
\int_{\mathbb{R}^d}
     \bigg( \frac{\vert x_1 - z \vert^2}{d(n_1)}  + \frac{\vert x_2 -
y \vert^2}{d(n_2)} + \frac{\vert x_3 - y' \vert^2}{d(n_3)}
\bigg)^{\frac{-3d + 2}{2}} \g (n_1,n_2) \nonumber \\
    & & \qquad \qquad \qquad \quad J ( z , y' )  V_{\epsilon} ( z - y )
dz dy dy', \nonumber
\end{eqnarray}
where $c_0$ was defined right after (3.1). Define
$$
X_{n_1,n_2,n_3}(q) = \sum_{i,j,k \in
I_q}{A_{n_1,n_2,n_3}^{\epsilon}(x_i,x_j,x_k)
    1 \! \! 1 \{ m_i =n_1, m_j= n_2, m_k = n_3 \} },
$$
and
$$
Y_{n_1,n_3}(q) = \sum_{n_2 \in \mathbb{N}}{X_{n_1,n_2,n_3}(q)}.
$$
Note that $\mathbb{A}_0 Y_{n_1,n_3}$ is given by
\begin{eqnarray}\label{cross}
& & \sum_{n_2 \in \mathbb{N}} \sum_{i,j,k \in I_q}{}\big( d(n_1)
\Delta_{x_i} + d(n_2) \Delta_{x_j}
+ d(n_3) \Delta_{x_k} \big) A_{n_1,n_2,n_3}^{\epsilon}(x_i,x_j,x_k)
\\
& & \qquad \qquad \qquad \qquad
1 \! \! 1 \{ m_i=n_1, m_j= n_2, m_k = n_3 \}  \nonumber \\
    & = & -  \sum_{n_2 \in \mathbb{N}} \sum_{i,j,k \in I_q}{\g (n_1,n_2)
    V_{\epsilon} ( x_i - x_j )
    J( x_i , x_k )
1 \! \! 1 \{ m_i=n_1,m_j=n_2,m_k=n_3 \} }. \nonumber
\end{eqnarray}
Note also that
$$
\mathbb{A}_C Y_{n_1,n_3} = I_{n_1,n_3} + G_{n_1,n_3} -  L_{n_1,n_3},
$$
where $I_{n_1,n_3}$ is given by
\begin{eqnarray}
    & &  \frac 12\sum_{k,l \in I_q}{} V_{\epsilon} ( x_k
- x_l ) \micp(m_k,m_l) \\
     & & \ \  \sum_{i,j \in I_q}{}\Big[ \frac{m_k}{m_k + m_l} A_{n_1,m_k +
     m_l,n_3}^{\epsilon}(x_i,x_k,x_j) +  \frac{m_l}{m_k + m_l} A_{n_1,m_k +
     m_l,n_3}^{\epsilon}(x_i,x_l,x_j) \nonumber \\
     & & \qquad \quad  -
     A_{n_1,m_k,n_3}^{\epsilon}(x_i,x_k,x_j)  -
A_{n_1,m_l,n_3}^{\epsilon}(x_i,x_l,x_j) \Big]
    1 \! \! 1 \{ m_i=n_1,m_j = n_3 \}. \nonumber
\end{eqnarray}
The term $G_{n_1,n_3}$ is set equal to
\begin{eqnarray}
     & &\frac 12 \sum_{M_1=1}^{n_1-1}{} \sum_{k,l \in I_q}{ V_{\epsilon} ( x_k
- x_l ) \micp(m_k,m_l) 1 \! \! 1 \{ m_k = M_1 , m_l = n_1 - M_1 \} }
\\
& & \qquad \sum_{n_2 \in \mathbb{N}}{} \sum_{i,j \in I_q}{} \Big[
\frac{m_k}{m_k + m_l}  A_{n_1,n_2,n_3}^{\epsilon}(x_k,x_i,x_j)
\nonumber \\
    & & \qquad \qquad \qquad +
\frac{m_l}{m_k + m_l}  A_{n_1,n_2,n_3}^{\epsilon}(x_l,x_i,x_j) \Big]
1 \! \! 1 \{ m_i=n_2,m_j = n_3 \}   \nonumber  \\
& + & \frac 12 \sum_{M_3=1}^{n_3-1}{} \sum_{k,l \in I_q}{ V_{\epsilon} ( x_k
- x_l ) \micp(m_k,m_l) 1 \! \! 1 \{ m_k = M_3 , m_l = n_3 - M_3 \} }
\nonumber \\
& & \qquad \sum_{n_2 \in \mathbb{N}}{} \sum_{i,j \in I_q}{} \Big[
\frac{m_k}{m_k + m_l}  A_{n_1,n_2,n_3}^{\epsilon}(x_i,x_j,x_k)
\nonumber \\
    & & \qquad \qquad \qquad +
\frac{m_l}{m_k + m_l}  A_{n_1,n_2,n_3}^{\epsilon}(x_i,x_j,x_l)  \Big]
1 \! \! 1 \{ m_i=n_1,m_j = n_2 \}.   \nonumber
\end{eqnarray}
The term $L_{n_1,n_3}$ may be written as a sum of $L^1_{n_1,n_3}$ and
$L^2_{n_1,n_3}$, where  $L^1_{n_1,n_3}$ is given by
\begin{eqnarray}
& &   \frac 12 \sum_{n_2 \in \mathbb{N}}{} \sum_{k,l \in I_q}{ V_{\epsilon} ( x_k
- x_l ) \micp(m_k,m_l)} \\
    & & \quad \Big[ \sum_{i \in I_q}{A_{n_1,n_2,n_3}^{\epsilon}(x_i,x_k,x_l)
    1 \! \! 1 \{ m_i=n_1,m_k=n_2,m_l = n_3 \} } \nonumber \\
    & & \qquad +  \sum_{i \in I_q}{A_{n_1,n_2,n_3}^{\epsilon}(x_k,x_i,x_l)
    1 \! \! 1 \{ m_k=n_1,m_i=n_2,m_l = n_3 \} }  \nonumber \\
    & & \qquad +  \sum_{i \in I_q}{A_{n_1,n_2,n_3}^{\epsilon}(x_k,x_l,x_i)
    1 \! \! 1 \{ m_k=n_1,m_l=n_2,m_i = n_3 \} }  \Big], \nonumber
\end{eqnarray}
and where  $L^2_{n_1,n_3}$ is given by
\begin{eqnarray}
& &   \frac 12 \sum_{n_2 \in \mathbb{N}}{} \sum_{k,l \in I_q}{ V_{\epsilon} ( x_k
- x_l ) \micp(m_k,m_l)} \\
    & & \quad \Big[
    \sum_{i,j \in I_q}{A_{n_1,n_2,n_3}^{\epsilon}(x_i,x_j,x_k)
    1 \! \! 1 \{ m_i=n_1,m_j=n_2,m_k = n_3 \} } \nonumber \\
    & & \qquad \ + \sum_{i,j \in I_q}{A_{n_1,n_2,n_3}^{\epsilon}(x_k,x_i,x_j)
    1 \! \! 1 \{ m_k=n_1,m_i=n_2,m_j = n_3 \} } \Big]. \nonumber
\end{eqnarray}
The term $L^1_{n_1,n_3}$ arises due to the loss at collision of
particles indexed by $k$ and $l$ in the contribution to $Y_{n_1,n_3}$
of summands where two of the three arguments of $A^{\epsilon}$ are
occupied by $x_k$ and $x_l$. The term $L^2_{n_1,n_3}$ carries those
summands where only one of the pair of particles, $x_k$ and $x_l$,
appears in the arguments of $A^{\epsilon}$. The case where only one of
the pair appears, and does so in the unrestricted second argument of
$A^{\epsilon}$, is not included here, but instead in the `internal'
term, $I_{n_1,n_3}$.

By the hypothesis (\ref{hypo2}), we
see from the form (\ref{estar}) of $A^{\epsilon}$ that $I_{n_1,n_3} \leq
    0$.
   From (\ref{cross}), the lemma will follow provided that we establish the
following assertion. \\
{\bf Claim}
There exist a collection of constants $C: \mathbb{N}^2 \to [0,\infty)$,
such that, for each pair $(n_1,n_3) \in \mathbb{N}^2$, and for all
smooth $J: \mathbb{R}^d \to
[0,\infty)$, and $\epsilon, T > 0$,
$$
 \mathbb{E}_N \int_0^T {}L_{n_1,n_3}(t) dt\le R_{n_1,n_3}, \ \ \ \ 
 - \mathbb{E}_N \int_0^T {} \mathbb{A}_0 Y_{n_1,n_3}(t) dt\le R_{n_1,n_3}, 
$$
where 
$$
R_{n_1,n_3}=    C_{n_1,n_3}
    \epsilon^{3(2-d)}
    \sum_{n_2 \in \mathbb{N}}{\int{A^{\epsilon}_{n_1,n_2,n_3} \big( x_1,x_2,x_3
    \big)h_{n_1}(x_1)h_{n_2}(x_2)h_{n_3}(x_3) dx_1 dx_2 dx_3}}.
$$
 We establish the claim by an induction on the parameters $n_1,n_3
\in \mathbb{N}$. We proceed to the generic step, at which the claim
may be assumed for parameter values $(M_1,M_3)$, where
$M_1 < n_1$, or $M_1 = n_1, M_3 < n_3$. This step includes the
initial case where $n_1=n_3=1$.
From
$$
\mathbb{E}_NY_{n_1,n_3}(T) = \mathbb{E}_NY_{n_1,n_3}(0) + \int_0^T{}
\mathbb{E}_N(A_0 Y_{n_1,n_3})(t)dt +  \int_0^T{}
\mathbb{E}_N(A_C Y_{n_1,n_3})(t)dt,
$$
and the negativity of the terms $- \mathbb{E}_N \int_0^T{}
L_{n_1,n_3}(t) dt $ and $ \mathbb{E}_N \int_0^T{}
I_{n_1,n_3}(t) dt $, follows
\begin{equation}\label{saq}
     - \int_0^T{}\mathbb{E}_N(A_0 Y_{n_1,n_3})(t)dt  \leq
     \mathbb{E}_NY_{n_1,n_3}(0) +  \mathbb{E}_N \int_0^T{} G_{n_1,n_3}(t) dt.
\end{equation}
Note that
\begin{displaymath}
\mathbb{E}_N \int_0^T{} G_{n_1,n_3}(t) dt
    \leq C_{n_1} \sum_{M_1 = 1}^{n_1
    -1}{\mathbb{E}_N\int_0^T{L_{M_1,n_3}(t)dt}} \ + \   C_{n_3} 
\sum_{M_3 = 1}^{n_3
    -1}{\mathbb{E}_N\int_0^T{L_{n_1,M_3}(t)dt}} ,
\end{displaymath}
where here we have used such facts as
$$
A_{n_1,n_2,n_3}^{\epsilon} \leq \max\Big\{
\frac{d(n_1)}{d(j)},\frac{d(j)}{d(n_1)}
\Big\}^{\frac{3d - 2}{2}} A_{j,n_2,n_3}^{\epsilon} \ \textrm{for $j \in \{ 1,
\ldots, n_1 - 1 \}$.}
$$
We do not require any hypotheses on the parameters in this last step,
since there are only finitely many possible values of mass for
particles that coagulate to form a new particle of some given mass,
and we can take care of each of these cases by choosing the value of
the constants $C_{n_1,n_3}$.
   From (\ref{saq}), the inductive hypotheses, the elementary inequalities
$$
A_{j,n_2,n_3}^{\epsilon} \leq \max\Big\{
\frac{d(n_1)}{d(j)},\frac{d(j)}{d(n_1)}
\Big\}^{\frac{3d - 2}{2}} A_{n_1,n_2,n_3}^{\epsilon} \ \textrm{for $j \in \{ 1,
\ldots, n_1 - 1 \}$,}
$$
$$
A_{n_1,n_2,j}^{\epsilon} \leq \max\Big\{
\frac{d(n_3)}{d(j)},\frac{d(j)}{d(n_3)}
\Big\}^{\frac{3d - 2}{2}} A_{n_1,n_2,n_3}^{\epsilon} \ \textrm{for $j \in \{ 1,
\ldots, n_3 - 1 \}$,}
$$
and the fact that the
particles are initially independent, we
find that
the expression
$$- \int_0^T{}\mathbb{E}_N(A_0 Y_{n_1,n_3})(t)dt,$$
is bounded above by
$$
     C_{n_1,n_3}
    \epsilon^{3(2-d)}
    \sum_{n_2 \in \mathbb{N}}{\int{A^{\epsilon}_{n_1,n_2,n_3} \big( x_1,x_2,x_3
    \big)h_{n_1}(x_1)h_{n_2}(x_2)h_{n_3}(x_3) dx_1 dx_2 dx_3}},
$$
which is the second part of the claim. To establish the first part,
note that
$$
     \int_0^T{}\mathbb{E}_N(L_{n_1,n_3})(t)dt  \leq
     \mathbb{E}_NY_{n_1,n_3}(0) +  \mathbb{E}_N \int_0^T{} G_{n_1,n_3}(t) dt.
$$
Hence, the term to which the first part of the claim refers satisfies
the same upper bound as the one that was just bounded. $\Box$

The following corollary to Lemma 3.2 will be used in Section 5.
\begin{lemma}\label{lem3.4}
There exists a constant $C$ such that for any $T \in [0,\infty)$,
$$
\epsilon^{d-2} \mathbb{E}_N \int_0^T dt \sum_{i,j \in I_q}{m_im_j\micp(m_i,m_j)}
    V_{\epsilon} \big( x_i - x_j \big) \leq C.
$$
\end{lemma}

{\bf Proof} On account of Lemma 3.2, it suffices to find a pair of
functions $(J,H)=(J^\e,H^{\epsilon})$ 
and a positive constant $\delta $ such that $J\ge 0$, $-\Delta J=H$,
$\epsilon^{d-2}H^{\epsilon}(x)\ge \delta $ whenever $V_\epsilon (x)\neq 0$ and
that the expression
\begin{displaymath} 
X_N=\epsilon^{2(d-2)}\mathbb{E}_N \sum_{i,j \in I_{q(0)}}{H^\epsilon (x_i - x_j)
m_i m_j}
\end{displaymath}
is bounded in $N$. For this, let us pick a nonnegative smooth function $W$ of
compact support such that $W(x)\ge 1$ whenever $|x|\le 2\eta $ where $\eta  $ is chosen 
so that the support of $V$ is contained in the ball $B_{\eta  }(0)$. We then set
$J^\epsilon(x)=\epsilon^{-d}W(x/\epsilon)$ and
$H^\epsilon(x)=\epsilon^{2-d}a(x/\epsilon)$ where the function $a$ solves $-\Delta
a=W$. More precisely,
$
a(x)=c_0(d)\int W(x-y)|y|^{2-d}dy$. By restricting the domain of 
$y$-integration to the case $|y|\le \eta  $, it is not hard to deduce that if $|x|\le \eta  $,
then $a(x)\ge (d-2)^{-1}\eta  $. Also, $a(x)\le \bar c |x|^{2-d}$ for a constant $\bar c$.
The proof of this is very similar to the proof of the first part of
Lemma 3.5 below and omitted. The boundedness
of the sequence $\{X_N\}$  is now a straightforward consequence of our assumption 
$\bar k*(-\Delta k)^{-1}(0)<\i$ on the initial density.   $\Box$

\subsection{Bounds on functionals of $u_{n,m}$}\label{sttwo}
We will verify the assertions presented in (\ref{jayeye}). The
following lemma provides the bounds on the behaviour of the functions
$\big\{ u^\epsilon_{n,m}: \mathbb{R}^d \to [0,\infty): (n,m) \in
\mathbb{N} \big\}$ and other functions that
will be used in this section.
Recall that $u^\e_{n,m}(x)=u^\e(x)=\e^{2-d}u(x/\e)$ where $u$ satisfies
(1.8). We choose the constant $C_0$ so that $V(x)=0$ whenever $|x|\ge C_0$.
Recall that $k=\sum_n nh_n$.
\begin{lemma}\label{phiest} There exists a collection of constants $C:
\mathbb{N}^2 \to (0,\infty)$ for which
the following bounds hold. \\
\begin{itemize}
\item for $x \in \mathbb{R}^d$,
$-u_{n,m}(x) \leq \frac{C_{n,m}}{\vert x \vert^{d-2}}$ and
$\big\vert \nabla u_{n,m} (x) \big\vert \leq \frac{C_{n,m}}{\vert x
\vert^{d-1}}$.
\item for $x \in \mathbb{R}^d$ satisfying $\vert x \vert \geq  \max \big\{ 2
\vert z \vert + C_0 \epsilon , 2 C_0 \epsilon \big\}$,
\begin{equation}\label{lept}
    \Big\vert u^\epsilon_{n,m} \big( x + z \big) - u^\epsilon_{n,m} \big( x
    \big) \Big\vert \leq \frac{C_{n,m} \vert z \vert}{\vert x \vert^{d-1}}
\end{equation}
and
\begin{equation}\label{lepta}
    \Big\vert \nabla u^\epsilon_{n,m} \big( x + z \big) - \nabla
u^\epsilon_{n,m} \big( x
    \big) \Big\vert \leq \frac{C_{n,m} \vert z \vert}{\vert x \vert^d}.
\end{equation}
\item letting $H=H_{n,m}: \mathbb{R}^d \to [0,\infty)$ be given by
\begin{displaymath}
    - \Delta H_{n,m} (x)  = u^{\epsilon}_{n,m} ( x + z) 1 \! \! 1 \big\{ \vert x
    \vert \leq \rho \big\},
\end{displaymath}
we have
\begin{displaymath}\label{lepti}
    \int{H_{n,m}(x_1-x_2)k(x_1)k(x_2)dx_1 dx_2} \leq C_{n,m}
    \big( \rho + \vert z \vert \big)^2.
\end{displaymath}
\item for $z \in \mathbb{R}^d$, let $\hat{H}_z \big( =
\hat{H}_{z,n,m} \big): \mathbb{R}^d \to
[0,\infty)$ be given by
\begin{displaymath}
    - \Delta \hat{H}_z (x) = \Big\vert \nabla u_{n,m}^\epsilon \big(
    x+z\big) \Big\vert 1 \! \! 1 \big\{ \vert x \vert \leq \rho \big\}.
\end{displaymath}
Then,
\begin{equation}\label{lepto}
\int{ \hat{H}_z (x_1-x_2)k(x_1)k(x_2)dx_1
dx_2} \leq C_{n,m} \big( \rho + \vert z
\vert \big)  .
\end{equation}
\item for any positive integers $n$ and $m$ and a smooth function
$\overline J$ of compact support, there exists
a constant $C_{n,m}(\overline J)$ such that for any given $z \in
\mathbb{R}^d$, the function
$A_{n_1,n_2,n_3}^{\epsilon}: \mathbb{R}^{3d} \to [0,\infty)$ defined
by
\begin{eqnarray}
    & & \big( d(n_1 )\Delta_{x_1} + d(n_2) \Delta_{x_2} + d(n_3)
\Delta_{x_3} \big)
    A_{n_1,n_2,n_3}^{\epsilon} (x_1,x_2,x_3) \\
    & & \qquad = - \gamma(n_1,n_2)  u_{n,m}^\e(x_1 - x_3 + z) \e^{-2} V
    \Big( \frac{x_1 - x_2}{\epsilon} \Big) 1\!\!1 ( |x_1-x_3|\le \rho )
\overline{J}(x_3) \nonumber
\end{eqnarray}
satisfies
\begin{eqnarray}\label{leptu}
    & &\sum_{n_1,n_2,n_3}\int{ A_{n_1,n_2,n_3}^{\epsilon} \big( 
x_1,x_2,x_3 \big)
    h_{n_1}(x_1)h_{n_2}(x_2)h_{n_3}(x_3) d x_1
    d x_2 d x_3  } \\
    & & \qquad \qquad\leq C_{n,m}(\overline J) 
\epsilon^{d-2}(\rho+|z|)^2. \nonumber
\end{eqnarray}
\end{itemize}
\end{lemma}
{\bf Proof}
Throughout the proof, we write $u$ for the function $u_{n,m}$ and
$\ald$ for the constant $\ald(n,m)$ of (\ref{aldo}). The dependence of the
constants on $n$ and $m$ arises from that of $\ald$, and is also omitted.
To prove the first part of the Lemma, we will show that
\begin{equation}\label{ussoon}
    0 \leq -u ( x )  \leq c \min \Big\{ \vert x \vert^{2-d} , 1 \Big\},\ \ \ \
\
    0 \leq |\nabla u ( x )|  \leq c \min \Big\{ \vert x \vert^{1-d} , 1 \Big\}.
\end{equation}
Note firstly that
\begin{displaymath}
    u \big( x \big) = - c_0 \ald\int_{\mathbb{R}^d}{\frac{V(y) }{\big\vert x -
    y \big\vert^{d-2}} \Big[ u \big( y \big) + 1 \Big] dy }
\end{displaymath}
By Theorem 6.1,  $u(x) \in [-1,0]$ for all $x \in \mathbb{R}^d$. From this we
find that
\begin{displaymath}
    -u \big( x \big) \leq c_0 \ald
\int_{\mathbb{R}^d}{\frac{V(y)}{\big\vert x -
    y \big\vert^{d-2}} dy}
\end{displaymath}
If $\vert x \vert$ is large, then $\vert x - y \vert \geq \vert x
\vert/2$ whenever $V (y) \not= 0$, so that
\begin{displaymath}
     \int_{\mathbb{R}^d}{\frac{V(y)}{\vert x - y \vert^{d-2}} dy} \leq
     2^{d-2}   \int_{\mathbb{R}^d}{\frac{V(y)}{\vert x  \vert^{d-2}} dy} =
     2^{d-2} \vert x \vert^{2-d}.
\end{displaymath}
If $\vert x \vert$ is small, then
\begin{displaymath}
     \int_{\mathbb{R}^d}{\frac{V(y)}{\vert x - y \vert^{d-2}} dy} \leq
     \int_{\vert x-y \vert \leq C}{\frac{V(y)}{\vert x - y \vert^{d-2}} dy} \leq
     \int_{0}^{C}{\rho d \rho} = \frac{C^2}{2},
\end{displaymath}
proving the first inequlity in (\ref{ussoon}).
In the same way we can
deduce the second inequality in (\ref{ussoon}). Indeed
\begin{displaymath}
    \nabla u(x) = c_0 \ald(2-d)\int_{\mathbb{R}^d}{ V(y)
     \frac{x-y}{\big\vert x - y \big\vert^d} \Big[  u \big(y\big) +
     1 \Big]} dy.
\end{displaymath}
If $\vert x \vert$ is large, then $\big\vert x - y \big\vert \geq
\vert x \vert / 2$, and
\begin{displaymath}
    \bigg\vert \frac{x - y}{\vert x - y \vert^d} \bigg\vert =
    \frac{1}{\vert x - y \vert^{d-1}} \le c \vert x \vert^{1-d},
\end{displaymath}
whereas, if $\vert x \vert$ is small, then
\begin{displaymath}
    \big\vert \nabla u (x) \big\vert \leq c \int_{\vert y \vert \leq
    c_1}{\frac{1}{\vert x - y \vert^{d-1}} dy} \leq c
    \int_0^{c_1}{\rho^{d-1} \rho^{1 - d} d \rho} \leq c_2.
\end{displaymath}
Thus,
\begin{displaymath}
\big\vert \nabla u (x) \big\vert \leq c_0 \min \Big\{ 1, \vert x
\vert^{1-d} \Big\},
\end{displaymath}
and so
$ \big\vert \nabla u (x) \big\vert \leq c \vert x \vert^{1-d}$.

To prove (\ref{lept}) in the second part of the lemma, note that
\begin{eqnarray}
    & & \Big\vert u^\epsilon \big( x + z \big) - u^\epsilon \big( x
    \big) \Big\vert
\\
    & \leq & C  \int_{\mathbb{R}^d} \epsilon^{-d} V \Big(
    \frac{y}{\epsilon}\Big) \bigg\vert \frac{1}{\vert x + z  - y \vert^{d-2}}
    - \frac{1}{\vert x -y \vert^{d-2}} \bigg\vert  dy  \nonumber \\
    & \leq & C \int_{\mathbb{R}^d} \epsilon^{-d} V \Big(
    \frac{y}{\epsilon}\Big) \left| \frac{\Big\vert  \vert x - y
    \vert^{d-2} -
    \vert x + z - y \vert^{d-2} \Big\vert}{\vert x + z  - y \vert^{d-2}
    \vert x -y \vert^{d-2}} \right| dy \nonumber
\end{eqnarray}
Note that
\begin{equation}\label{chon}
    \big\vert x + z - y \big\vert^{d-2} - \big\vert x - y \big\vert^{d-2}
    \leq  \Big( \big\vert x - y \big\vert  + \vert z \vert \Big)^{d-2} -
    \vert x - y \vert^{d-2},
\end{equation}
and that
\begin{equation}\label{chtwo}
    \big\vert x  - y \big\vert^{d-2} - \big\vert x + z - y \big\vert^{d-2}
    \leq  \vert x - y \vert^{d-2} -  \Big( \big\vert x - y \big\vert   - \vert z
    \vert \Big)^{d-2} .
\end{equation}
The right-hand-sides of (\ref{chon}) and (\ref{chtwo}) take the form $c
\alpha^{d-3} \vert z \vert$, for some $\alpha \in \big[ \vert x-y
\vert - \vert z \vert, \vert x- y \vert + \vert z \vert \big]$. Note that
if $V(y/\e)\neq 0$, then
$$
\vert x - y \vert - \vert z \vert \geq \vert x \vert - C_0 \epsilon -
    \vert z \vert\ge 0,
    $$
    by our assumption on $|x|$. As a result, we have that
$\vert x - y \vert + \vert z \vert \leq 2 \vert x-y \vert$,
so that
\begin{displaymath}
    \Big\vert \big\vert x - y \big\vert^{d-2} - \big\vert x + z - y
    \big\vert^{d-2} \Big\vert \leq C \vert z \vert \vert x-y \vert^{d-3}.
\end{displaymath}
Hence,
\begin{displaymath}
    \Big\vert u^\epsilon \big( x + z \big) - u^\epsilon \big( x
    \big) \Big\vert  \leq
     C \vert z \vert  \int_{\mathbb{R}^d}{ V^{\epsilon} \big(
    y \big) \frac{\vert x -y  \vert^{-1}}{\vert x + z - y
    \vert^{d-2}} dy } .
\end{displaymath}
   From the inequality $\vert x \vert \geq \max \big\{ 2 \vert z \vert +
C_0 \epsilon , 2 C_0 \epsilon \big\}$, we deduce that  $\big\vert x +
z - y \big\vert \geq \vert x -y \vert/2$ and $\vert x - y \vert \geq
\vert x \vert/2$. We obtain (\ref{lept}).

In seeking to prove (\ref{lepta}), note that
\begin{displaymath}
    \frac{x + z - y}{\big\vert x + z - y \big\vert^d} - \frac{x -
    y}{\big\vert x - y \big\vert^d} = \frac{ \big(x + z - y \big)
    \big\vert x - y \big\vert^d - \big( x - y \big) \big\vert x + z - y
    \big\vert^d}{ \big\vert x + z - y \big\vert^d \big\vert x - y \big\vert^d }.
\end{displaymath}
Note that, for any $a \in \mathbb{R}^d$,
\begin{equation}\label{cule}
    \Big\vert \big( a + z \big) \vert a \vert^d - \vert a + z \vert^d a
    \Big\vert \leq  \vert a \vert \Big\vert \vert a + z \vert^d -
    \vert a \vert^d \Big\vert + \vert z \vert \vert a \vert^d
     \leq c_1 \vert z \vert \vert a \vert^d,
\end{equation}
so long as $|z|\le|a|$. Given that
\begin{displaymath}
    \nabla \ue(x) = c_0 (2-d)\ald\int_{\mathbb{R}^d}{ V^{\epsilon}(y)
     \frac{x-y}{\big\vert x - y \big\vert^d} \Big[  u \big( y \big)
     + 1  \Big]} dy,
\end{displaymath}
we may apply (\ref{cule}) with the choice $a = x - y$ to obtain
\begin{displaymath}
     \Big\vert  \nabla\ue ( x + z ) - \nabla\ue (x)
     \Big\vert \leq C \vert z \vert \int_{\mathbb{R}^d}{
     \frac{V^{\e}(y)}{\big\vert x + z - y \big\vert^d}} dy.
\end{displaymath}
   From the inequality $\vert x \vert \geq \max \big\{ 2 \vert z \vert +
C_0 \epsilon , 2 C_0 \epsilon \big\}$, we deduce that  $\big\vert x +
z - y \big\vert \geq \vert x -y \vert/2$ and $\vert x - y \vert \geq
\vert x \vert/2$. We conclude that
\begin{displaymath}
     \Big\vert  \nabla\ue ( x + z ) - \nabla\ue (x)
     \Big\vert \leq  \frac{C \vert z \vert}{\vert x \vert^d} ,
\end{displaymath}
as required.

To prove the third part of the Lemma, note that
\begin{displaymath}
    |H (x)| \leq c_1 \int_{\mathbb{R}^d}{\big\vert x - y \big\vert^{2-d}
\big\vert y + z \big\vert^{2-d} 1
    \! \! 1 \big\{ \vert y \vert \leq \rho \big\} }
\end{displaymath}
by the first part of the Lemma. Hence,
\begin{eqnarray}
    & & \int{|H(x_1-x_2)|k(x_1)k(x_2)dx_1dx_2} \\
& \leq & c_1 \int_{\vert y \vert \leq
    \rho}{ \big\vert y + z  \big\vert^{2-d} \int{\vert x_1-x_2 - y
\vert^{2-d}k(x_1)k(x_2)dx_1dx_2}dy} \nonumber \\
    & = & c_2  \int_{\vert y \vert \leq
    \rho}{ \big\vert y + z  \big\vert^{2-d} \left(\bar k*(-\Delta)^{-1}k\right)
(-y)dy} \nonumber \\
    & \le & c_3    \int_{\vert y \vert \leq
    \rho}{ \big\vert y + z  \big\vert^{2-d} dy} 
    \leq  c_4  \big( \rho + \vert z \vert \big)^2, \nonumber 
\end{eqnarray}
where we used our assumption $\bar k*(-\Delta)^{-1}k\in L^\i_{loc}$ 
for the third inequality.
This establishes the third part of the Lemma.

To prove the fourth part of the Lemma, we firstly find a bound for
$\nabla\ue  (x) = \epsilon^{1-d} \nabla u \big( x/\epsilon
\big)$. Using the first part of the Lemma,
\begin{displaymath}
    \hat{H}_z (x) \leq c \int_{\mathbb{R}^d}{\big\vert x - y \big\vert^{2-d}
\big\vert y + z \big\vert^{1-d} 1
    \! \! 1 \big\{ \vert y \vert \leq \rho \big\} },
\end{displaymath}
it follows that
\begin{eqnarray}
    & &  \int{|\hat H_z(x_1-x_2)|k(x_1)k(x_2)dx_1dx_2} \\
& \leq & c_1 \int_{\vert y \vert \leq
    \rho}{ \big\vert y + z  \big\vert^{1-d} \int{\vert x_1-x_2 - y
\vert^{2-d}k(x_1)k(x_2)dx_1dx_2}dy} \nonumber \\
     & \le & c_2    \int_{\vert y \vert \leq
    \rho}{ \big\vert y + z  \big\vert^{1-d} dy} 
     \leq  c_3  \big( \rho + \vert z \vert \big).     \nonumber
\end{eqnarray}

We have deduced (\ref{lepto}).

As for the fifth part of the Lemma, let us write $J(a)$ for
$1\!\!1(|a|\le \rho)$ and define the quantity $I$
to be 
    \begin{eqnarray}
     &  & c_0(3d)\sum_{n_1,n_2,n_3} \int_{} dx_1 dx_2 dx_3
    \int_{\mathbb{R}^{3d}} \gamma(n_1,n_2) \Big( \frac{\vert x_1 - z' \vert^2}{d(n_1)}
+ \frac{\vert x_2 -
y \vert^2}{d(n_2)} + \frac{\vert x_3 - y' \vert^2}{d(n_3)}
\Big)^{\frac{-3d + 2}{2}}  \\
    & & \qquad h_{n_1}(x_1)h_{n_2}(x_2)h_{n_3}(x_3) u^\e( z' - y' +
    z) \e^{-2} V \Big( \frac{z' - y}{\epsilon} \Big) J (z'-y')
\overline{J}(y') dz' dy dy'. \nonumber
\end{eqnarray}
We write
\[
I=c_0(3d)\int_{\mathbb{R}^{3d}}  u^\e(z' - y' +
    z) \e^{-2} V \Big( \frac{z' - y}{\epsilon} \Big) J (z'-y')
    \overline{J}(y') G(z',y,y')
    dz' dy dy',
    \]
    where $G(z',y,y')$ is given by
    \[
   \sum_{n_1,n_2,n_3} \gamma(n_1,n_2)
    \int
    \bigg( \frac{\vert x_1 - z' \vert^2}{d(n_1)}  + \frac{\vert x_2 -
y \vert^2}{d(n_2)} + \frac{\vert x_3 - y' \vert^2}{d(n_3)}
\bigg)^{\frac{-3d + 2}{2}} h_{n_1}(x_1)h_{n_2}(x_2)h_{n_3}(x_3) dx_1 dx_2 dx_3.
\]
Using the elementary inequality $abc\le (a^2+b^2+c^2)^{3/2}$ we
deduce that $G(z',y,y')$ is at most
\begin{eqnarray}
    & &
\sum_{n_1,n_2,n_3}\gamma(n_1,n_2) {\left(d(n_1)d(n_2)d(n_3)\right)}^{d/2-1/3}
\int
    {\vert x_1 - z' \vert}^{-d+2/3}\ {\vert x_2 -
y \vert}^{-d+2/3} {\vert x_3 - y' \vert}^{-d+2/3} \\
    & & \qquad \qquad \qquad h_{n_1}(x_1)h_{n_2}(x_2)h_{n_3}(x_3) dx_1
dx_2 dx_3.\nonumber
\end{eqnarray}
   From our asumptions on $h_n$ we deduce that $G\in L^{\i}_{loc}$. Hence,
\begin{displaymath}
    I  \leq   C \int_{\mathbb{R}^{3d}}  u^\e(z' - y' +
    z) \e^{-2} V \Big( \frac{z' - y}{\epsilon} \Big) J (z'-y')
    \overline{J}(y') dz' dy dy'.
\end{displaymath}
Note that, for fixed $(z',y') \in \mathbb{R}^{2d}$,
\begin{displaymath}
    \int_{\mathbb{R}^d}{V \Big( \frac{z' - y}{\epsilon}\Big) dy} =
\epsilon^d.
\end{displaymath}
Thus,
    \begin{eqnarray}
    I & \leq &  C \epsilon^{d-2} \int_{\mathbb{R}^{2d}}  u^\e(z' - y' +
    z)  J (z'-y') \overline{J}(y') dz' dy' \\
    & \leq & C \epsilon^{d-2} \int_{K} \int u^\e(z' - y' +
    z) J(z'-y')dz' dy'  \nonumber \\
     & \leq & C \epsilon^{d-2} \int_{K} dy' \int \frac{1}{\vert z' - y'+z
     \vert^{d-2}}J(z'-y')dz' \nonumber \\
     &\leq& C \epsilon^{d-2}  \int_{|a|\le \rho+|z|} \frac{1}{\vert a
     \vert^{d-2}}da\leq C \epsilon^{d-2}(\rho+|z|)^2, \nonumber
\end{eqnarray}
where $K \subseteq \mathbb{R}^d$ denotes a compact set containing the
support $\overline{J}$, and where we made use of the first
part of the Lemma in the third inequality. This is the bound stated in
    (\ref{leptu}). $\Box$
\subsection{Estimating the terms}\label{stthree}
\subsubsection{The case of $H_{14}$}
Note that
\begin{eqnarray}
    \mathbb{E}_N \Big\vert \int_0^T H_{14} (t) dt \Big\vert & \leq &
    C \epsilon^{2(d-2)} \mathbb{E}_N \int_0^T \sum_{i,j \in
    I_q}{ \Big\vert \ue \big( x_i - x_j + z \big) - \ue
    \big( x_i - x_j \big) \Big\vert } \\
     & & \qquad \quad 1\!\!1 \Big\{ m_i = M_1, m_j = M_2 \Big\}, \nonumber
\end{eqnarray}
where the constant $C$ depends on the $L^{\infty}$ bounds satisifed by
$J,\overline{J}$ and their time derivatives.
Hence
$$
  \mathbb{E}_N \Big\vert \int_0^T H_{14} (t) dt \Big\vert
  \leq K_1 + K_2,
$$
where $K_1$ is given by
\begin{displaymath}
    C \epsilon^{2(d-2)} \mathbb{E}_N \int_0^T{} dt\ \sum_{i,j \in I_q:
\vert x_i - x_j \vert >
\rho}{\Big\vert \ue(
x_i - x_j + z ) - \ue(
x_i - x_j  ) }\Big\vert  1 \! \! 1 \{ m_i =
M_1 \} \, 1 \! \!  1 \{ m_j =
M_2 \}
\end{displaymath}
and $K_2$ is given by
\begin{displaymath}
   \epsilon^{2(d-2)} \mathbb{E}_N \int_0^T{} dt\ \sum_{i,j \in I_q:
\vert x_i - x_j \vert \leq
\rho}{\Big\vert \ue(
x_i - x_j + z ) - \ue(
x_i - x_j  ) } \Big\vert \, 1 \! \! 1 \{ m_i =
M_1 \} 1 \! \!  1 \{ m_j = M_2 \}  .
\end{displaymath}
Firstly, we treat $K_1$. Note that
\begin{eqnarray}
K_1  & \leq &  \frac{C \vert z \vert \epsilon^{2(d-2)}}{\rho^{d-1}}
\mathbb{E}_N \int_0^T{}
dt\ \sum_{i,j \in I_q: \vert x_i - x_j \vert > \rho}{ 1 \! \! 1 \{ m_i = M_1
\}   1 \! \! 1 \{ m_j = M_2 \} }
\\
& \leq & \frac{C \vert z \vert}{ \rho^{d-1}} \epsilon^{2(d-2)}
\mathbb{E}_N \int_0^T{} dt\ \sum_{i,j \in I_q: \vert x_i - x_j \vert >
\rho}{m_i m_j} \leq    \frac{C \vert z \vert T }{ \rho^{d-1}}
    \Big( \epsilon^{d-2} \sum_{i \in I_{q(0)}}{m_i} \Big)^2 \leq   \frac{C
    \vert z \vert T }{ \rho^{d-1}}  , \nonumber
\end{eqnarray}
where the first inequality follows from the second part of Lemma
\ref{phiest}, and the final one from an hypothesis on the initial conditions of
the system.
We now treat the term $K_2$. By writing,
\begin{eqnarray}
    K_2 & \leq & C \epsilon^{2(d-2)} \mathbb{E}_N \int_0^T{} dt\ \sum_{i,j \in
I_q: \vert x_i - x_j \vert \leq
\rho}{
    \bigg[ \ue \big(
x_i - x_j + z \big) + \ue \big(
x_i - x_j  \big) \bigg] }  \\
    & &   \qquad \qquad \qquad \qquad  \qquad \qquad 1 \! \! 1 \{ m_i =
M_1 \} 1 \! \!  1 \{ m_j = M_2 \}  , \nonumber
\end{eqnarray}
we obtain an expression on the right-hand-side which may be bounded by
applying Lemma \ref{lemtwop}. The relevant estimate is provided by
the third part of Lemma \ref{phiest}. We find that $K_2 \leq C \big(
\rho + \vert z \vert \big)^2$. Hence,
\begin{displaymath}
\mathbb{E}_N \left|\int_0^T{}  H_{14}(t) \right| dt  \leq  K_1 + K_2
\leq    \frac{C
    \vert z \vert T }{ \rho^{d-1}}   +  C \big( \rho + \vert z \vert \big)^2  .
\end{displaymath}
Making the choice $\rho = \vert z \vert^{\frac{1}{d+1}}$, we find that
\begin{displaymath}
\mathbb{E}_N \left|\int_0^T{}  H_{14}(t) \right| dt  \leq C   \vert z
\vert^{\frac{2}{d+1}}  .
\end{displaymath}
\subsubsection{The cases of $H_2$ and $H_3$}
The estimate of
$\mathbb{E}_N \int_{0}^{T}{\vert H_3 (t) \vert dt}$ is derived in an
identical fashion to that of
$\mathbb{E}_N \int_{0}^{T}{\vert H_2 (t) \vert dt}$.
Picking $\rho \in \mathbb{R}$ that satisfies $\rho \geq \max \big\{ 2 \vert z
\vert + C_0 \epsilon , 2 C_0 \epsilon \big\}$,
we write
$$
\int_{0}^{T}{\mathbb{E}_N{\vert H_2 (t) \vert} dt} \leq R_1 + R_2,
$$
where
\begin{eqnarray}
R_1 & = & 2 \epsilon^{2(d-2)} \mathbb{E}_N{} \int_{0}^{T}{\sum_{i,j \in
I_q: \vert x_i
- x_j \vert > \rho}{d(m_i)  \left|\overline{J}(x_j,m_j,t)\right|}} \\
    & & \qquad \quad \Big\vert \ue_x ( x_i - x_j + z ) -
    \ue_x ( x_i - x_j ) \Big\vert
      \big\vert J_x \big( x_i,m_i,t \big) \big\vert dt
, \nonumber
\end{eqnarray}
and
\begin{eqnarray}
R_2 & = & 2 \epsilon^{2(d-2)} \mathbb{E}_N \int_{0}^{T}{\sum_{i,j \in
I_q: \vert x_i
- x_j \vert \leq \rho}{d(m_i)  \big|\overline{J}(x_j,m_j,t)\big|}} \\
& & \qquad \quad
    \Big\vert \ue_x( x_i - x_j + z ) -
    \ue_x ( x_i - x_j ) \Big\vert
     \big\vert J_x \big( x_i,m_i,t \big) \big\vert dt . \nonumber
\end{eqnarray}
Firstly, we examine the sum $R_1$.
Recalling that we consider test functions $J$ and $\overline{J}$
respectively supported on particles of mass $M_1$ and $M_2$,
\begin{displaymath}
R_1  \leq  C \epsilon^{2(d-2)}   \frac{\vert
z \vert}{\rho^d} d(M_1)
\Big\vert \{ (i,j) \in I_q^2: m_i = M_1, m_j = M_2 \} \Big\vert,
\end{displaymath}
where the lower bound on $\rho$ allowed us to apply the second part of
Lemma \ref{phiest}. Thus, $R_1 \leq C  \vert z \vert/\rho^d$.

Secondly, we bound the sum $R_2$. Note that
\begin{eqnarray}\label{fhg}
R_2 & \leq & 2 \epsilon^{2(d-2)}  \vert \vert J_x \vert \vert \ \vert \vert
\overline{J} \vert \vert d(M_1)
    \  \mathbb{E}_N \int_0^T{} \sum_{i,j \in I_q: \vert x_i - x_j
\vert \leq \rho} \\
    & & \qquad   \Big\vert
    \ue_x(x_i - x_j + z ) - \ue_x(
x_i - x_j  ) \Big\vert  1 \! \! 1 \{ m_i = M_1, m_j = M_2 \} \nonumber \\
& \leq &  C \epsilon^{2(d-2)} \ \mathbb{E}_N \int_0^T{} \sum_{i,j \in I_q}{m_i
m_j \big( d(m_i) + d(m_j) \big)} \nonumber \\
    & & \qquad \qquad \quad  \bigg[ \Big\vert
    \ue_x(x_i - x_j + z ) \Big\vert + \Big\vert \ue_x(
x_i - x_j  ) \Big\vert   \bigg] , \nonumber
\end{eqnarray}
where $\|\cdot \|$ denotes the $L^\i$ norm and
the constant $C$ depends on the test functions $J$ and $\overline{J}$.
The expression (\ref{fhg}) is written in a form to which Lemma \ref{lemtwop}
may be applied. Doing so yields
$$
R_2 \leq  C \epsilon^{2(d-2)}   \mathbb{E}_N  \sum_{i,j
\in I_q}{m_i m_j \Big( \hat{H}_0 (x_i - x_j) + \hat{H}_z (x_i - x_j) \Big) },
$$
where the functions $\{ \hat{H}_z : z \in
\mathbb{R}^d \}$ appear in the fourth part of Lemma
\ref{phiest}. From that result and our assumptions on the initial
data, we deduce that
$R_2 \leq C \big( \rho + \vert z \vert \big)$.
Hence,
\begin{displaymath}
\int_0^T{\mathbb{E}_N{\vert H_2(t) \vert} dt}  \leq  R_1 + R_2
    \leq C  \frac{\vert z \vert}{\rho^d}  + C \big( \rho + \vert z \vert
    \big),
\end{displaymath}
for $\rho\le 1$. By making the choice $\rho =  \vert z
\vert^{\frac{1}{d+1}} $, we find
that
\begin{displaymath}
\int_0^T{\mathbb{E}_N{\vert H_2(t) \vert} dt}  \leq  C  \vert z
\vert^{\frac{1}{d+1}}.
\end{displaymath}
\subsubsection{The case of $H_4$}
Note that
\begin{eqnarray}
\int_0^T{\mathbb{E}_N{\big\vert H_4(t)  \big\vert dt}} & \leq & \Big(
d(M_1) \vert \vert
\Delta J \vert \vert \ \vert \vert  \overline{J} \vert \vert + d(M_2)
\vert \vert
    J \vert \vert \ \vert \vert  \Delta\overline{J} \vert \vert \Big)
    \epsilon^{2(d-2)} \nonumber \\
    & & \quad
\mathbb{E}_N \int_0^T{} dt \sum_{i,j \in
I_q}{1 \! \! 1 \{
m_i=M_1,m_j=M_2 \} } \Big\vert \ue(x_i - x_j
+ z ) - \ue(x_i - x_j ) \Big\vert
.
\end{eqnarray}
This quantity is bounded as in the case of $H_{14}$.
\subsubsection{The case of $G_z(1) - G_0(1)$}
We now estimate the term
$$
\int_{0}^{T}{ \mathbb{E}_N \big\vert G_z(1) - G_0(1) \big\vert (t) dt}.
$$
To ease the notation, we do not display the dependece of $J$ and 
$\overline J$ on the variable $t$.
Note that
\begin{equation}\label{pkkl}
\int_{0}^{T}{ \mathbb{E}_N \big\vert G_z(1) - G_0(1) \big\vert (t) dt}
\leq \sum_{i=1}^{8}{D_i},
\end{equation}
where
\begin{eqnarray}
D_1 & = &\frac 12 \mathbb{E}_N \int_0^T{} dt \sum_{k,l \in
I_q}{\micp(m_k,m_l)V_{\epsilon}(x_k - x_l)} \big\vert
J(x_k,m_k)\big\vert   \\
& & \quad \epsilon^{2(d-2)} \sum_{i \in I_q}{
    \big\vert\overline{J}(x_i,m_i)\big\vert   \Big\vert
\ue(x_k - x_i + z) - \ue(x_k - x_i )
\Big\vert} , \nonumber
\end{eqnarray}
each of the other seven terms on the right-hand-side of (\ref{pkkl})
differing from $D_1$ only in an inessential way.
Given this. the estimates involved for each of the eight cases are in
essence identical, and we examine only the case of $D_1. $We write $D_1 =
D^1 + D^2$,
where we have decomposed the inner $i$-indexed sum according to the
respective index sets
$$ \{ i \in I_q, i \not= k,l, \vert x_k - x_i \vert > \rho \} \
\textrm{and} \  \{ i \in I_q, i \not= k,l, \vert x_k - x_i \vert \leq \rho \}.
$$
Here, $\rho$ is a positive parameter that satisfies the bound $\rho
\geq \max \big\{ 2 \vert z \vert + C_0 \epsilon , 2 C_0 \epsilon \big\}$.
By the second part of Lemma \ref{phiest}, we have that
$$
D^1  \leq  \frac{C \vert z \vert
\epsilon^{d-2}}{\rho^{d-1}} \mathbb{E}_N \int_0^T{} dt \sum_{k,l \in
I_{q}}{\micp(m_k,m_l) V_{\epsilon}(x_k - x_l)},  \nonumber
$$
where we have also used the fact that the test functions $J$ and
$\overline{J}$ are each supported on the set of particles of
respective masses $M_1$ and $M_2$, and the fact that the total number
of particles living at any given time is bounded
above by $Z \epsilon^{2-d}$. From the bound on the collision that is
provided by Lemma
\ref{lembc}, follows
$$ D^1 \leq  \frac{C \vert z \vert}{\rho^{d-1}} .$$

To bound the term $D^2$, note that
$D^2$ is bounded above by
\begin{eqnarray}
    &  & C \mathbb{E}_N \int_0^T{} dt \sum_{k,l \in
I_{q}}{\micp(m_k,m_l) V_{\epsilon}(x_k - x_l) } 1 \! \! 1 \{ m_k =
M_1 \} \epsilon^{2(d-2)} \sum_{i \in I_q} \\
& & \quad {} 1 \! \! 1 \{ \vert x_i - x_k
\vert \leq \rho\} 1 \! \! 1 \{ m_i = M_2 \} \Big\vert
\ue(x_k - x_i + z) - \ue(x_k - x_i ) \Big\vert \big|J(x_i,m_i)\big| 
\nonumber \\
    & \leq & C \mathbb{E}_N \int_0^T{} dt \sum_{k,l \in
I_{q}}{\micp(m_k,m_l)V_{\epsilon}(x_k - x_l)} 1 \! \! 1 \{ m_k = M_1
\} \epsilon^{2(d-2)} \sum_{i \in I_q}\nonumber \\
& & \quad {} 1 \! \! 1 \{ \vert x_i - x_k
\vert \leq \rho \} 1 \! \! 1 \{ m_i = M_2 \}  \Big[
|\ue(x_k - x_i + z)| + |\ue(x_k - x_i )| \Big]\big|J(x_i,m_i)\big| . \nonumber
\end{eqnarray}
Note that the last expectation is bounded by
Lemma \ref{lemthreep} because by our assumption on $\micp$, we can find
$\g$ such that
$\micp\le \g$ and $\g$ satisfies the assumption of Lemma 3.3. The upper
bound provided by this Lemma in this
particular application is computed in the fifth part of Lemma
\ref{phiest}. We find that $D^2 \leq C (\rho+|z|)^2$.

Combining these estimates yields
$$
D_1 \leq D^1 + D^2 \leq C \frac{\vert z \vert}{\rho^{d-1}} +  C  (\rho+|z|)^2 .
$$
Making the choice $\rho = \vert z \vert^{\frac{1}{d+1}} $ leads to the
inequality $D_1 \leq   \vert z \vert^{\frac{2}{d+1}}$.
Since each of the cases of $\big\{ D_i : i \in \{ 1,\ldots, 8 \} \big\}$ may
be treated by a nearly verbatim
proof, we deduce that
$$
\int_{0}^{T}{ \mathbb{E}_N \big\vert G_z(1) - G_0(1) \big\vert (t) dt}
\leq   C \vert z \vert^{\frac{2}{d+1}}.
$$
\subsubsection{The case of $G_z(2)$} Recall that
\begin{displaymath}
G_z(2)  =  - \epsilon^{2(d-2)} \sum_{k,l \in
I_{q}}{ \micp(m_k,m_l) V_{\epsilon}(x_k - x_l)} \ue(x_k -
x_l + z ) J(x_k,m_k)
    \overline{J}(x_l,m_l).
\end{displaymath}
If $k,l \in I_q$ satisfy $V_\epsilon \big( x_k - x_l \big) \not= 0$,
then $\vert x_k -
x_l \vert \leq C_0 \epsilon$, and so
$$
\vert x_k - x_l + z \vert \geq \vert z \vert - C_0 \epsilon
\geq \vert z \vert/2 ,
$$
provided that $|z|\ge 2C_0\e$. This implies that
$$
    \ue(x_k - x_l + z ) \leq
    \frac{C}{\vert x_k - x_l + z \vert^{d-2}} \leq \frac{C}{\vert z 
\vert^{d-2}},
$$
where in the first inequality, we used the first part of Lemma \ref{phiest}.
Applying this bound, we find that
$$
\int_0^T{} \mathbb{E}_N \vert G_z (2) \vert dt \leq  \frac{C
\epsilon^{2(d-2)}}{\vert z
\vert^{d-2}} \mathbb{E}_N
\int_0^T{} \sum_{k,l \in I_q}{ \micp(m_k,m_l) V_{\epsilon}(x_k - x_l)} dt
$$
whose right-hand-side is bounded above by $C\epsilon^{d-2}/{\vert z
\vert^{d-2}}$, according to Lemma \ref{lembc}. As a result,
$$
\int_0^T{\vert G_z(2) \vert dt} \leq C \Big( \frac{\epsilon}{\vert z
\vert} \Big)^{d-2}.
$$
\subsubsection{The case of $\mathbb{E}_N \vert X_z - X_0 \vert $}\label{lak}

We now turn to $\mathbb{E}_N \vert X_z - X_0 \vert $. Assume that 
$|z|\ge C_0\e$.
Using the first and second part of Lemma 3.5 we have that
\begin{eqnarray}
\mathbb{E}_N \vert X_z - X_0 (0)\vert \le & & C|z|
\iint_{L^2} h_{M_1}(x)h_{M_2}(y)|x-y|^{1-d}
dxdy \nonumber\\
& & +C\iint_{L^2} h_{M_1}(x)h_{M_2}(y)|x-y|^{2-d}
1\!\!1({|x-y|\le 3|z|})dxdy,
\\
& & +C\iint_{L^2} h_{M_1}(x)h_{M_2}(y)|x-y+z|^{2-d}
1\!\!1({|x-y+z|\le 4|z|})dxdy , \nonumber
\end{eqnarray}
where $L$ is a bounded set that contains the support of $J$ and $\overline J$.
Using our assumptions on the initial data $h_n$ we obtain the bound 
$C|z|$ for the first term and
$C|z|^{2/3}$ for the other two terms on the right-hand-side. This evidently
implies our bound
on $\mathbb{E}_N \vert X_z - X_0 \vert(0)$ in (\ref{jayeye}).
\begin{subsection}{The martingale term}\label{stfour}
This section is devoted to proving the estimate (\ref{martest}).
Note that
$$
M_z(T) = X_z \big( q(T),T \big) - X_z \big( q(0),0 \big) -
\int_{0}^{T}{\left(\frac \partial {\partial t}+\mathbb{L} \right)X_z 
\big( q(t),t) \big) dt}
$$
is a martingale which satisfies
$$
\mathbb{E}_N \Big[  M_z(T)^2 \Big] = \sum_{i=1}^{3}{}
\mathbb{E}_N\int_{0}^{T}{A_i \big( q(t)
,t
\big) dt},
$$
where $A_1(q,t)$ is set equal to
\begin{displaymath}
  2\epsilon^{4(d-2)}\sum_{i \in I_q, m_i  = M_1}{d(M_1)}
  \bigg[ \nabla_{x_i}
\sum_{j \in I_q, m_j = M_2}{\ue(x_i - x_j + z)
J(x_i,M_1,t) \overline{J}(x_j,M_2,t) } \bigg]^2,
\end{displaymath}
with $A_2$ being given by the same expression with instances of $M_1$
and $M_2$ being interchanged, including those implicitly appearing in
the function $\ue = \ue_{M_1,M_2}$. The term $A_3$ is given by
\begin{eqnarray}\label{aatwo}
&  & \! \! \! \!  \! \! \! \! \! \! \! \! \! \! \! \!  \! \! \! \! \!
\! \!  \frac 12\epsilon^{4(d-2)}\sum_{i,j \in I_q}{\micp(m_i,m_j)
V_{\epsilon} \big( x_i - x_j \big)} \nonumber \\
     & & \qquad \quad \bigg\{ \sum_{k \in I_q}{} \Big[
     \frac{m_i}{m_i + m_j} \ue (x_i - x_k+z )  J(x_i,m_i +
     m_j) \overline{J}(x_k,m_k) \nonumber \\
    & & \qquad \qquad \qquad + \
     \frac{m_i}{m_i + m_j} \ue (x_k - x_i + z)  J(x_k,m_k)
\overline{J}(x_i,m_i + m_j) \nonumber \\
    & & \qquad \qquad \qquad + \
     \frac{m_j}{m_i + m_j} \ue (x_j - x_k + z) J(x_j,m_i +
     m_j) \overline{J}(x_k,m_k) \nonumber \\
    & & \qquad \qquad \qquad + \
     \frac{m_j}{m_i + m_j} \ue (x_k - x_j + z)  J(x_k,m_k)
\overline{J}(x_j,m_i + m_j) \\
    & & \qquad \qquad \qquad - \  \ue (x_i - x_k + z)
     J(x_i,m_i) \overline{J}(x_k,m_k) \, - \,  \ue (x_k - x_i + z)
     J(x_k,m_k) \overline{J}(x_i,m_i) \nonumber \\
      & & \qquad \qquad \qquad - \   \ue (x_j
     - x_k + z)   J(x_j,m_j) \overline{J}(x_k,m_k) \, - \,  \ue (x_k
     - x_j + z)  J(x_k,m_k) \overline{J}(x_j,m_j) \Big] \nonumber \\
     & & \qquad \qquad \qquad \qquad \qquad \qquad \qquad - \  \ue (x_i
     - x_j + z) J(x_i,m_i) \overline{J}(x_j,m_j)
      \bigg\}^2
. \nonumber
\end{eqnarray}
Recall that we by our convention, we do not display the dependence of 
$J$ and $\overline J$on the
$t$-variable. To bound these terms, we require two variants of Lemma 
\ref{lemthreep}:
\begin{lemma}\label{tpln}
There exists a collection of constants $C: \mathbb{N}^2 \to
(0,\infty)$ such that,
for any continuous functions $t,v,a_1,a_2,a_3: \mathbb{R}^d \to 
[0,\infty)$ and any $z\in \R^d$,
\begin{eqnarray}
& & \mathbb{E}_N\int_0^T{}dt \sum_{i,j,k \in I_q(t)}{}\g(m_i,m_j)t 
\Big( \frac{x_i -
x_j+z}{\epsilon} \Big) v \Big( \frac{x_i - x_k+z}{\epsilon} \Big)
\\
    & & \qquad \qquad \qquad a_1(x_i)
a_2(x_j) a_3(x_k) 1 \! \! 1 \big\{ m_i = n_1 ,
m_k = n_3 \big\}   \nonumber \\
& & \qquad \qquad \quad  \leq C_{n_1,n_3} \epsilon^{3(2-d)}\sum_{n_2}
\mathbb{E}_N
\sum_{i,j,k \in I_{q(0)} }A_{n_1,n_2,n_3}^{\epsilon}(x_i,x_j,x_k)  , \nonumber
\end{eqnarray}
where $A_{n_1,n_2,n_3}^{\epsilon}: \mathbb{R}^{3d} \to [0,\infty)$ is
given by
\begin{eqnarray}
     A_{n_1,n_2,n_3}^{\epsilon}(x_1,x_2,x_3) & = &
c_0(3d)\g(n_1,n_2)\int_{\mathbb{R}^{3d}} \bigg(
\frac{\vert x_1 - z' \vert^2}{d(n_1)} + \frac{\vert x_2 - y
\vert^2}{d(n_2)} + \frac{\vert x_3 - y' \vert^2}{d(n_3)}
\bigg)^{\frac{-3d + 2}{2}} \nonumber \\
    & & \quad  t \Big( \frac{z' - y+z}{\epsilon} \Big) v \Big(
\frac{z' - y'+z}{\epsilon} \Big) a_1(z') a_2(y) a_3(y') dz' dy dy,
\end{eqnarray}
with $\g$ as in Lemma 3.3.
\end{lemma}
\begin{lemma}\label{fpln}
There exists a collection of constants $C: \mathbb{N}^3 \to
[0,\infty)$ such that,
for any $z\in \R^d$,
  any continuous functions $v,w: \mathbb{R}^d \to [0,\infty)$ and
another $(a_1,a_2,a_3): \mathbb{R}^{3d} \to [0,\infty)$,
\begin{eqnarray}
& & \mathbb{E}_N\int_0^T{}dt \sum_{k,l,i,j \in I_q}{}\g(n_i,n_j)V_{\e}( {x_i -
x_j}) v \Big( \frac{x_i - x_k+z}{\epsilon} \Big) w \Big(
\frac{x_i - x_l+z}{\epsilon} \Big) a_1( x_i)a_2(x_k)a_3(x_l) \\
    & & \qquad \qquad \qquad  1 \! \! 1 \big\{ m_i = n_1 ,  m_k = n_3,
m_l = n_4 \big\}   \nonumber \\
& & \qquad \quad \leq C_{n_1,n_3,n_4} \epsilon^{4(2-d)}
\sum_{n_2}\mathbb{E}_N \sum_{i,j,k,l
\in I_{q(0)}}{B_{m_i,m_j,m_k,m_l}^{\epsilon}(x_i,x_j,x_k,x_l)} \nonumber
\\
& & \qquad \qquad \qquad \qquad 1 \! \! 1 \big\{ m_i \leq n_1 , m_j
\leq n_2, m_k \leq n_3,
m_l \leq n_4 \big\} , \nonumber
\end{eqnarray}
where  $B_{n_1,n_2,n_3,n_4}^{\epsilon}: \mathbb{R}^{4d} \to [0,\infty)$ is
given by
\begin{eqnarray}
    & & B_{n_1,n_2,n_3,n_4}^{\epsilon} \big( x_1,x_2,x_3,x_4 \big) \\
    & = & \int_{\mathbb{R}^{4d}} \bigg(
\frac{\vert x_1 - \hat{z} \vert^2}{d(n_1)} + \frac{\vert x_2 - z'
\vert^2}{d(n_2)} + \frac{\vert x_3 - y \vert^2}{d(n_3)} +
\frac{\vert x_4 - y' \vert^2}{d(n_4)}
\bigg)^{-2d +1} \nonumber \\
    & & \qquad \g(n_1,n_2)V_\e( {\hat{z} - z'}) v \Big(
\frac{\hat{z} - y+z}{\epsilon} \Big) w \Big( \frac{z' - y'+z}{\epsilon} \Big)
    a_1( \hat{z})a_2(y)a_3( y') d \hat{z}
dz' dy dy', \nonumber
\end{eqnarray}
where the function $\g:\mathbb{N}^2\to (0,\i)$ satisfies
\begin{displaymath}
n_2 \g \big( n_1, n_2 + n_3 \big) \max{ \Big\{ 1, \Big[ \frac{d (
n_2 + n_3 )}{d(n_2)}
\Big]^{2d-1}\Big\}} \leq \big( n_2 + n_3 \big) \g(n_1,n_2),
\end{displaymath}
\end{lemma}

The proof of Lemma 3.5 is identical to that of Lemma 3.3. The proof
of Lemma 3.6 is very similar
to the proof of Lemma 3.3 and is omitted.

We now bound the three terms. Of the first two, we treat only $A_1$,
the other being bounded by an identical argument.
By multiplying out the brackets appearing in the definition of $A_1$
and using the facts that
\begin{displaymath}
     \ue  (x_i - x_j + z)  =  \epsilon^{2 - d}
    u \Big( \frac{x_i - x_j + z}{ \epsilon} \Big)
,
\end{displaymath}
and
\begin{displaymath}
    \nabla_{x_i} \ue (x_i - x_j + z)  =
\epsilon^{1 - d}  \nabla u  \Big( \frac{x_i - x_j + z}{\epsilon} \Big),
\end{displaymath}
we obtain that $A_1$ is bounded above by
\begin{eqnarray}\label{yhj}
    & & C \epsilon^{2(d-3)} \sum_{i,j,k \in I_q} \big\vert  \nabla u
    \big\vert \Big( \frac{x_i - x_j + z}{\epsilon} \Big)   \big\vert
    \nabla u \big\vert \Big( \frac{x_i -
x_k + z}{\epsilon} \Big)  J^2(x_i,m_i) \\
& &  \qquad \qquad
    \big|\overline{J}(x_j,m_j)\big| \big|\overline{J}(x_k,m_k)\big|  1 
\! \! 1 \big\{
m_i = M_1 , m_j = m_k = M_2
    \big\} \nonumber \\
    & + & C \epsilon^{2(d-2)}   \sum_{i,j,k \in I_q}  u \Big( \frac{x_i
- x_j + z}{\epsilon} \Big)   u \Big( \frac{x_i -
x_k + z}{\epsilon} \Big) |\nabla J(x_i,m_i)|^2 \nonumber \\
    & &  \qquad \qquad
    \big|\overline{J}(x_j,m_j)\big| \big|\overline{J}(x_k,m_k)\big|  1 
\! \! 1 \big\{
m_i = M_1 , m_j = m_k = M_2
    \big\}  . \nonumber
\end{eqnarray}
Let us assume that $z=0$ because this will not affect out arguments.
We are required to bound the quantity appearing in
the statement of Lemma \ref{tpln}, for each of the following cases:
\begin{equation}\label{twocases}
    \big( t,v,a_1,a_2,a_3 \big)  \in  \Big\{ \big( \vert\nabla u
\vert , \vert \nabla u \vert,
    J^2, \big|\overline{J}\big|,  \big|\overline{J}\big|\big) \ , \   \big(
  u,  u,  \vert \nabla J \vert^2 ,  \big|\overline{J}\big|,  \big|\overline{J}
\big|
    \big) \Big\}.
\end{equation}
Recall that each  of the test functions $J$, $\overline{J}$, and
their gradients,
is assumed to be uniformly bounded with compact support. To each of
the two cases, Lemma \ref{tpln} applies. For either of them, the
right-hand-side of the inequality in Lemma \ref{tpln} may be written as a
finite sum of the expectations appearing there, with the sum being
taken over triples of given masses $n_1,n_2$ and $n_3$. Such an expectation is
bounded above by
\begin{eqnarray}\label{yun}
    & & C \epsilon^{3(2-d)}\sum_{n_2}
    \int \int_{K^3} \bigg(
\frac{\vert x_1 - z' \vert^2}{d(n_1)} + \frac{\vert x_2 - y
\vert^2}{d(n_2)} + \frac{\vert x_3 - y' \vert^2}{d(n_3)}
\bigg)^{\frac{-3d + 2}{2}} \\
    & & \qquad \qquad t \Big( \frac{z' - y}{\epsilon} \Big) v \Big(
\frac{z' - y'}{\epsilon} \Big) h_{n_1}(x_1)h_{n_2}(x_2)h_{n_3}(x_3) dz'
dy dy'  d x_1 d x_2 d x_3 , \nonumber
\end{eqnarray}
where $K=\{x:|x|\le \ell\} \subseteq \mathbb{R}^d$ is chosen to
contain the support of $J$ and $\overline{J}$. We bound
(\ref{yun}) for the first of the two cases in (\ref{twocases}), the
other three being similar. In this case, $t = v = \vert \nabla u
\vert$, the expression (\ref{yun}) equals
\[
    C \epsilon^{3(2-d)} \bigg[ \int_{K^3}{}   \big\vert \nabla u
\big\vert  \Big( \frac{z' - y}{\epsilon} \Big)  \big\vert \nabla u
\big\vert \Big(
\frac{z' - y'}{\epsilon} \Big) G(z',y,y')   dz' dy dy' \bigg] ,
\]
where $G(z',y,y')$ is defined to be
\[
\sum_{n_2}\int  \bigg(
\frac{\vert x_1 - z' \vert^2}{d(n_1)} + \frac{\vert x_2 - y
\vert^2}{d(n_2)} + \frac{\vert x_3 - y' \vert^2}{d(n_3)}
\bigg)^{\frac{-3d + 2}{2}}h_{n_1}(x_1)h_{n_2}(x_2)h_{n_3}(x_3)d x_1 d
x_2 d x_3.
\]
As in the proof of the fifth part of Lemma 3.5 we can show that $G\in 
L^\i_{loc}$.
To bound the expression
\begin{equation}\label{polk}
    \int_{K^3}{}   \big\vert \nabla u
\big\vert \Big( \frac{z' - y}{\epsilon} \Big)  \big\vert \nabla u
\big\vert \Big(
\frac{z' - y'}{\epsilon} \Big) dz' dy dy' ,
\end{equation}
note that, for $z' \in K$,
$$\label{uone}
    \int_{K^2}{dy dy'  \big\vert \nabla u
\big\vert  \Big( \frac{z' - y}{\epsilon} \Big)  \big\vert \nabla u
\big\vert \Big(
\frac{z' - y'}{\epsilon} \Big)  }   \leq
    \bigg[ \int_{|a|\le 2\ell}  { \big\vert \nabla u
\big\vert \Big( \frac{a}{\epsilon} \Big)} da \bigg]^2 .
$$
Using  $\vert \nabla u \vert (x) \leq C\vert x \vert^{1-d}
$ yields a bound $
    C \ell^2\epsilon^{2(d-1)} $
for (\ref{polk}).
Hence, the expression appearing in (\ref{yun}) is bounded above by
$ C \epsilon^{3(2-d)} \epsilon^{2(d-1)} $.

Similar reasoning applies in the second case listed in
(\ref{twocases}). The bound on the term in (\ref{yun}) in this case is
a constant multiple of $\epsilon^{3(2-d)} \epsilon^{2
(d-2)}$. Applying Lemma \ref{tpln} to the expression in (\ref{yhj}), we
find that
\begin{displaymath}
\mathbb{E}_N\int_0^T{A_1 \big( q(t) \big) dt} \leq
    C  \epsilon^{2(d-3)}
\epsilon^{3(2 - d)}  \epsilon^{2(d-1)} +
C  \epsilon^{2(d - 2)}
\epsilon^{3(2 - d)}  \epsilon^{2(d - 2)}.
\end{displaymath}
These bounds correspond to the two cases appearing in
(\ref{twocases}). Thus,
$$
\mathbb{E}_N\int_0^T{A_1 \big( q(t) \big) dt} \leq C \epsilon^{d-2} .
$$
We must treat the third term, $A_3$.
An application of the inequality
$$
\Big( a_1 + \ldots + a_n \Big)^2 \leq n \Big( a_1^2 + \ldots + a_n^2 \Big)
$$
to the bound on $A_3$ provided in (\ref{aatwo}) implies that
\begin{equation}\label{aathr}
    A_3(q) \leq \frac 92\epsilon^{4(d-2)} \sum_{i,j \in I_q}{\micp (m_i,m_j)
    \epsilon^{-2} V \Big( \frac{x_i - x_j}{\epsilon} \Big) \bigg[
    \sum_{n=1}^{8}{\Big( \sum_{k \in I_q}Y_n \Big)^2 } + Y_9^2  \bigg]},
\end{equation}
where $Y_1$ is given by
\begin{displaymath}
     \frac{m_i}{m_i + m_j} \ue (x_i - x_k + z ) J(x_i,m_i +
     m_j) \overline{J}(x_k,m_k),
\end{displaymath}
and where $\big\{ Y_i : i \in \{ 2, \ldots, 8 \} \big\}$ denote the
other seven expressions in (\ref{aatwo}) that appear in a sum over $k
\in I_q$, while $Y_9$ denotes the last term in (\ref{aatwo})
that does not appear in this sum.
There are nine cases to consider. The first eight are practically
identical, and we treat only the fifth. Note that
\begin{eqnarray}\label{zvc}
    & &  \epsilon^{4(d-2)} \sum_{i,j \in I_q}{\micp (m_i,m_j)  V_\e (
{x_i - x_j}) \Big(
    \sum_{k \in I_q}Y_5 \Big)^2 } \\
    & =& C \epsilon^{2d-4} \sum_{i,j \in I_q}{\micp(m_i,m_j)
    V_\e( x_i - x_j )} \nonumber \\
    & &  \qquad \Big[ \sum_{k,l \in I_q}{u \Big( \frac{x_i - x_k +
    z}{\epsilon} \Big) u \Big( \frac{x_i - x_l + z}{\epsilon} \Big)
    J^2(x_i,m_i) \overline{J}(x_k,m_k) \overline{J}(x_l,m_l) } \Big]. \nonumber
\end{eqnarray}
In the sum with indices involving $k,l \in I_q$, we permit the
possibility that these two may be equal, though they must be distinct
from each of $i$ and $j$ (which of course must themselves be distinct
by the overall convention).

Note that the expression (\ref{zvc}) appears in the statement of
Lemma \ref{fpln},
    provided that the choice
$$
\Big( v,w ,a_1,a_2,a_3\Big) = \Big(   u, u, J^2 ,\big| \overline{J}\big|,
    \big|\overline{J}\big|\Big)
$$
is made. Again we  set $z=0$ because this does not
affect the estimates.
Given that the support of the functions $a_1,a_2,a_3: \mathbb{R}^{4d} \to
    [0,\infty)$ are bounded,
we must bound
\begin{eqnarray}
    & & \sum_{n_2}\int  \int_{L^4} \Big(
\frac{\vert x_1 - \hat{z} \vert^2}{d(n_1)} + \frac{\vert x_2 - z'
\vert^2}{d(n_2)} + \frac{\vert x_3 - y \vert^2}{d(n_3)} +
\frac{\vert x_4 - y' \vert^2}{d(n_4)}
\Big)^{-2d+1} V \Big( \frac{\hat{z} - z'}{\epsilon} \Big) \\
    & & \qquad   u \Big(
\frac{\hat{z} - y}{\epsilon}  \Big)  u \Big(
\frac{z' - y'}{\epsilon}  \Big)h_{n_1}(x_1)h_{n_2}(x_2)h_{n_3}(x_3)h_{n_4}(x_4)
    d \hat{z} dz' dy dy' d x_1 d x_2 d x_3  dx_4 , \nonumber
\end{eqnarray}
for  a compact set $L$. This expression is bounded above by
\begin{eqnarray}
     & & \int_{L^4}  V \Big( \frac{\hat{z} - z'}{\epsilon} \Big) u \Big(
\frac{\hat{z} - y}{\epsilon}  \Big)  u \Big(
\frac{z' - y'}{\epsilon}  \Big)  d \hat{z} dz' dy dy' \\
    & & \qquad  \sum_{n_2}{ \int_{K^4}  \bigg(
\frac{\vert x_1 - \hat{z} \vert^2}{d(n_1)} + \frac{\vert x_2 - z'
\vert^2}{d(n_2)} + \frac{\vert x_3 - y \vert^2}{d(n_3)} +
\frac{\vert x_4 - y' \vert^2}{d(n_4)}
    \bigg)^{-2d+1}  }  \nonumber \\
    & & \qquad \qquad
    h_{n_1}(x_1)h_{n_2}(x_2)h_{n_3}(x_3)h_{n_4}(x_4)dx_1 dx_2 dx_3 
dx_4  , \nonumber
\end{eqnarray}
which is less than
$$
C  \int_{L^4}  V \Big( \frac{\hat{z} - z'}{\epsilon} \Big) u \Big(
\frac{\hat{z} - y}{\epsilon}  \Big)  u \Big(
\frac{z' - y'}{\epsilon}  \Big)  d \hat{z} dz' dy dy'.
$$
The proof of this follows the proof of the fifth part of Lemma 3.5;
we use the elementary
inequality $abcd\le (a^2+b^2+c^2+d^2)^2$ and the fact that the function
$$
\hat k(x)=\sum_nd(n)^{d/2-1/4}\int h_n(y)|x-y|^{1/2-d}dy
$$
is locally bounded. Since
$\int_{L}{u \Big( \frac{\hat{z} - y}{\epsilon} \Big) dy}$
is bounded above by $C \epsilon^{d-2}$, we have that
$$
     \int_{L^4}  V \Big( \frac{\hat{z} - z'}{\epsilon} \Big) u \Big(
\frac{\hat{z} - y}{\epsilon}  \Big)  u \Big(
\frac{z' - y'}{\epsilon}  \Big)  d\hat{z} dz' dy dy'
     \leq  C \epsilon^{2(d-2)} \int_{L^2}  V \Big( \frac{z -
z'}{\epsilon} \Big)   d \hat{z} dz'.
    $$
This is at most $ C\e^{3d-4}$.  Applying Lemma \ref{fpln},
    we find that the contribution to
$$\mathbb{E}_N{\int_{0}^{T}{A_3 \big(q(t)\big) dt}}$$
arising from the fifth term in (\ref{aathr}) is at most
$C \epsilon^{2d- 4} \epsilon^{4(2-d)} \epsilon^{3d - 4} = \epsilon^{d}$.

We now treat the ninth term, as they are classified in
(\ref{aathr}). It takes the form
\begin{displaymath}
     \epsilon^{4d-8} \sum_{i,j \in I_q}{\micp(m_i,m_j)
     V_\e\Big( {x_i - x_j} \Big)} \ue \big(x_i - x_j + z \big)^2
      J \big( x_i,m_i \big)^2 \overline{J} \big( x_j,m_j \big)^2.
\end{displaymath}
This is bounded above by
$$
C\epsilon^{2d-4} \sum_{i,j \in I_q}{\micp(m_i,m_j)
     V_\e\Big( {x_i - x_j} \Big)},
     $$
  because $u^\e\le C\e^{2-d}$ by (3.20).
The expected value of the integral on the interval of time $[0,T]$
of this last expression is bounded above by
$$
C\e^{2d-4} \mathbb{E}_N \int_0^T{dt}
\sum_{i,j \in I_q}{\micp(m_i,m_j)
    V_{\e}(x_i - x_j)}\le C\e^{d-2} ,
$$
where we used Lemma \ref{lembc} for the last inequality.
This completes the proof of (\ref{martest}).
\end{subsection}
\subsection{Using the estimates}\label{stfive}
We now use the estimates derived in the preceding sections to prove
the Stosszahlansatz. The precise form of the \stoss that we derive is
now stated.
\begin{prop}\label{szprop}
Recall from (\ref{sfd}) that
\begin{equation}\label{qudef}
Q(0) = \epsilon^{d-2} \sum_{(i,j) \in
I_{q}}{\micp(m_i,m_j) V_{\epsilon}(x_i - x_j) J(x_i,m_i)
    \overline{J}(x_j,m_j)},
\end{equation}
where, as before,  we simply write $J(x_i,m_i)$ and $\overline J(x_i,m_i)$
for $J(x_i,m_i,t)$ and $\overline J(x_i,m_i,t)$. Let $\eta:
\mathbb{R}^d \to [0,\infty)$ denote a smooth function of compact
support for which $\int_{\mathbb{R}^d}{\eta(x)dx} = 1$.  
We have that
\begin{eqnarray}\label{szrec}
    & & \int_0^T  Q (0)(t) dt \\
    & = &
\int_0^T{} dt \int_{\mathbb{R}^d}
d\omega  \, \macp(M_1,M_2)  J(\omega,M_1)
\overline{J}(\omega,M_2)  \bigg[ \epsilon^{d-2} \sum_{i \in I_q: m_i = M_1}{
\delta^{-d} \eta \Big( \frac{x_i - \omega}{\delta} \Big) }
\bigg] \nonumber \\
& & \qquad \quad \bigg[ \epsilon^{d-2} \sum_{j \in I_q; m_j = M_2}{
\delta^{-d} \eta \Big(\frac{x_j - \omega}{\delta} \Big)
} \bigg]  \, + \rm{Err} \big( \epsilon, \delta
\big), \nonumber
\end{eqnarray}
where the constants $\macp: \mathbb{N}^2 \to [0,\infty)$ were
defined in (\ref{recp}), and
where the function $\rm{Err}$ satisfies
$$
  \lim_{\delta \downarrow 0}\limsup_{\epsilon\downarrow 0} \mathbb{E}_N
  \big\vert Err \big( \epsilon, \delta \big) \big\vert = 0.
$$
\end{prop} 
\noindent{\bf Proof.} Note firstly that
the inequalities (\ref{abc}), (\ref{jayeye}) and (\ref{martest}) imply
that, for large $T$,
$$
\lim_{\vert z \vert \to 0}\limsup_{\epsilon \downarrow
0}\mathbb{E}_N   \bigg\vert\int_0^T H_{11} \big( t \big) dt  +
\int_0^T H_{13} \big( t \big) dt\bigg\vert  = 0.
$$
That is,
\begin{eqnarray}
  & & \! \!  \! \! \! \! \! \! \! \! \! \! \! \!  \! \!  \! \! \! \! \!
  \! \! \! \!  \! \!  \! \! \! \! \! \! \! \lim_{\vert z \vert \to 
0}{\limsup_{\epsilon \downarrow
0}} \    \epsilon^{2(d-2)} \mathbb{E}_N  \bigg\vert
     \int_{0}^{T}{dt \sum_{i,j \in
    I_{q(t)}}{\micp(m_i,m_j) J(x_i,m_i) \overline{J}(x_j,m_j)} } \nonumber \\
& & \qquad \qquad \qquad  \qquad \qquad \qquad  \Big[ V^{\epsilon} 
(x_i - x_j + z) -
V^{\epsilon} ( x_i - x_j ) \Big] \\
& & \qquad \qquad \quad  + \  \int_{0}^{T}{dt}\sum_{i,j \in
I_{q(t)}} \micp(m_i,m_j)
     J(x_i,m_i) \overline{J}(x_j,m_j) \nonumber \\
& & \qquad  \qquad \qquad \qquad \qquad \qquad \quad
V_{\epsilon}(x_i - x_j + z) \ue(x_i - x_j + z)
     \ \bigg\vert  \   = \ 0. \nonumber
\end{eqnarray} 
This implies that
\begin{eqnarray}\label{otwdfw}
& &
    \epsilon^{2(d-2)} \int_{0}^{T} dt {\sum_{i,j \in
I_{q(t)}}{\micp(m_i,m_j)
    V^{\epsilon} ( x_i - x_j ) J(x_i,m_i)}
    \overline{J}(x_j,m_j) } \\
& = & \epsilon^{2(d-2)} \int_{0}^{T} dt {\sum_{i,j \in
I_{q(t)}}{\micp(m_i,m_j)
    V^{\epsilon}( x_i - x_j + z ) J(x_i,m_i)
    \overline{J}(x_j,m_j)} }\nonumber \\
& & \qquad \qquad \qquad  \ \Big[ 1 +
     \epsilon^{d-2} \ue(x_i - x_j
+ z)   \Big]  \ + \ err(\epsilon,z) \nonumber
\\
& = & \epsilon^{2(d-2)} \int_{0}^{T} dt {\sum_{i,j \in
I_{q(t)}}{\micp(m_i,m_j)
\epsilon^{-d}  U_{m_i,m_j} \Big( \frac{x_i -
x_j + z}{\epsilon} \Big) }} \,  J(x_i,m_i)
    \overline{J}(x_j,m_j)   \  + \ err(\epsilon,z) , \nonumber
\end{eqnarray}
where
\begin{equation}\label{umm}
U_{m_1,m_2}(x) = V(x) \bigg[ 1 +  u_{m_1,m_2} \big( x
\big) \bigg],
\end{equation}
and where $err$ satisfies
\begin{equation}\label{lati}
    \lim_{\vert z \vert \to 0}{\limsup_{\epsilon \downarrow
    0}}{ \, \mathbb{E}_N \big\vert err
    \big( \epsilon,z \big) \big\vert} = 0.
\end{equation}
Writing
\begin{equation}\label{quzdef}
    \overline{Q} (z)  =  \epsilon^{d-2}
    \sum_{i,j \in I_q}{\micp(m_i,m_j) \epsilon^{-2} U_{m_i,m_j}\Big(
    \frac{x_i - x_j + z
}{\epsilon} \Big) J(x_i,m_i)
    \overline{J}(x_j,m_j) },
\end{equation}
it follows from (\ref{otwdfw}) that
\begin{equation}\label{weakstoss}
\int_0^T  Q (0)(t) dt  =
\int_0^T \overline{Q} (z)(t) dt  + err\big(\epsilon,z \big).
\end{equation}
Note that
\begin{eqnarray}\label{stra}
  & & \bigg\vert \overline{Q}\big(z_2 - z_1 \big) -   \epsilon^{d-2}
    \sum_{i,j \in I_q}\micp(m_i,m_j) \epsilon^{-2} U_{m_i,m_j}\Big(
    \frac{x_i - x_j + z_2 - z_1 }{\epsilon} \Big)  J(x_i - z_1,m_i)
    \overline{J}(x_j - z_2,m_j)  \bigg\vert \nonumber \\
  & \leq &  \epsilon^{d-2}
    \sum_{i,j \in I_q}\micp(m_i,m_j) \epsilon^{-2} U_{m_i,m_j}\Big(
    \frac{x_i - x_j + z_2 - z_1 }{\epsilon} \Big) \nonumber \\
  & & \qquad \quad \Big\vert \big( J(x_i -
    z_1,m_i) - J(x_i,m_i) \big)
    \overline{J}(x_j - z_2,m_j) \, + \, J(x_i,m_i)  \big( J(x_j -
    z_2,m_j) - J(x_j,m_j) \big)   \Big\vert \\
  & \leq &  \epsilon^{d-2}
    \sum_{i,j \in I_q}\micp(m_i,m_j) \epsilon^{-2} U_{m_i,m_j}\Big(
    \frac{x_i - x_j + z_2 - z_1 }{\epsilon} \Big)  \big( \vert z_1
    \vert + \vert z_2 \vert \big)
  K(x_i,m_i)
    \overline{K}(x_j,m_j), \nonumber
\end{eqnarray}
for small $z_1,z_2$,
where $K,\overline{K}:\mathbb{N} \times \mathbb{R}^d \times [0,\infty)
\to [0,\infty)$ are two smooth compactly supported test functions
satisfying
$$
  K , \overline{K}\le 2 \big( \vert\vert \nabla J\vert\vert + \vert\vert
  J\vert\vert\big) 1 \! \! 1 \big\{ x \in \mathbb{R}^d : d (x,supp J\cup
  supp \overline{J} )
  \leq 1 \big\}.
$$
Note that the final term in (\ref{stra}) is equal to $(|z_1| + |z_2|)
\overline{Q}_K (z_2 - z_1)$, where $Q_K$ and $\overline{Q}_K$
are defined by (\ref{qudef}) and (\ref{quzdef}) with $K$ and
$\overline{K}$ replacing $J$ and $\overline{J}$.
By (\ref{weakstoss}) applied with $Q_K(0)$ and $\overline{Q}_K(z)$ in place
of $Q(0)$ and $\overline{Q}(z)$, and Lemma \ref{lembc}, we find that
\begin{equation}\label{strb}
\epsilon^{d-2}
   \mathbb{E}_N \Big\vert \int_0^T dt  \sum_{i,j \in 
I_q}\micp(m_i,m_j) \epsilon^{-2} U_{m_i,m_j}\Big(
    \frac{x_i - x_j + z_2 - z_1 }{\epsilon} \Big)  K(x_i,m_i)
    \overline{K}(x_j,m_j) \Big\vert \leq Z \vert\vert K \vert \vert \,
    \vert \vert \overline{K} \vert\vert.
\end{equation}
By (\ref{stra}) and (\ref{strb}),
\begin{eqnarray}\label{strc}
  \overline{Q}(z_2 - z_1) & = & \epsilon^{d-2}
    \sum_{i,j \in I_q}{\micp(m_i,m_j) \epsilon^{-2} U_{m_i,m_j}\Big(
    \frac{x_i - x_j + z_2 - z_1 }{\epsilon} \Big)} \\
   & & \qquad \qquad \qquad J(x_i - z_1,m_i)
    \overline{J}(x_j - z_2,m_j)  \, + \, err(\epsilon,z_1,z_2), \nonumber
\end{eqnarray}
where we have that
$
\mathbb{E}_N \big\vert err(\epsilon,z_1,z_2)\big\vert \leq
  Z \vert\vert K \vert \vert \,
    \vert \vert \overline{K} \vert\vert \big( \vert z_1 \vert + \vert
    z_2 \vert \big).$
By (\ref{weakstoss}) and (\ref{strc}),
\begin{eqnarray}
    & & \int_0^T  Q (0)(t) dt  \label{bhap} \\
    & = & \epsilon^{d-4} \int_{0}^{T}dt
\sum_{i,j \in I_{q(t)}}  J(x_i - z_1,m_i)
\overline{J}(x_j - z_2,m_j) \alpha(m_i,m_j) \nonumber \\
    & & \quad \int_{\mathbb{R}^d}{\int_{\mathbb{R}^d}{U_{m_i,m_j}
    \bigg( \frac{
(x_i - z_1) - (x_j - z_2)}{\epsilon} \bigg) {\delta}^{-d}
\eta \Big( \frac{z_1}{\delta} \Big) {\delta}^{-d} \eta \Big(
\frac{z_2}{\delta} \Big) dz_1
dz_2 }} \, + \, err\big(\epsilon,\delta\big) \nonumber \\
    &  =  &  \epsilon^{d-4} \int_{0}^{T} dt
\int_{\mathbb{R}^{2d}} d \omega_1 d \omega_2 \sum_{i,j
\in I_{q(t)}}{}U_{m_i,m_j} \Big( \frac{\omega_1 -
\omega_2}{\epsilon} \Big)  \micp(m_i,m_j)  J(\omega_1,m_i)
\overline{J}(\omega_2,m_j) \nonumber \\
    & & \qquad \quad
\delta^{-d} \eta \Big( \frac{x_i - \omega_1}{\delta} \Big)
\delta^{-d} \eta \Big( \frac{x_j - \omega_2}{\delta} \Big)
     \, + \, err \big( \epsilon, \delta \big) \nonumber \\
& = & \epsilon^{-d}
\int_0^T{} dt \int_{\mathbb{R}^{2d}} 
d\omega_1 d\omega_2 U_{M_1,M_2} \Big( \frac{\omega_1 -
\omega_2}{\epsilon} \Big) \micp(M_1,M_2)  J(\omega_1,M_1)
\overline{J}(\omega_2,M_2) \nonumber \\
& & \qquad \ \   \bigg[ \epsilon^{d-2} \sum_{i \in I_q: m_i = M_1}{
\delta^{-d} \eta \Big( \frac{x_i - \omega_1}{\delta} \Big) }
\bigg]  \, \bigg[ \epsilon^{d-2} \sum_{j \in I_q; m_j = M_2}{
\delta^{-d} \eta \Big(\frac{x_j - \omega_2}{\delta} \Big)
} \bigg]  \, + \, err \big( \epsilon, \delta
\big) \nonumber
\end{eqnarray}
where the function $err$ satisfies
$$
  \lim_{\delta \downarrow 0}\limsup_{\epsilon_\downarrow 0} \mathbb{E}_N
  \big\vert err \big( \epsilon, \delta \big) \big\vert = 0,
$$
and where in the last equality, we made use of the fact that the test
functions $J$ and $\overline{J}$ take non-zero values only on
particles of a given mass, respectively $M_1$ and $M_2$.
Recalling (\ref{recp}) and (\ref{umm}). and making use of the continuity of
these test functions,
we obtain the statement in the proposition with an error term of
the form $\rm{Err}(\epsilon,\delta) = err(\epsilon,\delta)
 + err'(\epsilon,\delta)$, where the latter error $err' = O\big(
 \epsilon \delta^{-2d-1}\big)$. This completes the proof of
 Proposition \ref{szprop}.

\end{section}

\begin{section}{Uniform boundedness of the macroscopic densities}\label{unifb}

Before discussing the main goals of this section, let us review what
has been achieved
so far and what remains to be done in order to complete the proof of
the main Theorem 1.1.
Our \stoss, stated in Proposition \ref{szprop}, allows us to replace the microscopic coagulation
term with an expression
that involves $(f^\e_n*\eta^\d)( f_m^\e*\eta^\d)$ with an error that
goes to zero as $\e$
and $\d$ go to zero in this order. Here $f^\e_n$ is the microscopic density
of particles of size $n$ and $\eta^\d$ is an approximation to the identity.
   After passing to the limit $\e\to 0$,
we will have expressions of the form
\setcounter{equation}{0}
\begin{equation}
\int (f_n*\eta^\d)( f_m*\eta^\d) J \ dx+ Err(\d),
\end{equation}
where  $f_n$ now represents the macroscopic density, $J$ is a bounded
continuous test function, and
$Err(\d)$ represents an error term that goes to $0$ as $\d\to 0$. 
Recall that by the Lebesgue's theorem, the
sequence $\{f_n*\eta^\d\}$ converges to $f_n$ almost everywhere. To
derive the PDE
(1.6), we  need to show that in (4.1) we can pass to the limit
$\d\to 0$ and replace $f_n*\eta^\d$
with $f_n$. This can be achieved if we can establish the uniform
integrability of the family
\begin{equation}
\{(f_n*\eta^\d)( f_m*\eta^\d):\d>0\}.
\end{equation}
   To prepare for this, let us consider the transformation
   \begin{equation}
   q(\cdot)\mapsto g_n(dx,dt)=\e^{d-2}\sum_{i\in 
I_{q(t)}}\d_{x_i(t)}(dx)1\!\!1\big(m_i(t)=n\big)dt
   \in \mathcal M,
   \end{equation}
   where $\mathcal M$ denotes the space of measures $\mu_n(dx,dt)$ on
the set $\R^d\x [0,T]$
   such that $\mu_n\ge 0$ for
   every $n$, with
   $\mu_n\left(\R^d\times [0,T]\right)\le TZ $.
We use this transformation to define the probability measure
$\mathcal P_N$ on $\mathcal M^\N$
as the pull-back of $\mathbb P_N$. We equip $\M$ with the topology of
vague convergence
and $\mathcal M^{\N}$ with the product topology so that $\mathcal
M^\N$ is a compact metrizable space.
    We now outline our strategy for the proof
of the uniform integrability of (4.2). Let $\mathcal P$ denote a
limit point of the sequence
$\mathcal P_N$.
\begin{enumerate}
\item  We show that $\P$ is concentrated on the space of
measure sequences
$\mu=(\mu_n:n\in\N)$ such that
$\mu_1(dx,dt)=f_1(x,t)dxdt$ with $f_1(x,t)\le k_1$ almost everywhere,
where $k_1=\max_xh_1(x)$.
\item Assume that for every $l<n $, there exists a constant $k_l$
such that the measure $\P$ is concentrated on the space of $\mu$ with
$\mu_l(dx,dt)=f_l(x,t)dx dt$ for a function $f_l$ satisfying
$f_l(x,t)\le k_l$ almost everywhere
for every $l<n$. We then deduce that $\mu_n(dx,dt)=f_n(x,t)dxdt$
with a function $f_n$ weakly satisfying
\begin{equation}\label{prepde'}
\frac{\partial{f_n}}{\partial t}(x,t)  = d(n)\Delta f_n (x,t)  +
Q^n_1(f) (x,t) -Y_n,
\end{equation}
where $Q^n_1$ is defined as in (1.2) and $Y_n\ge 0$. In other words, we
can determine the form of the coagulation
gain term but all we can say about the loss term is that $Y\ge 0$.

\item  Under the assumptions of the previous step, we show that the
function $f_n$ satisfies
\[
\int_0^T \int f_n(x,t)^p\exp(-|x|) dx dt<\i  ,
\]
for every $p>2$, $\P$-almost surely.
\item  Under these same assumptions, we show that with
probability $1$
with respect to $\P$, the density $f_n$ satisfies
\begin{equation}\label{prepde}
\frac{\partial{f_n}}{\partial t}(x,t)  = d(n)\Delta f_n (x,t)  +
Q^n_1(f) (x,t) -
\hat Q^n_2(f) (x,t),
\end{equation}
where $Q^n_1$ is defined as in (1.2) and
\begin{displaymath}
\hat  Q^n_2(f) =   f_n (x,t) \sum_{m=1}^{n} 
    \macp(m,n) f_m (x,t)+X_n.
\end{displaymath}
with $X_n\ge 0$. That is, at this stage we can only determine
the form of the coagulation
loss term up to mass $n$ and for the rest $X_n$, we can only assert that $X_n\ge 0$.
\item
The equation (4.5) is good enough to deduce that there exists a
constant $k_n$ such that
   $f_n(x,t)\le k_n$, $\P$-almost surely.
\item  Applying an induction on $n$ we show that
$\mu_n(dx,dt)=f_n(x,t) dxdt$  with $f_n$ solving
(1.6) weakly and that $f_n(x,t)\le k_n$ for every $n$, $\P$-almost surely.
\end{enumerate}

The main objectives of this section are Steps 1, 3 and 5 of the above outline.
The remaining Steps 2, 4 and 6 will be carried out in Section 5.

Firstly, we establish the claim we made in Step 1.
Throughout this Section, we
let $\psi: [0,\infty) \to [0,\infty)$ denote a convex and increasing
function with $\psi(0)=0$.
\begin{definition}
For $n \in \mathbb{N}$, write
$f^{\delta}(n,x,t)=f^{\delta}(n,x;q(t))$ where
$$
f^{\delta}(n,x;q) = \epsilon^{d-2} \sum_{i \in I_q}{\eta^{\delta}(x -
x_i) 1 \! \! 1 \big\{ n_i = n \big\} },
$$
where $\eta:\R^d\to[0,\i)$ is a smooth function of compact support
with $\int\eta dx=1$ and $\eta^\d(z)=\d^{-d}\eta(z/\d)$.
\end{definition}
Note that
   $g_n*_x\eta^\d=f^\d(n,x,t)dt$ where $g_n$ was defined in (4.3).
\begin{lemma}\label{ltem}
Suppose that $\psi: [0,\infty) \to [0,\infty)$ is, in addition to its
other properties, chosen so that  $\psi'': [0,\infty) \to [0,\infty)$
exists and has compact support. Then, for any $T,s \in
[0,\infty)$ satisfying $T \geq s$,
\begin{displaymath}
    \mathbb{E}_N \int_{\mathbb{R}^d} \psi \big( f^{\delta}(1,x,T) \big) dx
    \leq  \mathbb{E }_N \int_{\mathbb{R}^d} \psi \big( f^{\delta}(1,x,s)
    \big) dx \, + \, \big( T - s \big)
    \frac{\epsilon^{d-2}}{\delta^{d+2}} \int_{\mathbb{R}^d}{\big| \nabla
    \eta \big|^2 dx }.
\end{displaymath}
\end{lemma}
{\bf Proof}
  From the identity ${\mathbb E}_N F(q(T))={\mathbb E}_N
F(q(s))+{\mathbb E}_N\int_s^T
{\mathbb L} F(q(t))dt$, we deduce that
\begin{equation}\label{abi}
    \mathbb{E}_N \int_{\mathbb{R}^d} \psi \big( f^{\delta}(1,x,T) \big) dx
     =   \mathbb{E}_N \int_{\mathbb{R}^d} \psi \big( f^{\delta}(1,x,s)
    \big) dx  + A_1 +  \sum_{n=1}^{\infty}A_{1,n},
\end{equation}
where
\begin{displaymath}
    A_1 =  d(1) \mathbb{E}_N \int_s^T dt \int_{\mathbb{R}^d} \sum_{i \in
I_{q(t)}}{\Delta_{x_i}  \psi \big( f^{\delta}(1,x,t)
    \big) 1 \! \! 1 \big\{ m_i = 1 \big\} dx }
\end{displaymath}
and where, for example,
\begin{displaymath}
    A_{1,1} =   \micp(1,1) \mathbb{E}_N\int_s^T dt \int_{\R^d}\,  \sum_{i,j \in
I_{q(t)}} V_{\epsilon}(x_i-x_j)
    D_{i,j}(x,t)   1 \! \! 1 \Big\{
    m_i = m_j = 1 \Big\}dx,
\end{displaymath}
where $D_{i,j}(x,t)$ equals
\[
\psi\left(f^\d(1,x,t)-\e^{d-2}\eta^{\delta} ( x - x_i)
-\e^{d-2}\eta^{\delta}( x - x_j )\right)-\psi\left(f^\d(1,x,t)\right).
\]
  From the monotonicity of $\psi$, we know that $D_{i,j}\le 0$. This in
turn implies that, for each $n \in \mathbb{N}$, we have that
$A_{1,n} \le 0$.

In bounding $A_1$, we make use of the equality
\begin{equation}\label{abg}
\Delta_{x_i} \psi \big( f^{\delta} \big) = \psi''  \big(
f^{\delta} \big) \big\vert \nabla_{x_i} f^{\delta} \big\vert^2
+  \psi'  \big(
f^{\delta} \big)  \Delta_{x_i} f^{\delta}.
\end{equation}
Since $\psi''$ is of compact support, we find that, for any given $t
\in [s,T]$,
\begin{eqnarray}\label{abf}
    & & \int_{\mathbb{R}^d}   \mathbb{E}_N \sum_{i \in I_{q(t)}}{\psi''
\big( f^{\delta}(1,x,t)
    \big)  \big\vert \nabla_{x_i} f^{\delta} \big( 1, x,t \big)
    \big\vert^2 1 \! \! 1 \big\{ m_i = 1 \big\} dx } \\
    & \leq & C\mathbb{E}_N  \int_{\mathbb{R}^d} \sum_{i \in I_q} \big\vert
\nabla_{x_i} f^{\delta} \big( 1, x,t \big)
    \big\vert^2 1 \! \! 1 \big\{ m_i = 1 \big\} \nonumber \\
    & = & C\frac{\epsilon^{2(d-2)}}{\delta^{2(d+1)}} \int_{\mathbb{R}^d} dx
    \sum_{i \in I_q} \big| \nabla \eta \big|^2 \Big( \frac{x -
    x_i}{\delta} \Big) 1 \! \! 1 \big\{ m_i = 1 \big\}
    \leq  C \frac{\epsilon^{d-2}}{\delta^{d+2}} \int_{\mathbb{R}^d} \big|
    \nabla \eta \big|^2 (x) dx, \nonumber
\end{eqnarray}
where in the last inequality, we used the fact that the total number
of particles in the model is bounded above by $Z \epsilon^{2-d}$.
Note also that, for $t \in [s,T]$,
\begin{eqnarray}\label{abe}
    & & \int_{\mathbb{R}^d} \mathbb{E}_N \sum_{i \in I_q}{\psi' \big(
f^{\delta}(1,x,t)
    \big)  \Delta_{x_i} f^{\delta} \big( 1, x,t \big) 1 \! \! 1\big\{ m_i
    = 1 \big\} dx } \\
    & = &  \int_{\mathbb{R}^d} \mathbb{E}_N \psi' \big( f^{\delta}(1,x,t)
    \big)  \Delta_x f^{\delta} \big( 1, x,t \big) dx
    = -  \int_{\mathbb{R}^d} \mathbb{E}_N \psi'' \big( f^{\delta}(1,x,t)
    \big)  \big\vert \nabla_x f^{\delta} \big( 1, x,t \big) \big\vert^2
    dx \leq 0 \nonumber
\end{eqnarray}
Applying the equality (\ref{abg}) to the term $A_1$ in (\ref{abi}) and
using the bounds (\ref{abf}) and (\ref{abe}), we deduce the statement
of the Lemma, each of the terms $\big\{ A_{1,n}: n \in \mathbb{N}
\big\}$ being non-positive. $\Box$ \\

We are now ready to establish the claim of the first step in the outline.
\begin{lemma}\label{step1}
Let the measure $\P$ be a limit point of the sequence $\{\P_N\}$.
Then the measure $\P$ is concentrated on the space of $\mu$ such that
$\mu_1(dx,dt)=f_1(x,t)dxdt$ with $f_1(x,t)\le k_1$ almost everywhere,
where $k_1=\max_xh_1(x)$.
\end{lemma}
{\bf Proof} Fix a continuous function $\g:[0,\i)\to [0,\i)$ of
compact support with $\int\g dt=1$.
  It is straightforward to
show that the map
\begin{equation}\label{ilsc}
   \nu\mapsto F(\nu):=\int_{\mathbb{R}^d} \psi \left(
\int\eta^\d(x-y)\g(t)\nu(dy,dt)\right)
   dx,
\end{equation}
is lower semicontinuous with respect to the vague topology.
By Lemma 4.1,
\begin{equation}\label{lrhs}
\int {\mathcal P}_N(d\mu) \int_{\mathbb{R}^d} \psi \left(
\int\eta^\d(x-y)\g(t)\mu_1(dy,dt)\right)
   dx\nonumber
  \leq \mathbb E_N \int \psi\left(f^\d(1,x,0)\right) dx
+C\frac{\e^{d-2}}{\d^{d+2}},
\end{equation}
where $\psi$ satisfies the conditions of Lemma 4.1.
Choose $k \in \mathbb{R}$ to satisfy $k >k_1=\max_x h_1(x)$. By
approximation, we may choose
   $\psi: \mathbb{R} \to [0,\infty)$  to be $\psi
(z) = (z-k)_{+}$. With this choice, we can show that the right-hand-side of
(\ref{lrhs}) goes to $0$ in the limit $N\to\i$. The left-hand-side
being non-negative, it
has zero limit also. Using the lower semicontinuity of (\ref{ilsc}) we deduce
\[
\int {\mathcal P}(d\mu)\int_{\mathbb{R}^d} \psi \left(
\int\eta^\d(x-y)\g(t)\mu_1(dy,dt)\right)
   dx =0.
   \]
Since the function  $\g$ is arbitrary we deduce that ${\mathcal
P}(d\mu)$-almost surely,
the measure $\mu_1$ has a density $f_1$ and that
$f_1(x,t) \leq k_1$, for almost all $(x,t) \in
\mathbb{R}^d\x [0,\i)$.   ${\Box}$\\

We now turn to Steps 3 and 5 of the outline. The second part of Lemma 
4.3 below will take care of
Step 5. A microscopic version of the first part of Lemma 4.3 is out 
Step 5 and will be carried out in
Lemma 4.4. Indeed Lemma 4.3
is a variation of a theorem of Wrzosek \cite{wrz} regarding the 
boundedness of the solutions of (1.6).
We give a detailed proof of Lemma 4.3 even though the proof is a straightforward
adaptation of the arguments in \cite{wrz}.
\begin{lemma}\label{lemunifb}
Let $J:\R^d\to (0,\i)$ be a smooth function with $\int J\ dx<\i$ such 
that for a constant $R>\frac 12$,
we have that $|\nabla J(x)|\le RJ(x)$
for every $x$.
\begin{itemize}
\item Assume that $\big\{ f_j: \mathbb{R}^d \times [0,\infty) \to
[0,\infty), j \in \{1,2,\dots,n\} \big\}$ denotes a collection of functions
satisfying (4.4) with
$Y$ nonnegative. Assume also that
there exist constants $\{k_1,\dots,k_{n-1}\}$ such that
$\sup_{x,t}f_n(x,t)\le k_j$ for $j<n$. Then for every $p>2$,
\begin{equation}\label{unbo'}
\sup_t \int f_n(x,t)^p J(x)dx <\i.
\end{equation}
\item
In addition to the above assumptions, suppose that $f_n$ satisfies (4.5).
Then
\begin{equation}\label{unbo}
\sup_{x,t}f_n(x,t)=:k_n<\i.
\end{equation}
\end{itemize}
\end{lemma}
{\bf Proof} We first establish (\ref{unbo'}).
Let us write $\psi_p(z)$ for $z^p$.
Note that
\begin{eqnarray}
   & & \frac{d}{dt} \int_{\mathbb{R}^d} \psi_p \big( f_n(t,x) \big)
J(x) dx \nonumber \\
   & = & \int_{\mathbb{R}^d} \psi'_p \big( f_n(t,x) \big) \bigg[ d(n) \big(
   \Delta f_n \big) (t,x)  \nonumber \\
   & & \qquad \qquad \quad  + \frac 12\sum_{m=1}^{n-1}{\hat{\alpha}(m,n-m)
f_m(t,x) f_{n-m}(t,x)} \, - \, Y \bigg] J(x) dx , \nonumber
\end{eqnarray}
We have that
\begin{eqnarray}
   & & \int_{\mathbb{R}^d} \psi'_p \big( f_n(t,x) \big) \big( \Delta f_n \big)
   \big( t,x \big) J(x) dx \label{werone} \\
   & = & - \int_{\mathbb{R}^d} \psi''_p  \big( f_n(t,x) \big) \big\vert
   \nabla f_n \big\vert^2 \big( t,x \big) J(x) dx \nonumber \\
   & & - \int_{\mathbb{R}^d} \psi'_p  \big( f_n(t,x) \big)
   \nabla f_n  \big( t,x \big) \cdot \nabla J(x) dx. \nonumber
\end{eqnarray}
We find that
\begin{eqnarray}
   & & \Big\vert \int_{\mathbb{R}^d} \psi'_p  \big( f_n(t,x) \big)
   \nabla f_n \big( t,x \big) \cdot \nabla J(x) dx \Big\vert
   \nonumber \\
    & \leq & R \int_{\mathbb{R}^d} \frac{ \psi'_p  \big( f_n(t,x)
    \big) }{\psi_p'' \big( f_n(t,x) \big)}
   \big\vert \nabla f_n \big( t,x \big) \big\vert \psi''_p \big(
f_n(t,x) \big)  J(x)  dx
   \nonumber \\
    & \leq & \frac{R}{2} \int_{\mathbb{R}^d} R^2 \frac{ \psi'_p  \big( f_n(t,x)
    \big) ^2}{\psi''_p \big( f_n(t,x) \big)}
    J(x) dx \, + \,  \frac{R}{2} \int_{\mathbb{R}^d}
   \frac{\big\vert  \nabla f_n (t,x)
     \big\vert^2}{R^2} \psi''_p \big( f_n(t,x) \big)
    J(x)  dx
   \nonumber \\
   & \leq & {R^3} \int_{\mathbb{R}^d}   \psi_p  \big(
   f_n(t,x) \big)  J(x) dx \, + \, \frac{1}{2R}
   \int_{\mathbb{R}^d}{\big\vert  \nabla   f_n (t,x)
    \big\vert^2 \psi''_p \big( f_n(t,x) \big) J(x) dx}. \label{wertwo}
\end{eqnarray}
In the first inequality, we applied the bound $\big\vert \nabla J
\big\vert \leq R J(x)$. In the third, we made use of the inequality
\begin{displaymath}
   \frac{\big\vert \psi'_p \big\vert^2}{\psi_p''} 
=\frac{p}{p-1}\psi_p\leq 2  \psi_p
.\end{displaymath}

  From (\ref{werone}) and (\ref{wertwo}), we learn that
\begin{eqnarray}
   \int_{\mathbb{R}^d} \psi'_p \big( f_n(t,x) \big) \big( \Delta f_n \big)
   \big( t,x \big) J(x) dx & \leq & \Big( \frac{1}{2R} - 1 \Big)
\int_{\mathbb{R}^d} \psi''_p  \big( f_n(t,x) \big) \big\vert
   \nabla f_n \big\vert^2 \big( t,x \big) J(x) dx   \nonumber \\
   & & +  R^3 \int_{\mathbb{R}^d}   \psi_p  \big(
   f_n(t,x) \big) J(x)dx \nonumber \\
      & \leq &   R^3 \int_{\mathbb{R}^d}   \psi_p  \big(
   f_n(t,x) \big) J(x)dx.\label{neweq}
\end{eqnarray}
We find that
\begin{eqnarray}
   & & \frac{d}{dt} \int_{\mathbb{R}^d} \psi_p \big( f_n(t,x) \big)
J(x) dx \nonumber \\
   & \leq &\frac 12 \sum_{m=1}^{n-1} \hat{\alpha} \big( m,n-m \big)
\int_{\mathbb{R}^d} \psi'_p \big( f_n(t,x) \big)
   f_m(t,x) f_{n-m}(t,x)  J(x) dx
   \nonumber \\
   & & +  \, R^3 d(n) \int_{\mathbb{R}^d}  \psi_p 
\big( f_n(t,x) \big)
    J(x) dx
   \nonumber \\
   & \leq  &   A  \int_{\mathbb{R}^d} \psi'_p \big(
   f_n(t,x) \big) J(x) dx \nonumber \\
   & &  +  R^3 d(n)\int_{\mathbb{R}^d}  \psi_p  \big( f_n(t,x) \big)
    J(x) dx,
   \nonumber
\end{eqnarray}
where the uniform bounds $\sup f_i 
\leq k_i, i \in \{1,\ldots, n-1 \}$
provided by the hyptohesis, were applied in the second
inequality, and where the constant $A$ was defined by
\[
A=\frac 12\sum_{m=1}^{n-1}
   \hat{\alpha} \big( m,n-m \big) k_m k_{n-m}.
   \]
   Note that
$$
\psi'_p (z) \leq {p}\psi_p(z)+p,
$$
implying that
\begin{eqnarray}
   & & \frac{d}{dt} \int_{\mathbb{R}^d} \psi_p \big( f_n (t,x) \big)
J(x) dx \nonumber \\
   & \leq & B \int_{\mathbb{R}^d} \psi_p \big( f_n (t,x) \big) J(x)dx \
+ C, \nonumber
\end{eqnarray}
where
$$
B = {p A} + R^3 d(n),$$
and
$$
C = p  A \int_{\mathbb{R}^d} J(x) dx.$$
This implies that
$$
\sup_{t \in [0,T]}{ \int_{\mathbb{R}^d} \psi_p \big( f_n (t,x)
\big)J(x) dx } < \infty,
$$
and concludes the first part of the Lemma.

We now turn to the second part of the Lemma.
Note that
\begin{equation}\label{nsat}
   \frac{\partial}{\partial t} f_n \big( x,t \big)   =  d(n) \Delta f_n
   \big(x,t \big)  + \frac 12\sum_{m=1}^{n-1}
   \hat{\alpha} \big( m,n-m \big) f_m \big( x,t \big) f_{n-m} \big( x,t \big) \,
   - \hat{\alpha}\big(n,n\big) f_n \big( x,t\big)^2 \, - \, \hat X,
\end{equation}
where $\hat X \geq 0$. We now choose $\psi_p: [0,\infty) \to
[0,\infty)$ according to  $\psi_p (z) = \big[(z-1)^+\big]^p$.
A repetition of the proof of (\ref{neweq}) yields
\begin{eqnarray}
   \int_{\mathbb{R}^d} \psi'_p \big( f_n(t,x) \big) \big( \Delta f_n \big)
   \big( t,x \big) J(x) dx & \leq & \Big( \frac{1}{2R} - 1 \Big)
\int_{\mathbb{R}^d} \psi''_p  \big( f_n(t,x) \big) \big\vert
   \nabla f_n \big\vert^2 \big( t,x \big) J(x) dx   \nonumber \\
   & & +  R^3 \int_{\mathbb{R}^d}   \psi_p  \big(
   f_n(t,x) \big) J(x)dx \nonumber \\
      & \leq &   R^3 \int_{\mathbb{R}^d}   \psi_p  \big(
   f_n(t,x) \big) J(x)dx.
\end{eqnarray}
It follows that
\begin{eqnarray}
   & &
\frac{d}{dt} \int_{\mathbb{R}^d} \psi_p \big( f_n (x,t) \big) J(x) dx
\label{befdefa} \\
   & \leq & A \int_{\mathbb{R}^d} \psi'_p
   \big( f_n (x,t) \big) J(x) dx +d(n)R^3 \int_{\mathbb{R}^d} \psi_p
   \big( f_n (x,t) \big) J(x) dx \nonumber\\
   &&- \hat{\alpha} \big(n,n\big)
\int_{\mathbb{R}^d} f_n(x,t)^2  \psi'_p
   \big( f_n (x,t) \big) J(x) dx. \nonumber
\end{eqnarray}
Note that
$$
   \psi'_p (a) \leq r \psi_p(a) + \frac{p^p}{r^{p-1}},
$$
implying that
\begin{eqnarray}
    \frac{d}{dt} \int_{\mathbb{R}^d} \psi_p \big( f_n(x,t) \big) J(x)
   dx  & \leq &  \big(A r+d(n)R^3\big) \int_{\mathbb{R}^d}
   \psi_p \big( f_n (x,t) \big) J(x) dx   \label{lop} \\
    & &  - \, \hat{\alpha}\big(n,n\big) \int_{\mathbb{R}^d}
   f_n (x,t)^2  \psi'_p \big( f_n (x,t) \big) J(x) dx. \nonumber
   \\
   && +  A \frac{p^p}{r^{p-1}}
   \int_{\mathbb{R}^d}{J(x)dx}.\nonumber
\end{eqnarray}
On the other hand, $a^2\psi'(a)\ge p\psi(a)$.
Hence, the right-hand-side of (\ref{lop}) is bounded above by
$$
   \left(A r+d(n)R^3\right) \int_{\mathbb{R}^d}
   \psi_p \big( f_n \big(x,t\big) J(x)dx   +  A \frac{p^p}{r^{p-1}}
   \int_{\mathbb{R}^d}{J(x)dx} \, - \,  p \hat{\alpha}\big(n,n\big)
\int_{\mathbb{R}^d}
    \psi_p \big( f_n \big( x,t \big) \big) J(x)dx.
$$
It suffices to establish (\ref{unbo'}) for sufficiently large $p$.
We choose $p$ large enough so that our choice $r = A^{-1}\big( p 
\hat{\alpha}\big(n,n\big)-d(n)R^3)$
for $r$ is positive. Such a choice yields
$$
\frac{d}{dt} \int_{\mathbb{R}^d} \psi_p \big( f_n (x,t) \big) J(x) dx \leq
\frac{A^p p^p}{\left(p\hat{\alpha}\big(n,n\big)-d(n)R^3\right)^{p-1}} 
\int_{\mathbb{R}^d}{J(x)dx} .
$$
  From this, it follows that
$$
\int_{\mathbb{R}^d} \psi_p \big( f_n (x,T) \big) J(x) dx  \leq
\int_{\mathbb{R}^d} \psi_p \big( h_n (x) \big) J(x) dx \, + \, T
\frac{A^p p}{\left(\hat{\alpha}\big(n,n\big)-d(n)R^3\right)^{p-1}} 
\int_{\mathbb{R}^d}{J(x)dx}.
$$
By taking the $p$-th root of this inequality and then taking the value
of $p$ to infinity, we deduce that
\begin{equation}\label{ghg}
   \sup_{x,t} f_n(x,t) \leq \sup_x h_n(x)+
\frac{A}{\hat{\alpha}\big(n,n\big)}.
\end{equation}
  $\Box$

   We now concentrate on the third step of the outline.
   Recall the random measure $g_n$ that was defined by (4.3).
   For our purposes we need to regard the microscopic density
   as a member of the smaller space 
$\L=L^\i\big([0,T];\M_Z\big(\R^d\big)\big)\subset \M$
   where $\M_Z\big(\R^d\big)$
   denotes the space of nonnegative measures $\g$ with 
$\g\big(\R^d\big)\le Z$. More precisely, the measure
   $g_n(dx,dt)\in \L$ is absolutely continuous with respect to the 
Lebesgue measure for almost all $x$
   and if $g_n(dx,dt)=g_n(dx,t)dt$, then $g_n\big(\R^d,t\big)\le Z$ 
for almost all $t\in [0,T]$ and every
   $n$. Since we have a uniform bound $Z$ on the total measure of 
$g(dx,t)$ for almost every $t$,
   onee can readily show that indeed $\L$ is a closed subset of the 
compact metric
   space $\M$.
     Consider the transformation
   \begin{equation}\label{aug}
   q\mapsto \big(g(dx,t),X(dx,dt)\big)\in \hat \M =\L^{\N}\x \mathcal N^\N,
   \end{equation}
   where the random measures $(g,X)=(g_n,X_n:n\in\N)$  are given by
\begin{eqnarray}
g_n(t,dx)& =&\frac{\e^{d-2}}2 \sum_{i \in I_{q(t)}}\d_{x_i(t)}(dx)
1\!\!1\big(m_i(t)=n\big),\\
X_n(dx, dt)& =&\e^{d-2} \sum_{i,j \in I_{q(t)}}\micp(m_i(t),m_j(t))
V_{\epsilon}(x_i(t) - x_j(t))
1\!\!1 \big(m_i(t)+m_j(t)=n\big)\nonumber\\
& &\qquad\qquad \frac 
1n\left(m_i(t)\d_{x_i(t)}+m_j(t)\delta_{x_j(t)}\right)(dx)\ 
dt,\nonumber
\end{eqnarray}
and $\mathcal N$ denotes the space of nonnegative measures on $\R^d\x [0,T]$.
The transformation (\ref{aug}) defines a probability measure
$\hat\P_N$ for each $N$ on
the space $\hat\M$. This measure is an augmentation of the measure
$\P_N$ that is defined on
the space $\M^\N$. We equip $\hat \M$ with the product toplogy with
$\L$ and $\cN$ both equipped
with the vague topology. Certainly $\L^\N$ is a compact metric space
but $\cN$ is only a
separable complete metric space. As an immediate consequence of Lemma 3.1,
we learn that the sequence $\big\{\hat\P_N\big\}$ is indeed tight.
More precisely, if $\cN_l$
denotes the space of measures $X\in \cN$ such that $X\big(\R^d\x 
[0,T]\big)\le l$, then
$\sup_N \P_N\big(\cN_l^c\big)\ge 1-C/l$, for a constant $C$.

\begin{lemma}\label{huno}
Let $J$ be as in Lemma 4.3.
Let $\hat \P$ be a limit point of the sequence $\big\{\hat\P_N\big\}$
and write $\P$ for the
$g$-marginal of $\hat\P$.
Assume that there exist constants $k_1,k_2,\dots,k_{n-1}$,
such that the measure $ \P$ is concentrated on the set of measure
sequences $g=(g_n:n\in\N)$
such that $g_j(dx,t)=f_j(x,t)dxdt$ with $\sup_{x,t}f_j(x,t)\le k_j$
for $j=1,2,\dots,n-1$.
Then $g_n(dx,dt)=f_n(x,t)dxdt$ for $\P$-almost surely and
$$
\int \left\{ \sup_{t \in [0,T]}  \int_{\mathbb{R}^d}
f_n(x,t)^p J(x)dx \right\}\ \P(dg)< \infty.
$$
Moreover,  we have that $X_n(dx,dt) = \g_n(x,t) dx dt$, with
$$
\sup_{(x,t)} \g_n(x,t) \leq
    \sum_{m=1}^{n-1} \macp \big(m,n-m\big) ,$$
$\hat{\mathcal P}$-almost surely.
\end{lemma}
{\bf Proof}
The statement of the Lemma asserts a microscopic analogue of the
first part of Lemma 4.3.
Again we write $\psi_p(z)$ for $z^p$.
Let $\ell:\R\to [0,\i)$ be a smooth function of compact support
with $\int\ell=1$ and set $\hat f^\d(n,x,t)=\int_0^\i f^\d(n,x,s)\ell(t-s)ds$.
A straightforward calculation yields,
\begin{eqnarray}
\frac {d}{dt}\hat f^\d(n,x,t)\le& -&f^\d(n,x,0)\ell(t)+\int_0^\i d(n)
\Big( \sum_{i \in I_q} \Delta_{x_i} f^{\delta}
    (n,x,s) \Big)\ell(t-s)ds \\
    &+&\int_0^\i \frac {\e^{d-2}}2\sum_{m=1}^{n-1} {\micp} \big( m,n-m
    \big)\sum_{i,j \in I_{q(s)}}
    V_{\epsilon} \big( x_i(s) - x_j(s) \big)\nonumber\\
  & &  \quad\quad
    \left(\frac{m_i}n\eta^{\delta} \big( x - x_i(s) 
\big)+\frac{m_j}n\eta^{\delta} \big( x - x_j(s)\big)
    \right)
\nonumber\\
   & &\quad\quad 1 \! \! 1 \big\{ m_i(s) = m, m_j(s) = n-m \big\}
\ell(t-s)ds\nonumber\\
    &+&\int_0^\i M_s(x)\ell'(t-s)ds.\nonumber
    \end{eqnarray}
Here we have an inequality because we only used the gain term in applying
the coagualtion operator $\mathbb A_C$, and $M_s$ is a martingale satisfying
   \begin{displaymath}
\mathbb{E}_N \Big[  M_s(x)^2 \Big] =  \mathbb{E}_N\int_{0}^{s}{A_1 \big( q(t)
,x
\big) dt} +  \mathbb{E}_N\int_{0}^{s}{A_2 \big( q(t),x \big) dt},
\end{displaymath}
where $A_1(q,x)$ and $A_2(q,x)$ are respectively set equal to
\begin{displaymath}
A_1 \big( q ,x\big) = \e^{2(d-2)} \sum_{i \in I_q}  d (n)
    \big\vert \nabla_{x} \eta^\d \big\vert^2 \big(x- x_i
\big)1\!\!1\big(m_i=n\big),
\end{displaymath}
and
\begin{eqnarray}
    A_2 \big( q ,x\big)= &  & \frac {\epsilon^{2(d-2)} }2\sum_{i,j \in I_q}  \micp
    \big( m_i , m_j \big) V_{\epsilon} \big( x_i - x_j \big)
    \\
     & & \qquad \bigg[ \frac{m_i}{m_i + m_j} K \big(x- x_i, m_i + m_j \big)
    +  \frac{m_j}{m_i + m_j} K \big(x- x_j, m_i + m_j \big)
      \nonumber \\
     & & \qquad \qquad \quad   - K \big(x- x_i, m_i \big) - K \big(
     x- x_j,m_j\big) \bigg]^2, \nonumber
\end{eqnarray}
where $K(x-y,m)=\eta^\d(x-y)1\!\!1(m=n)$. We can readily show,
\begin{eqnarray}
\int_{\mathbb{R}^d} A_1 \big( q,x \big) dx &\leq& C \epsilon^{2(d-2)}
\d^{-2d-2}\sum_{i \in I_q}  d (n)
1\!\!1(m_i=n)
\leq C\e^{d-2}
\d^{-d-2}, \\
\mathbb E_N\int_0^s\int_{\mathbb{R}^d}  A_2 \big( q(t),x \big)dxdt  &\leq&
C\epsilon^{2(d-2)}\d^{-2d} \mathbb E_N\int_0^s\sum_{i,j \in I_{q(t)}}  \micp
    \big( m_i(t) , m_j(t) \big) V_{\epsilon} \big( x_i(t) - x_j(t) \big)
\nonumber\\
& &\qquad\qquad\quad 1\!\!1(m_i(t)+m_j(t)\le n)dt\nonumber\\
&\le& C\e^{d-2}
\d^{-2d},\nonumber
     \end{eqnarray}
    where we have used Lemma 3.1 for the last inequality.

As a result,
\begin{equation}\label{mart}
   \mathbb{E}_N\int_{\mathbb{R}^d} \Big[  M_s(x)^2 \Big]dx \le 
C\e^{d-2}\d^{-2d-2}.
   \end{equation}
   We certainly have,
   \begin{displaymath}
   \int_{\mathbb{R}^d}\frac{d}{dt} \psi_p \big(\hat f^{\delta}(n,x,t) \big) J(x)
dx \le
B_1+ B_2+B_3+B_4,
\end{displaymath}
where
\begin{eqnarray}
B_1&=& -\ell(t)\int_{\R^d}\psi'_p \big(\hat f^{\delta}(n,x,t) \big)
f^\d(n,x,0)J(x)
dx \\
    B_2& =& d(n)  \int_{\mathbb{R}^d} \psi'_p \big(
    \hat f^{\delta} (n,x,t) \big) \Big( \sum_{i \in I_q} \Delta_{x_i} 
\hat f^{\delta}
    (n,x,t) \Big) J(x) dx,\nonumber
\end{eqnarray}
   \begin{eqnarray}
    B_3 & =  &\epsilon^{d-2}  \int_0^\i ds\ \ell(t-s)  \int_{\mathbb{R}^d}
\psi'_p \big(
    \hat f^{\delta} (n,x,t) \big)   \sum_{m=1}^{n-1} {\micp} \big( m,n-m
    \big) \nonumber \\
    & & \qquad  \sum_{i,j \in I_{q(s)}}
    V_{\epsilon} \big( x_i(s) - x_j(s) \big)  \eta^{\delta} \big( x-x_i(s) \big)
    1 \! \! 1 \big\{ m_i(s) = m, m_j(s) = n-m \big\}  J(x) dx,
\end{eqnarray}
and
\begin{eqnarray}
\mathbb{E}_NB_4&\le& \mathbb{E}_N \int_0^\i ds\ \ell'(t-s)  \int_{\mathbb{R}^d}
\psi'_p \big(
    \hat f^{\delta} (n,x,t) \big) |M_s(x)|J(x)dx \\
    &\le& C\d^{-pd+d} \int_0^\i ds\ \ell'(t-s) \left\{
\int_{\mathbb{R}^d}\mathbb{E}_N
\big[M_s(x)\big]^2J(x)dx\right\}^{1/2}\nonumber\\
&\le&C \ell(t)\e^{d/2-1}\d^{-1-pd},
    \nonumber
\end{eqnarray}
where we used $\hat f^\d\le Z\d^{-d}$ for the the second inequality.

On the other hand,
\begin{displaymath}
    B_2 =  d(n)    \int_{\mathbb{R}^d} \psi'_p \big(
    \hat f^{\delta} (n,x,t) \big) \Delta_x \hat f^{\delta}
    (n,x,t)  J(x) dx  = B_{21} + B_{22},
\end{displaymath}
where
\begin{displaymath}
    B_{21} =  - d(n)  \int_{\mathbb{R}^d} \psi''_p \big(
    \hat f^{\delta} (n,x,t) \big) \big\vert \nabla_x \hat f^{\delta}
    (n,x,t) \big\vert^2  J(x) dx
\end{displaymath}
and
\begin{displaymath}
    B_{22} =  - d(n) \int_{\mathbb{R}^d} \psi'_p \big(
    \hat f^{\delta} (n,x,t) \big) \nabla_x \hat f^{\delta}
    (n,x,t) \cdot \nabla_x  J(x) dx.
\end{displaymath}
Note that
\begin{eqnarray}
    \big\vert B_{22} \big\vert & \leq &  R d(n)
    \int_{\mathbb{R}^d} \frac{ \psi'_p \big(
    \hat f^{\delta} (n,x,t) \big) }{\psi''_p\big(
    \hat f^{\delta} (n,x,t) \big) } \big\vert \nabla_x \hat f^{\delta}
    (n,x,t) \big\vert  \psi''_p\big(
    \hat f^{\delta} (n,x,t) \big)  J(x)  dx \\
    & \leq &  \frac{R d(n)}{2}
    \int_{\mathbb{R}^d} R^2 \frac{\big\vert \psi'_p \big(
    \hat f^{\delta} (n,x,t) \big) \big\vert^2}{\psi''_p\big(
    \hat f^{\delta} (n,x,t) \big) }  J(x)  dx \nonumber \\
     & & + \,  \frac{R d(n)}{2}
    \int_{\mathbb{R}^d} \frac{1}{R^2} \big\vert \nabla_x \big(
    \hat f^{\delta} (n,x,t) \big) \big\vert^2 \psi''_p\big(
    \hat f^{\delta} (n,x,t) \big)   J(x)  dx \nonumber
\end{eqnarray}
It follows that
\begin{eqnarray}\label{cvt}
   \left| B_2\right| & \leq & \Big( 1 - \frac{1}{2R} \Big) B_{21} +
    {R^3}     \int_{\mathbb{R}^d}  \psi_p \big(
     \hat f^{\delta} (n,x,t) \big)   J  (x)   dx
    \\
    & \leq &   {R^3}     \int_{\mathbb{R}^d}
\psi_p \big(
    \hat f^{\delta} (n,x,t) \big)   J  (x)   dx , \nonumber
\end{eqnarray}
where in the second inequality, the fact that $B_{21} \leq 0$ was used.

Note that
\begin{displaymath}
    B_{3} \leq   \int_0^\i ds\ \ell(t-s) \int_{\mathbb{R}^d}
\psi'_p \big( \hat f^{\delta}(n,x,t) \big)
     X^{\delta}(x,n;q(s))J(x) dx
\end{displaymath}
where
\begin{eqnarray}
X^{\delta}(x,n;q)& =& \frac{\epsilon^{d-2}}2 \sum_{m=1}^{n-1} \micp \big( m,n-m
\big) \sum_{i,j \in I_q}V_{\epsilon}
    \big( x_i - x_j \big)   1 \! \! 1 \big\{ m_i= m, m_j= n-m \big\} \\
    & &\qquad\qquad\quad\left(\frac {m_i}n
    \eta^{\delta}(x - x_i)+\frac {m_j}n
    \eta^{\delta}(x - x_j)\right) .
\nonumber
\end{eqnarray}
In summary, the expression
\[
   \frac{d}{dt}\int_{\mathbb{R}^d} \psi_p \big(\hat f^{\delta}(n,x,t) \big) J(x)
dx
\]
is bounded above by,
\begin{eqnarray}
\label{summa}
   & & -\ell(t)\int_{\R^d}\psi'_p \big(\hat f^{\delta}(n,x,t) \big)
f^\d(n,x,0)J(x)
dx \\
& &\  +R^3    \int_{\mathbb{R}^d}
\psi_p \big(
    \hat f^{\delta} (n,x,t) \big)   J  (x)   dx\nonumber\\
    & &\ \  +  \int_0^\i ds\ \ell(t-s) \int_{\mathbb{R}^d}
\psi'_p \big( \hat f^{\delta}(n,x,t) \big)
     X^{\delta}(x,n;q(s))J(x) dx
+\hat M^\e, \nonumber
        \end{eqnarray}
with $\hat M^\e$ satisfying
\begin{equation}\label{martin}
\mathbb{E}_N\left|\hat M^\e\right|\le C\e^{d/2-1}\d^{-1-pd}.
\end{equation}
Before sending$\e\to 0$ in (\ref{summa}), recall the space
$\cN_l$ which was defined right before the Lemma and
let us observe that the functional
$\G:\M\x \cN_l\to [0,\i)$, defined by
$$
\Gamma ( g , X) = \int_0^\i\iint \psi' \left( \int_0^\i \int
   \eta^{\delta}(x-y)\ell(t-\theta) g(dy,d\theta)\right)
\eta^{\delta}(x-z)J(x)\ell(t-s)X(dz,ds)dx,
$$
is continuous with respect to the weak topology.
Indeed, if $g^r\rightarrow g$ and $X^r\rightarrow X$,
then $F^r\to F$ and $Y^r\to Y$, where
\begin{eqnarray}
F(x)&=&\int_0^\i \int
   \eta^{\delta}(x-y)\ell(t-\theta) g(dy,d\theta) \\
   F^r(x)&=&\int_0^\i \int
   \eta^{\delta}(x-y)\ell(t-\theta) g^r(dy,d\theta)\nonumber\\
  Y(x)&=&\int_0^\i \int
   \eta^{\delta}(x-y)\ell(t-\theta) X(dy,d\theta)\nonumber\\
   Y^r(x)&=&\int_0^\i \int
   \eta^{\delta}(x-y)\ell(t-\theta) X^r(dy,d\theta).\nonumber
   \end{eqnarray}
   Since both sequences $\{F^r\}$ and $\{Y^r\}$ are uniformly bounded, 
we deduce that
   $$
   \lim_{r\to\i}\int\psi'(F^r(x))Y^r(x) J(x)dx=\int\psi'(F(x))Y(x) J(x)dx.
   $$

We now fix $\d$ and send $\e$ to $0$.
  From (\ref{summa}) we deduce that the measure  $\hat{\mathcal P}$
     is concentrated on the space of $\big(g_n , X_n \big)$ for
    which
\begin{equation}\label{pbd}
\frac{d}{dt} \int{} \psi'_p \big( \hat f^{\delta}_n \big) J dx\le -\ell(t)
\int \psi_p \big( \hat f^{\delta}_n \big) h_n^\d J dx
+R^3\int{} \psi_p \big( \hat f^{\delta}_n \big) J dx+\int{} \psi_p'
    \big(\hat f^{\delta}_n \big) X^\d_n Jdx,
\end{equation}
where,
\[
   \hat f^\d_n =g_n*_x\eta^\d*_t\ell,\ \ \   h_n^\d=h_n*_x\eta^\d,\ \ \
   X_n^\d=X*_x\eta^\d*_t\ell.
   \]
   We now claim that
\begin{equation}\label{mcla}
X_n^\d \leq \sum_{m=1}^{n-1} \macp \big( m,n-m \big) k_m k_{n-m},
\end{equation}
$\hat \P$--almost surely. By the \stoss, as it is recorded
in (\ref{szrec}),
\begin{eqnarray}
  & &
\lim_{\delta \to 0}\limsup_{N\to \i}{\int_{\mathbb{R}^d} \bigg|
\int_0^\i\int_{\mathbb{R}^d} K(x,t)} \\
  & & \qquad \qquad \Big( 2X_n(dx,dt) - \sum_{m=1}^{n-1}
\macp \big( m,n-m \big)
    g_m*_x\eta^{\delta}(x,t)g_{n-m}*_x\eta^{\delta}(x,t) \Big)
    dxdt \bigg| d \hat{\mathcal P}_N = 0, \nonumber
\end{eqnarray}
for any  smooth function $K$ of compact support.
 From this, the continuity of the integrand with respect to $(g,X)$, 
and our assumptions on $g_j$ for $j<n$
we deduce,
$$
\lim_{\delta \to 0}{\int{} \left[ \int_0^\i\int{} K(x,t) \left( 
2X_n(dx,dt) -  \sum_{m=1}^{n-1}
\macp \big( m,n-m\big ) k_m k_{n-m} \right) dxdt\right]^{+}} d
\hat{\mathcal P} = 0,
$$
which implies that
$$
2\int_0^\i\int_{\mathbb{R}^d} K dX_n  \leq \sum_{m=1}^{n-1} \macp 
\big(m,n-m\big)
k_m k_{n-m} \int_0^\i\int_{\mathbb{R}^d} K(x,t)dxdt
$$
for $\hat{\mathcal P}$- almost all realizations of the random measure
$X_n$. So, $X_n =
\g_n dx dt$,
    with $\g_n$ bounded, implying (\ref{mcla}) and the second part of the
Lemma. Applying (\ref{mcla}) in
    (\ref{pbd}), using $\psi_p'(z)\le p\psi_p(z)+p$ and Gronwall's 
inequality, we find that
$$
\sup_{t \in [0,T]}{\int_{\mathbb{R}^d} \psi_p \big(
    \hat f^{\delta}(n,x,t) \big)J(x) dx }\le c,
$$
for a constant $c$ that is independent of $\ell$. From this we can 
readily deduce,
$$
\sup_{t \in [0,T]}\sup_{\d>0}{\int_{\mathbb{R}^d} \psi_p \big(
    g_n*_x\eta^{\delta}(x,t) \big)J(x) dx } <
    \infty,
$$
as required to derive the first part of the Lemma. $\Box$
\end{section}

\begin{section}{Deriving the PDE}\label{derpde}

The aim of this section is to carry out the second and fourth steps of
the outline for deriving the system of PDE (\ref{syspde}) described in
Section \ref{unifb}. We then apply an induction on the number $n$ to 
complete the proof of
Theorem 1.1. In each case, we
choose a test function $\bar K: [0,\infty)
\times \mathbb{R}^d  \to \mathbb{R}$ and consider the
functional
\begin{displaymath}
    Y(q,t) = \epsilon^{d-2} \sum_{i \in I_q}{ K \big( t,x_i,m_i \big) },
\end{displaymath}
where $K(t,x,m)=\bar K(x,t)1\!\!1(m=n)$. Note that
\setcounter{equation}{0}
\begin{eqnarray}\label{yevol}
    Y\big(q(T), T \big) & = & Y\big( {q(0)}, 0 \big) + \int_0^T \mathbb{A}_C
    ( Y )\big ({q}(t),t\big) dt \\
   & & \qquad    + \int_0^T \mathbb{A}_0
    ( Y ) \big({q}(t),t\big) dt +\int_0^T Y_t\big ({q}(t),t\big) dt\, +
    \, M_T, \nonumber
\end{eqnarray}
where
\begin{eqnarray}
    \mathbb{A}_C Y(q,t) & = & \frac{\epsilon^{d-2}}2 \sum_{i,j \in I_q} \micp
    \big( m_i,m_j \big) V_{\epsilon} \big( x_i - x_j \big) \\
    & & \qquad \bigg[ \frac{m_i}{m_i + m_j} K \big( t,x_i, m_i + m_j \big)
    +  \frac{m_j}{m_i + m_j} K \big( t,x_j, m_i + m_j \big)
     \nonumber \\
     & & \qquad \qquad \quad   - K \big( t , x_i, m_i \big) - K \big(
     t , x_j,m_j \big) \bigg], \nonumber
\end{eqnarray}
and where
\begin{eqnarray}
    \mathbb{A}_0 Y(q,t)& = &\epsilon^{d-2} \sum_{i \in I_q} d (m_i) \Delta_x K
    \big( t, x_i ,m_i \big),
    \\
    Y_t(q,t)& = &\epsilon^{d-2} \sum_{i \in I_q}  K_t
    \big( t, x_i ,m_i \big).\nonumber
    \end{eqnarray}
The term $M_T$ is a martingale satisfying
\begin{displaymath}
\mathbb{E}_N \Big[  M_T^2 \Big] =  \mathbb{E}_N\int_{0}^{T}{A_1 \big( q(t),t
\big) dt} +  \mathbb{E}_N\int_{0}^{T}{A_2 \big( q(t),t \big) dt},
\end{displaymath}
where $A_1(q,t)$ and $A_2(q,t)$ are respectively set equal to
\begin{displaymath}
A_1 \big( q ,t\big) = \e^{2(d-2)} \sum_{i \in I_q}  d (m_i)
    \big\vert \nabla_{x_i} K \big\vert^2 \big(t, x_i,m_i \big),
\end{displaymath}
and
\begin{eqnarray}
    A_2 \big( q ,t\big)= &  & \frac{\epsilon^{2(d-2)} }2\sum_{i,j \in I_q}  \micp
    \big( m_i , m_j \big) V_{\epsilon} \big( x_i - x_j \big) \
    \\
     & & \qquad \bigg[ \frac{m_i}{m_i + m_j} K \big( t,x_i, m_i + m_j \big)
    +  \frac{m_j}{m_i + m_j} K \big( t,x_j, m_i + m_j \big)
      \nonumber \\
     & & \qquad \qquad \quad   - K \big( t , x_i, m_i \big) - K \big(
     t , x_j,m_j \big) \bigg]^2. \nonumber
\end{eqnarray}
We can readily show
\begin{eqnarray}
A_1 \big( q ,t\big)  &\leq& C \epsilon^{2(d-2)} \sum_{i \in I_q}  d (m_i)
1\!\!1(m_i=n)
\leq C\e^{d-2},
\\
    A_2 \big( q,t \big)  &\leq&  C\epsilon^{2(d-2)} \sum_{i,j \in I_q}  \micp
    \big( m_i , m_j \big) V_{\epsilon} \big( x_i - x_j \big)
1\!\!1(m_i+m_j\le n)\le C\e^{d-2}, \nonumber
     \end{eqnarray}
    where we have used Lemma 3.1 for the last inequality.
   From these inequalities and Doob's inequality, we deduce that the
martingale tends to zero
uniformly, in the $\epsilon \downarrow 0$ limit.

We rewrite the terms of (\ref{yevol}) in terms of the empirical
measures introduced in (4.23). We have that
\begin{equation}
    Y_q \big( T \big) = \sum_{m \in \mathbb{N}}{\int_{\mathbb{R}^d} \bar K
    \big( T,x \big)  g_n(dx,T)}
\end{equation}
and that
\begin{equation}
    \int_0^T \mathbb{A}_0 Y_q \big( t \big) dt  =\int_0^T
    d (n) \int_{\mathbb{R}^d}  \Delta_x \bar K
    \big( t,x \big) g_n(dx,t).
\end{equation}
The Stosszahlansatz enables us to rewrite the time-averaged action of
the collision operator on $Y \big( q(t), t \big)$ in terms of these
empirical measures. That is, by Proposition \ref{szprop},
\begin{eqnarray}\label{fineps}
   \int_0^T  dt \mathbb{A}_C Y\big(q(t), t \big)
  & = &
         \int_0^T dt \int_{\mathbb{R}^d} dx \bar K (t,x) \nonumber \\
    & & \
     \left[ \frac 12\sum_{m=1}^{n-1}
    \macp \big( m, n-m \big)\big( g_m *_x \eta^{\delta}\big) \big( x,t \big)
\big(g_{n-m} *_x \eta^{\delta}\big) \big( x,t \big)
\right. \\
    & & \ \ \         \,  - \, \left. \sum_{r=1}^{\infty}{ \macp
    \big( n, r  \big)  \big( g_n *_x \eta^{\delta}\big) \big( x,t \big)
\big(g_{r} *_x \eta^{\delta}\big) \big( x,t \big)
     } \right] \, + \, Err \big( \epsilon , \delta \big)\ . \nonumber
\end{eqnarray}
To carry out the second step of the outline,
note that by the assumption, the first $n-1$ marginals in the mass
variable of any weak limit of the measures
${\mathcal P}_N$ have bounded densities $\big\{ f_m:
[0,\infty) \times \mathbb{R}^d  \to
     [0,\infty) : m \in \big\{ 1, \ldots, n-1 \big\}  \big\}$
     with respect to Lesbesgue measure .
   By
     passing to the limit in low $\epsilon$ in (\ref{fineps}), we find
   that any weak limit $\P$ is concentrated on the space of 
$g=(g_m:m\in\N)$ such that,
\begin{eqnarray}
    & &  \int_{\mathbb{R}^d} g_n \big(dx, T \big) \bar K \big(
T,x \big) dx \\
& = &
     \int_{\mathbb{R}^d} g_n \big(  dx,0 \big) \bar K \big( 0,x \big) 
dx +d(n) \int_0^T
     dt \int_{\mathbb{R}^d} g_n \big( dx,t \big) \Delta \bar K \big(
t,x \big) dx
    \nonumber \\
    & & \, + \, \int_0^T  dt \int_{\mathbb{R}^d} dx  \bar K (t,x )\left[
\frac 12\sum_{m=1}^{n-1}
    \macp \big( m, n - m \big) \big( g_m  *_x
\eta^{\delta}\big) \big( x,t \big)  \big(g_{n-m} *_x 
\eta^{\delta}\big) \big( x,t \big)-X_{n,\delta}\right]
  \nonumber \\
    & & \qquad \qquad \quad +
     Err \big( \delta \big), \nonumber
\end{eqnarray}
where $X_{n,\delta} \geq 0$ and the $\P$-expectation of $Err \big( 
\delta \big)$
goes to zero
as $\delta \downarrow 0$. We know by the third step in the outline
that $g_m(dx,t)=f_m(x,t)dx$ with $f_m$ being uniformly bounded for 
$m<n$.  We find that
\begin{eqnarray}
    \int_{\mathbb{R}^d} g_n \big( dx,T \big) \bar K \big( T,x
    \big) dx
    & = &
     \int_{\mathbb{R}^d} g_n \big( dx,0 \big) \bar K \big( 0,x \big) 
dx + d(n)\int_0^T
    ds \int_0^s dt \int_{\mathbb{R}^d} g_n \big( dx,t \big) \Delta \bar K \big(
t,x \big) dx \nonumber \\
    & & \, + \, \int_0^T  dt \bar K \big( t,x \big) \left[\frac 12
\sum_{m=1}^{n-1}
    \macp \big( m, n - m \big) f_m \big( x,t \big)  f_{n-m}\big( x,t \big) -X_n
  \right],
\end{eqnarray}
where $X_n \geq 0$.
 From this we would like to conclude that the measure $g_n=f_n dx$ and that
the sequence $(f_m:m\le n)$
  solves (\ref{prepde}) in weak form.
Since we do not know the absolute continuity of $g_n$ with respect to the
  Lebesgue measure yet, let us observe that if we set $\hat 
f_n=g_n*_x\eta^\d$, then the sequence
  $(f_n:n<m; \hat f_n)$ does satisfy an equation that is similar to 
(4.5). This equation is close
  enough to (4.5) so that we can apply Lemma 4.3(i) to deduce that
  $$
  \sup_\d \sup_t \int\left(g_n*_x \eta^\d\right)^p J\ dx <\i,
  $$
  for every $p>2$. From this we deduce that indeed $g_n =f_n dx$ with 
$f_n$ satisfying the same $L^p$-bounds.

To carry out the fourth step, note that at the end of the third, we know
that  any weak
limit of the measures
${\mathcal P}_{N}$ is concentrated on those measures on $\mathbb{R}^d$
     having densities with respect to Lebesgue measure $\big\{ f_m:
[0,\infty) \times \mathbb{R}^d  \to
     [0,\infty) : m \in \mathbb{N} \big\}$ with $f_m$ uniformly 
bounded for $m<n$ and
     and that $f_n \in L^p
(\mathbb{R}^d \times [0,T])$ for each $p> 2$ and $T>0$. Thus, the
families $ \big\{ \big( f_{n} * \eta^{\delta} \big) \big( f_m *
\eta^{\delta} \big): \mathbb{R}^d \times [0,\infty) \to
[0,\infty), \d>0 \big\}$ are uniformly integrable for $m\le n$. When 
we pass to the limit in
     low $\epsilon$ of (\ref{fineps}) in this case,
   we find
     that,
\begin{eqnarray}
    & &  \int_{\mathbb{R}^d} f_n \big( T,x \big) \bar K \big(
T,x\big) dx \\
& = &
     \int_{\mathbb{R}^d} f_n \big( 0 , x \big) \bar K \big( 0,x \big) 
dx + d(n)\int_0^T
     dt \int_{\mathbb{R}^d} f_n \big( x,t \big) \Delta \bar K \big(
t,x \big) dx
    \nonumber \\
    & & \, + \, \int_0^T  dt \int_{\mathbb{R}^d} dx  \bar K (t,x )\left[\frac 12
\sum_{m=1}^{n-1}
    \macp \big( m, n - m \big)  f_m  \big( x,t \big)  f_{n-m} \big( x,t \big)
\right. \nonumber \\
    & & \qquad \qquad \quad - \, \left. \sum_{m=1}^{n}{ \macp
    \big( n, m  \big)   f_n
    \big( x,t \big)  f_m  \big( x,t \big)  }-Y_n \right]
    \, + \, Err \big( \delta \big), \nonumber
\end{eqnarray}
where $Y_n\ge 0$ and the $\P$-expectation of $Err \big( \delta \big)$
goes to zero
as $\delta \downarrow 0$.

By induction on $n$ we deduce that the functions $\big\{ f_n : n \in 
\mathbb{N} \big\}$ being known to be
uniformly
bounded ,$\P$--almost surely. We then  apply Lebesgue's
theorem to deduce that, for each $n
\in \mathbb{N}$, $f_n * \eta^{\delta} \to f$ almost surely, as $\delta
\downarrow 0$. We find that the weak limit ${\mathcal P}$ of the sequence
${\mathcal P}_{N}$ is supported on those functions $f$ satisfying
the system
(\ref{syspde}) in weak form.
\end{section}

\begin{section}{Potential theory}
\label{sec6}

The purpose of this section is twofold.  Firstly, we show the existence of the
function $u$ that was used in the proof of the \stoss.
Secondly, we establish a
bound on $\macp$ and show that if $\micp \to \i$ then $\macp$
converges to
macroscopic coagulation rate in the case of a hard sphere interaction.
The limit can be expressed in terms of
the Newtonian
capacity of the support of the function $V$.  We start with the
statements of the main
results of this section.  Let $V: {\R}^d \to {\R}$ be a continuous
function of compact support with $V \ge
0$ and $\ds{\int_{\mathbb{R}^d} V(x)dx = 1}$.  We also write $K_0$ 
for the topological closure
of $U_0$ where
$U_0 = \{x: V(x) \ne 0\}$.
\begin{theorem}
\label{th6.1}
For every $\b > 0$, there exists a unique function $u \in
C^2({\R}^d)$ such that
$\ds{\lim_{|x| \to \i} u(x) = 0}$ and
\setcounter{equation}{0}
\begin{equation}
\label{eq6.2}
\D u - \b Vu = \b V.
\end{equation}
Moreover $-1 \le u(x) \le 0$ for every $x \in {\R}^d$.
\end{theorem}

We have applied this theorem in Section~3 with $\ds{\b = \frac
{\micp(n,m)}{d(n)+d(m)}}$.  Also recall that $\macp(n,m)$ is defined by
\begin{equation}
\label{eq6.3}
  \macp (n,m) = \micp(n,m) \int_{\mathbb{R}^d} V(x)(1+u(x))dx.
\end{equation}
More generally, we define $F: (0,\i) \to (0,\i)$ by $\ds{F(\b) = \b 
\int_{\mathbb{R}^d}
V(1+u^{\b})dx}$ with $u = u^{\b}$ as in Theorem \ref{th6.1}.  Recall that the
(Newtonian) capacity of a compact set $K \subseteq {\R}^d$ is given by
\begin{equation}
\label{eq6.4}
\mbox{Cap}(K) = \inf\left\{ \frac {1}{2} \int_{\mathbb{R}^d} 
|\nabla\psi|^2dx: \psi
\in C^1({\R}^d),\
\psi \ge 1 \mbox{ on a neighborhood of $K$}\right\},
\end{equation}
or equivalently,
\begin{equation}
\label{eq6.5}
\mbox{Cap}(K) = \sup\left\{ \mu(K): \supp \mu \subset K,\ c_0 
\int_{\mathbb{R}^d}
|x-y|^{2-d}\mu(dy) \le
1 \mbox{ for every $x$}\right\},
\end{equation}
where $c_0=c_0(d)=(d-2)^{-1}\o_d^{-1}$ with $\o_d$ denoting the surface area of
the unit sphere. We are now ready to state our second main result.

\begin{theorem}
\label{th6.2}
We always have $F(\b) \le \mbox{Cap}(K_0)$.  Moreover $\ds{\lim_{\b \to \i}
F(\b) =
\mbox{Cap}(K_0)}$.
\end{theorem}

\bigskip
\noindent
{\bf Proof of Theorem \ref{th6.1}} \\
\noindent{\bf Step 1.}  Let $J$ be a
bounded continuous
function with $J > 0$, $\ds{\int_{\mathbb{R}^d} J(x)dx = \i}$ and

\noindent{
$\ds{\int_{|x| \ge 1}
J(x)|x|^{4-2d}dx <
\i}$}.  Define
\begin{equation}
\label{eq6.6}
\cH = \left\{ u : \mbox{ $u$ is measurable and }\int_{\mathbb{R}^d} 
u^2(x) J(x) dx <
\i\right\}.
\end{equation}
We then define $\cF: \cH \to \cH$ by
\begin{equation}
\label{eq6.7}
\cF(u)(x) = c_0 \int_{\mathbb{R}^d} |x-y|^{2-d}V(y)u(y)dy.
\end{equation}
Let us verify that $\cF(u) \in \cH$ for $u \in \cH$.  To see this, write
\begin{equation}
\label{eq6.8}
\G(x) = c_0 \int_{\mathbb{R}^d} |x-y|^{2-d}V(y)dy.
\end{equation}

As in (\ref{ussoon}) we can readily show
\begin{equation}
\label{eq6.9}
0 \le \G(x) \le c_1 \min(|x|^{2-d},1)
\end{equation}
for a constant $c_1$.  Also, we may use H\"older's inequality to assert
\begin{eqnarray}\label{eq6.10}
(\cF(u)(x))^2 &= &\left[ \G(x) \int_{\mathbb{R}^d} c_0|x-y|^{2-d}V(y)u(y) \frac
{dy}{\G(x)} \right]^2
\\
&\le &\G(x)\int_{\mathbb{R}^d} c_0|x-y|^{2-d}V(y)u^2(y)dy. \nonumber
\end{eqnarray}
    From this and \eqref{eq6.9} we deduce
\[
\int_{\mathbb{R}^d} (\cF(u)(x))^2J(x)dx \le c_0c_1 
\int_{\mathbb{R}^d} V(y)u^2(y) \left[ \int_{\mathbb{R}^d}
\min(|x|^{2-d},1)|x-y|^{2-d}J(x)dx\right].
\]
If $V(y) \ne 0$ then $|y| \le R_0$ for a suitable $R_0$. If
\[
I(y) = \int_{\mathbb{R}^d} \min(|x|^{2-d},1)|x-y|^{2-d}J(x)dx,
\]
then we have $\ds{\sup_{|y| \le R_0} I(y) < \i}$ because
\begin{eqnarray*}
I(y) &\le &\int_{|x| \le 2R_0} + \int_{|x| > 2R_0}
\min(|x|^{2-d},1)|x-y|^{2-d}J(x)dx \\
&\le &C \int_{|x-y| \le 2R_0} |x-y|^{2-d}dx + C \int_{|x| > 2R_0}
|x|^{2-d}|x|^{2-d}J(x)dx  < \i
\end{eqnarray*}
by our assumption on $J$.  As a result,
\[
\int_{\mathbb{R}^d}(\cF(u)(x))^2 J(x)dx \le C \int_{\mathbb{R}^d} 
V(y)u^2(y)dy < \i
\]
because $V$ is of compact support.  This shows that $\cF(u) \in \cH$
whenever $u \in \cH$. \\
\noindent{\bf Step 2.}  Observe that $\cH$ is a Hilbert space with respect to
the inner product
\[
\<u,v\> = \int_{\mathbb{R}^d} u(x)v(x) J(x)dx.
\]
Note that if $u$ solves \eqref{eq6.2}, then
\[
u +\b(-\D)^{-1}\left((1+u)V\right)=0.
\]
As a result,
\begin{equation}
\label{eq6.11}
(id + \b\cF)(u) = g
\end{equation}
where $g(x) = -\b\G(x)$ with $\G$ as in \eqref{eq6.8} and $id$ means
the identity transformation.
The formulation
\eqref{eq6.11} allows us to use some standard arguments to establish
the existence of
a solution of \eqref{eq6.2}.  Note that our assumption on $J$
implies that $\G \in
\cH$ because of \eqref{eq6.9}.  Hence as a strategy for the
existence, we show that
the range of the operator $id + \b\cF$ is the space $\cH$.  By
the
Fredholm Alternative Theorem, it suffices to show that the operator
$\cF$ is compact
and that the operator $id + \b\cF$ has a trivial kernel. \\
\noindent{\bf Step 3.}  In this step we verify the injectivity of the operator
$id + \b\cF$.  Let us assume that for some $u \in \cH$, $\b\cF u = -u$.
    From this we deduce
$\D u = \b uV$ weakly by standard arguments.  Using Sobolev's
inequalities and a
bootstrap we can readily show that $u \in C^2$ because $V$ is continuous.  As a
result, $\D u = \b uV$ in the classical sense.  We also have
\begin{equation}
\label{eq6.12}
u(x) =- c_0\b \int_{\mathbb{R}^d} u(y)|x-y|^{2-d}V(y)dy.
\end{equation}
We now assert
\begin{equation}
\label{eq6.13}
\lim_{|x| \to \i} u(x)|x|^{d-2} = -c_0\b \int_{\mathbb{R}^d} V(y)u(y)dy,
\end{equation}
\begin{equation}
\label{eq6.14}
\lim_{|x| \to \i} |\nabla u(x)||x|^{d-1} \le c_0(d-2)\b 
\int_{\mathbb{R}^d} V(y)|u(y)|dy
    .
\end{equation}
The proof of \eqref{eq6.13} follows from
\[
\lim_{|x| \to \i} |x|^{d-2}|x-y|^{2-d} = 1
\]
for every $y$, and the elementary inequality
\[
|x|^{d-2}|x-y|^{2-d} \le 2^{d-2}
\]
for $|x| \ge 2R_0$ (because $|y| \le R_0$ and $\ds{|x-y| \ge |x|/2}$
whenever $|x| \ge 2R_0$).
    The proof of \eqref{eq6.14} is
similar.  We now
choose $R > R_0$ and use $\D u =\b  Vu$ to write
\[
    \int_{|x| \le R} u\D u dx = \b\int_{|x| \le R} Vu^2dx.
\]
After an integration by parts we obtain
\[
- \int_{|x| \le R} |\nabla u|^2dx +  \int_{|x| = R} u \nabla u
\cdot n d\s =\b
\int_{|x| \le R} Vu^2dx,
\]
where $\ds{n = \frac {x}{|x|}}$ is the normal vector and $d\s$ is the
surface measure on $|x| = R$.  From \eqref{eq6.13}--\eqref{eq6.14} we
can readily
deduce
\[
\int_{|x| = R} u \nabla u \cdot n d\s = O(R^{2-d}).
\]
As a result,
\[
- \int_{\mathbb{R}^d} |\nabla u|^2dx =\b \int_{\mathbb{R}^d} Vu^2dx.
\]
    From this we deduce $\ds{\int_{\mathbb{R}^d} |\nabla u|^2dx = 
\int_{\mathbb{R}^d} Vu^2dx = 0}$.  This in
turn implies $u \equiv 0$.\\
\noindent{\bf Step 4.}  In this step we verify the compactness of the operator
$\cF$.  Set
$\phi(x) = c_0|x|^{2-d}$ so that $\cF(u) = (Vu) * \phi$.  We then
have $\nabla \cF(u)
= (Vu) * \nabla \phi$.  Let $B_R$ denote the ball $\{x: |x| \le R\}$.
We then use the
Young's inequality to write
\begin{equation}
\label{eq6.15}
\|\nabla\cF(u)\|_{L^2(B_R)} \le \|Vu\|_{L^2(B_R)} \|\nabla\phi\|_{L^1(B_R)} \le
CR\|Vu\|_{L^2(B_R)},
\end{equation}
\begin{equation}
\label{eq6.16}
\|\cF(u)\|_{L^2(B_R)} \le \|Vu\|_{L^2(B_R)} \|\phi\|_{L^1(B_R)} \le CR^2
\|Vu\|_{L^2(B_R)},
\end{equation}
for some constant $C$.  Now if $\{u_n\}$ is bounded in $\cH$, then by
positivity of
$J$ we deduce that $\ds{\sup_n\|u_n\|_{L^2(B_R)} < \i}$.  This and
\eqref{eq6.15}-\eqref{eq6.16} imply
\begin{equation}
\label{eq6.17}
\sup_n (\|\cF(u_n)\|_{L^2(B_R)} + \|\nabla\cF(u_n)\|_{L^2(B_R)}) < \i.
\end{equation}
    From this and Rellich's theorem we learn that
$\{\cF(u_n)\}$ has a
convergent subsequence in $L_{{loc}}^2({\R}^d)$.  As in Step~1 we have
\[
\int_{|x| \ge l} (\cF(u_n)(x))^2 J(x)dx \le C \left[ \int_{|x| \ge l}
|x|^{4-2d}J(x)dx\right] \cdot \int_{\mathbb{R}^d} V(y)u_n^2(y)dy,
\]
for $l\ge 2 R_0$. As a result,
\[
\lim_{l \to \i} \sup_n \int_{|x| \ge l} (\cF(u_n)(x))^2J(x)dx = 0.
\]
This and the precompactness of $\{\cF(u_n)\}$ in
$L_{{loc}}^2({\R}^d)$ imply the
precompactness of $\{\cF(u_n)\}$ in $\cH$.\\
\noindent{\bf Step 5.}  So far we have shown the existence of a unique
solution $u \in \cH$ of
$u + \b \cF(u) = g$.  As in Step~3 we can readily show that $u$ is a
weak solution of
\eqref{eq6.2}.  Again we can use the Sobolev's inequality and a
bootstrap to show that
in fact $u \in C^2({\R}^d)$ and satisfies \eqref{eq6.2} in the 
classical sense.\\
\noindent{\bf Step 6.}  It remains to show that for the unique solution $u$ of
\eqref{eq6.2}, we have
$-1 \le u(x) \le 0$ for every $x$.  First take a smooth function
$\varphi_{\d}: {\R}
\to (-\i,0]$ such that $\varphi'_{\d} \ge 0$ and
\[
\varphi_{\d}(r) = \begin{cases}
0 &r > -1, \\
1 + r &r < -1-\d.
\end{cases}
\]
We then have
\begin{equation}
\label{eq6.18}
-\int_{\mathbb{R}^d} \varphi'_{\d}(u)|\nabla u|^2dx = 
\int_{\mathbb{R}^d} \varphi_{\d}(u)\D u dx = \b \int_{\mathbb{R}^d}
V(1+u)\varphi_{\d}(u)dx
\end{equation}
by an integration by parts.  The formula \eqref{eq6.18} is
established with the aid of
\eqref{eq6.13} as in Step~3.  Since the left-hand side of
\eqref{eq6.18} is negative and $(1+u)\varphi_{\d}(u) \ge 0$ we deduce
\[
\int_{\mathbb{R}^d} \varphi'_{\d}(u)|\nabla u|^2dx = 
\int_{\mathbb{R}^d} V(1+u)\varphi_{\d}(u)dx = 0.
\]
We now send $\d \to 0$ to deduce
\[
0 = \int_{\mathbb{R}^d} |\nabla u|^21\!\!1(1+u \le 0)dx = 
\int_{\mathbb{R}^d} V(1+u)^21\!\!1(1+u \le 0)dx.
\]
As a result, on the set $A = \{x: 1 + u(x) < 0\}$ we have $\nabla u =
0$.  Hence
$u$ is constant on
each component $B$ of $A$.  But this constant can only be $-1$
because on the boundary of $A$ we
have $1+u=0$.  This is impossible unless $A$ is empty and we deduce
that $u \ge -1$ everywhere.
Finally, from
\[
u(x)=c_0\int_{\mathbb{R}^d} V(y)(1+u(y))|x-y|^{2-d}dy,
\]
and $1+u\ge 0$ we deduce that $u\le 0$.
\qed

\bigskip
\noindent
{\bf Proof of Theorem \ref{th6.2}}  Let us write $u^{\b}$ for the
unique solution of
\eqref{eq6.2}.  Note that by Theorem~\ref{th6.1},
\begin{equation}
\label{eq6.19}
-1 \le u^{\b}(x) \le 0
\end{equation}
for every $x$.  Since $\D u^{\b} = \b(1+u^{\b})V$, we also have
\begin{equation}
\label{eq6.20}
u^{\b}(x) =  -c_0 \int_{\mathbb{R}^d}
|x-y|^{2-d}\mu^{\b}(dy).
\end{equation}
where $\mu^{\b}(dy) = \b(1+u^{\b}(y))V(y)dy$.
    From this, \eqref{eq6.19} and \eqref{eq6.5} we learn
\begin{equation}
\label{eq6.21}
F(\b) = \int_{\mathbb{R}^d} \b(1+u^{\b}(y))V(y)dy = \mu^{\b}({\R}^d) 
\le \mbox{Cap}(K_0).
\end{equation}

It remains to show $\ds{\lim_{\b \to \i} F(\b) = \mbox{Cap}(K_0)}$.  By
\eqref{eq6.21}
the family $\{\mu^{\b}\}$ is precompact with respect to the weak
topology.  Let $\mu$
be a limit point and put
\begin{equation}
\label{eq6.22}
u(x) = -c_0 \int_{\mathbb{R}^d} |x-y|^{2-d}\mu(dy).
\end{equation}
Assume $\ds{\lim_{\b' \to \i} \mu^{\b'} = \mu}$ where $\{\mu^{\b'}\}$
is a subsequence
of $\{\mu^{\b}\}$.  It is not hard to see
\[
\limsup_{\b' \to \i} u^{\b'}(x) \le
u(x)
\]
by replacing the integrand in \eqref{eq6.22} and \eqref{eq6.20} with a bounded
cutoff function and passing to the limit.  Since $u^{\b} \ge -1$, we
deduce that $u\ge -1$.  Also \eqref{eq6.22} implies that $u \le 0$.  Hence $-1
\le u \le 0$.

In fact it is well known that $\ds{\limsup_{\b \to\i} u^{\b'}(x) =
u(x)}$ except on a
set of zero capacity.  (See for example page~161 or the proof of
Theorem~$5.7.1$ in
\cite{cpt}.)  Since $\D u^{\b} = \mu^{\b}$, $\D u = \mu$ and $\mu^{\b}
\Rightarrow \mu$, we
can readily deduce that $u^{\b} \to u$ weakly.  As a result,
\[
\lim_{\b \to \i} \int_{\mathbb{R}^d} (1+u^{\b})Vdx = 
\int_{\mathbb{R}^d} (1+u)Vdx.
\]
Since $\ds{\b \int_{\mathbb{R}^d} (1+u^{\b})Vdx}$ is bounded in $\b$, we deduce
\[
\int_{\mathbb{R}^d} (1+u)Vdx = 0.
\]
As a result, $u \equiv -1$ almost everywhere in $U_0 = \{x: V(x) \ne
0\}$.  Since
$\ds{u(x) = \limsup_{y \to x} u(y)}$ except for a set of zero
capacity, we learn $u
\equiv -1$ in $U_0$ except for a set of zero capacity.  (See, for example,
Theorem~$3.7.5$ of \cite{cpt}.)  Let $\nu$ denote the equilibrium measure of
$K_0$. That is, $\supp(\nu) \subset K_0$, $\nu(K_0)=Cap(K_0)$, and
that the function $\ds{v(x)
= c_0 \int_{\mathbb{R}^d} |x-y|^{2-d}\nu(dy)}$ satisfies $v \equiv 1$ 
in $K_0$ except
for a set of zero
capacity.  Also, since $u$ is finite everywhere, we learn that
$\mu(A)=0$ whenever $Cap(A)=0$.
The same is also true for $\nu$. (See for example Theorem 5.1.9 of
[1].) As a result,
\begin{equation}
\label{eq6.23}
\mu(K_0) = \int_{\mathbb{R}^d} vd\mu = c_0 \int_{\mathbb{R}^{2d}} 
|x-y|^{2-d}\nu(dy)\mu(dx) = -\int_{\mathbb{R}^d}
u(y)\nu(dy).
\end{equation}
But $\nu(A) = 0$ if $\mbox{Cap}(A) = 0$ and $-u \equiv 1$ on $K_0$
except for a set of
zero capacity.  This and \eqref{eq6.23} imply
\[
\mu(K_0) = \nu(K_0) = \mbox{Cap}(K_0).
\]
\qed
\end{section}
\bibliographystyle{plain}
\bibliography{biblion}
\end{document}